\documentclass[a4paper,oneside]{report}
\usepackage[utf8]{inputenc}
\usepackage[T1]{fontenc}
\usepackage[english]{babel}
\usepackage{amsmath} 
\usepackage{amssymb} 
\usepackage{cancel,amsfonts,graphics,mathdots}
\usepackage{braket} 
\usepackage{hyperref}
\usepackage{xfrac} 

\newcommand{\numberset}{\mathbb}
\newcommand{\PP}{\numberset{P}}
\newcommand{\Gr}{\numberset{G}}
\newcommand{\C}{\numberset{C}}
\newcommand{\N}{\numberset{N}}

\DeclareMathOperator{\Hess}{Hess}
\DeclareMathOperator{\hess}{hess} 
\DeclareMathOperator{\Ann}{Ann}
\DeclareMathOperator{\rank}{rank}
\DeclareMathOperator{\sol}{Sol}
\DeclareMathOperator{\tr}{tr}
\DeclareMathOperator{\Sing}{Sing}
\DeclareMathOperator{\Soc}{Soc}
\DeclareMathOperator{\Cat}{Cat}
\DeclareMathOperator{\Jac}{Jac}
\DeclareMathOperator{\Imm}{Im}
\DeclareMathOperator{\codim}{codim}

\usepackage[style=alphabetic,backend=biber,maxbibnames=10]{biblatex}
\usepackage{amsthm}
\theoremstyle{plain}
\newtheorem{theorem}{Theorem}[chapter]
\newtheorem{lemma}[theorem]{Lemma}
\newtheorem*{lemma*}{Lemma}
\newtheorem*{theorem*}{Theorem}
\newtheorem{proposition}[theorem]{Proposition}
\newtheorem{corollary}[theorem]{Corollary}

\theoremstyle{definition}
\newtheorem*{definition*}{Definition}
\newtheorem{definition}[theorem]{Definition}

\theoremstyle{remark}
\newtheorem*{example*}{Example}
\newtheorem{example}[theorem]{Example}
\newtheorem*{remark*}{Remark}
\newtheorem{remark}[theorem]{Remark}

\newtheorem*{notation}{Notation}
\newtheorem*{important}{Important hypothesis}

\newcommand{\pause}{\vspace{5pt}}
\newcommand{\I}{\mathcal{I}}
\newcommand{\G}{\mathcal{G}}
\newcommand{\D}{\mathcal{D}}

\usepackage{chngcntr}
\counterwithout{footnote}{chapter}

\usepackage{graphicx,color,xcolor,csquotes}
\definecolor{MyDarkGreen}{cmyk}{0.7,0,1,0}

\newcommand\Tstrut{\rule{0pt}{3ex}}
\newcommand\Bstrut{\rule[-1.6ex]{0pt}{0pt}}

\newcommand{\xrightarrowdbl}[2][]{%
    \leftarrow\mathrel{\mkern-14mu}\xrightarrow[#1]{#2}
}

\usepackage[font=small,format=hang,labelfont=bf]{caption}
\usepackage{float}

\usepackage{stackengine}

\author{Luca Fiorindo}
\usepackage{verbatim}

\begin{document}
    \thispagestyle{empty}

    \hfill
    \begin{center}
        \Large
        {\bf\LARGE\textsc{Polynomials with vanishing Hessian and Lefschetz properties}}
    \end{center}
    \bigskip\bigskip
    
    \begin{center}
        {\Large\textit{Luca Fiorindo}}\\
        email:luca.fiorindo@dima.unige.it
    \end{center}
    \smallskip
    
    \begin{center}
        Dipartimento di Matematica e Geoscienze\\
        Università degli Studi di Trieste\\
        Via Valerio 12/1, 34127, Trieste, Italy\\
    \end{center}
    \hfill

    \begin{abstract}
        This dissertation is the author's Master thesis in Mathematics. This thesis has been written under the supervision of Prof. Emilia Mezzetti, and it was defended in Trieste on July $22^{nd}$ 2022. The thesis has been corrected of some misprints, but the original document can be found in \cite{tesimia}. The aim is to study \emph{Perazzo hypersurfaces} $X=V(F)\subseteq\PP(K^5)$, defined by $F(x_0,x_1,x_2,u,v)=p_0(u,v)x_0+p_1(u,v)x_1+p_2(u,v)x_2+g(u,v)$, where $p_0,p_1,p_2$ are algebraically dependent, but linearly independent forms of degree $d-1$ in $u,v$, and $g$ is a form in $u,v$ of degree $d$. These hypersurfaces are the "building blocks" for all possible hypersuface in $\PP^4$ with vanishing Hessian. A minimal and a maximal Hilbert vector is found for the associated Artinian Gorenstein $K$-algebras $A_F$: in the minimal case they satisfy the Weak Lefschetz property, but in the maximal case they don't. Furthermore, we classify all Perazzo $3$-folds with minimal $h$-vector. We also summarise basic knowledge and already known results about hypersurfaces with vanishing Hessian and their geometry in low dimension, and also about Artinian Gorenstein $K$-algebras.
    \end{abstract}
	
	\tableofcontents
	\thispagestyle{empty}
	\clearpage
	\cleardoublepage\phantomsection\setcounter{page}{1}
	\pagestyle{plain}
	
	\addcontentsline{toc}{chapter}{Sunto esteso}
	\chapter*{Sunto esteso}
	    La matrice Hessiana di una ipersuperficie e il suo determinante sono dei concetti matematici ben noti anche agli studenti dei primi anni di studio universitario. È noto che i coni sono una classe elementare di ipersuperfici con Hessiano nullo, ma non tutte queste particolari ipersuperfici sono necessariamente dei coni. Perazzo studiò in \cite{perazzo} alcuni tipi di ipersuperfice con Hessiano nullo che sono state poi nominate \emph{ipersuperfici di Perazzo} (Definizione \ref{perazzo}).
	    
	    \pause
	    Lo scopo di questa tesi è studiare le ipersuperfici di Perazzo in $\PP^4$ ($3$-varietà di Perazzo) della forma $$f=p_0(u,v)x_0+p_1(u,v)x_1+p_2(u,v)x_2+g(u,v)\in K[x_0,x_1,x_2,u,v]_d,$$con $K$ campo algebricamente chiuso di caratteristica zero, e in particolare le algebre Artiniane di Gorenstein ad esse associate (Definizione \ref{definizionedigorenstein}). Definito l'anello degli operatori differenziali $S=K\left[\frac{\partial}{\partial x_0},\frac{\partial}{\partial x_1},\frac{\partial}{\partial x_2},\frac{\partial}{\partial u},\frac{\partial}{\partial v}\right]$, si definisce l'algebra graduata Artiniana di Gorenstein $A$ associata a $f$ il quoziente $$A=\frac{S}{\Ann_S(f)}=\frac{S}{\left\{g\in S\,\Big\vert\, g\left(\frac{\partial}{\partial x_0},\dots,\frac{\partial}{\partial x_n}\right)f(\mathbf{x})=0\right\}}.$$
	    
	    Con questo obiettivo in mente, studiamo il vettore di Hilbert associato a $A$, Sezione \ref{constructionofthehilbertvector}; in particolare, nella Proposizione \ref{hilbert function}, diamo una formula esplicita che coinvolge le matrici cataletticanti (Definizione \ref{catalecticant}) dei polinomi $p_0,p_1,p_2$. In questo modo siamo in grado di dare una maggiorazione e una minorazione dell'$h$-vettore, rispettivamente nei Teoremi \ref{upper} e \ref{lower}. Un esempio di $3$-varietà di Perazzo con vettore di Hilbert massimo è data nell'Esempio \ref{esempiodellavita} in ogni grado; mentre, una classificazione delle $3$-varietà di Perazzo con vettore di Hilbert minimale è data nel Teorema \ref{lowcharacterisationreprise}.
	    
	    \pause
	    Successivamente andremo a studiare la proprietà debole di Lefschetz usando il teorema di Watanabe (Teorema \ref{watanabe}). Questo teorema connette le matrici Hessiane superiori (Definizione \ref{higherhessiansdefinizione}) alle mappe di moltiplicazione per un dato elemento $L\in A_1$. Proveremo che se l'$h$-vettore è minimo allora questa proprietà è sempre verificata da $A$; mentre se l'$h$-vettore è massimo, allora $A$ fallisce sempre questa proprietà (Proposizioni \ref{main1} e \ref{main2}). Infine, studieremo alcune proprietà geometriche delle ipersuperfici con Hessiano nullo usando come riferimento \cite{franchetta}, \cite{russo},\cite{gondimrusso},\cite{gondimrussostagliano},\cite{cilibertorussosimis}.
	    
	    \pause
	    Vediamo ora nel dettaglio la composizione di questa tesi. Nel \textbf{Capitolo 1} illustreremo l'articolo di M. Watanabe e M. de Bondt \cite{watanabe2019theory}. Qui, gli autori studiano e formalizzano i risultati classici di Gordan e Noether contenuti in \cite{gordannoether}. In particolare, nei Teoremi \ref{binary}, \ref{ternary}, \ref{quaternary}, e \ref{quinary}, viene data la classificazione delle ipersuperfici con Hessiano nullo in $\PP^n$ con $n\le 4$.
	    
	    \pause
	    Nel \textbf{Capitolo 2}, andremo ad introdurre la teoria delle algebre Artiniane di Gorenstein, del vettore di Hilbert e delle proprietà di Lefschetz. Verranno introdotti risultati che ci permetteranno di studiare la proprietà debole di Lefschetz e, in particolare, introdurremo gli Hessiani superiori e il teorema di Watanabe (Teorema \ref{watanabe}).
	    
	    \pause
	    Il \textbf{Capitolo 3} rappresenta la parte centrale di questa tesi. Nella Sezione \ref{catal}, introdurremo le matrici cataletticanti di una certa forma $h\in K[u,v]_t$ e daremo il loro legame con la posizione di $[h]$ in $\PP(K[u,v]_t)$ rispetto alla curva razionale normale $C_t$ e alle sue varietà secanti $\sigma_r(C_t)$. Studieremo poi le algebre Artiniane di Gorenstein associate alle $3$-varietà di Perazzo. Nella Sezione \ref{constructionofthehilbertvector}, studieremo il loro vettore di Hilbert stabilendo un limite superiore ed inferiore (Proposizioni \ref{upper} e \ref{lower}). Daremo la classificazione delle $3$-varietà di Perazzo con vettore di Hilbert minimale (Teorema \ref{lowcharacterisationreprise}) e esempi di $3$-varietà di Perazzo con vettore di Hilbert massimo (Esempio \ref{esempiodellavita}). Andremo poi a studiare la proprietà debole di Lefschetz dimostrando che: se il vettore di Hilbert è massimale, allora l'algebra di Gorenstein fallisce la proprietà (Teorema \ref{main1}); mentre se è minimale, allora sempre la soddisfa (Teorema \ref{main2}). Nella parte finale del capitolo cercheremo di generalizzare quanto dimostrato alle ipersuperfici di $\PP^4$ con  Hessiano nullo o meno.
	    
	    \pause
	    Infine, nel \textbf{Capitolo 4}, andremo a definire alcuni oggetti geometrici quali la varietà duale e la mappa polare di un'ipersuperficie per poter dare una descrizione geometrica delle ipersuperfici con Hessiano nullo. Nel risultato principale, Teorema \ref{cfrs}, viene riassunto il lavoro presente in \cite{franchetta} e \cite{cilibertorussosimis} che descrive le ipersuperfici con Hessiano nullo in $\PP^4$. Alla fine di questo capitolo faremo alcuni calcoli espliciti andando a esemplificare il Teorema \ref{cfrs} nel caso delle $3$-varietà di Perazzo con vettore di Hilbert minimale.

	\cleardoublepage\phantomsection
	\addcontentsline{toc}{chapter}{Introduction}
	\chapter*{Introduction}
	    The Hessian matrix of a hypersurface and its determinant, the Hessian, are basic mathematical concepts that are studied by first years students. It is an elementary fact that a cone has always vanishing Hessian, i.e. the determinant of its Hessian matrix is identically zero. The study of the opposite implication is not a trivial fact, and this topic has been studied  over the last centuries until the present day. The first mathematicians who had studied hypersurfaces with vanishing Hessian were O. Hesse, M. Pasch, P. Gordan, and M. Noether. Since 1852 O. Hesse claimed twice in \cite{hesse1} and \cite{hesse2} that the only possible hypersurfaces with zero Hessian are the cones.
	    
	    \pause
	    The first counterexample to that claim was stated and studied by P. Gordan, and M. Noether in 1876 \cite{gordannoether}. They noticed that being a form with vanishing Hessian is related to the algebraic dependence of the partial derivatives of that form. On the contrary, a form being a cone is related to the linear dependence of the partial derivatives of that form. From this point of view, it seems plausible that Hesse's claim may not be true even in forms involving a low number of variables. However, they proved that given a form $f\in K[x_0,\dots,x_n]$ with vanishing Hessian, then, if $n\le 3$, $X=V(f)\subset\PP^n$ is indeed a cone, proving Hesse's claim in these cases. Furthermore, for $n=4$ they gave a complete description of such forms.
	    
	    \pause
	    The classical paper \cite{gordannoether} has been studied by many mathematicians trying to analyse, formalise and reinterpreting their theory. Among them, the works of J. Watanabe in \cite{watanabeerror}, and \cite{watanabe2019theory}, the last one together with M. de Bondt, are some of the most noticeable ones. Inspired by the notes of H. Yamada, Watanabe's professor, he introduced the so called "self-vanishing systems" which are tuples $(h_0,\dots,h_n)$ of forms in $n+1$ variables of the same degree such that every form in it satisfies the partial differential equation
	    $$h_0\frac{\partial f}{\partial x_0}+\dots+h_n\frac{\partial f}{\partial x_n}=0.$$
	    In this way he formalised Gordan-Noether's results giving them a more algebraic interpretation. All his work was done under the hypothesis of $K$ to be an algebraically closed field of characteristic zero.
	    
	    \pause
	    Watanabe was interested in forms with vanishing Hessian mainly because of their implications in the theory of Artinian Gorenstein algebras and, in particular, in their relations with the Strong Lefschetz property. We recall the definition of such mathematical objects: given a form $f\in K[x_0,\dots,x_n]$, we define the annihilator of $f$ in the ring of differential operators $S=K\left[\frac{\partial}{\partial x_0},\dots,\frac{\partial}{\partial x_n}\right]$ as 
	    
	    $$\Ann_S(f):=\left\{g\in S\,\Big\vert\, g\left(\frac{\partial}{\partial x_0},\dots,\frac{\partial}{\partial x_n}\right)f(\mathbf{x})=0\right\}.$$
	    
	    We then define the Artinian  \emph{Gorenstein algebra} associated to $f$ to be the quotient $A=S/\Ann_S(f)$. The quotient $A=\oplus_{i=0}^d A_i$ results to be a graded Artinian $K$-algebra with standard grading. $A$ is said to have the \emph{Strong Lefschetz property} (respectively, the \emph{Weak Lefschetz property}), if there exists an element $L\in A_1$ such that the multiplication map $\times L^k:A_s\to A_{s+k}$ has full rank for every integer $k$ and $s$ (respectively, the map $\times L:A_s\to A_{s+1}$ has full rank for every integer $k$).
	    
	    \pause
	    Watanabe in \cite{watanabequeen} and, together with T. Maeno, in \cite{maenowatanabe} proved, under the hypothesis of $K$ to be an algebraically closed field of characteristic zero, that the matrix associated to the map
	    $$\times L^{d-2}:A_1\longrightarrow A_{d-1}$$
	    is the Hessian matrix of $f$ up to a non-zero constant. Furthermore, they introduced the \emph{higher Hessians} to describe all the possible maps of the type
	    $$\times L^{d-2k}:A_k\longrightarrow A_{d-k}.$$
	    The rank's maximality  of all these maps turns out to be a necessary and sufficient condition for $A$ having the Strong Lefschetz property.
	    
	    \pause
	    Due to this important result, hypersurfaces with vanishing Hessian give rise to Artinian Gorenstein algebras which lack the SLP. Thus, it is a natural question to study the Weak Lefschetz property on such class of forms. This problem was first studied by R. Gondim in \cite{gondim} where examples of algebras both with and without the WLP were given. In particular, he proved that, for $d=3$, every such Gorenstein algebra fails the WLP; intead, for $d=4$, it always satisfies the WLP.
	    
	    \pause
	    Another important object studied by Gondim in \cite{gondim} is a special class of hypersurfaces with vanishing Hessian: the \emph{Perazzo hypersurfaces}. Following his notations, a Perazzo hypersurface of degree $d$ in $\PP^N$ is defined by a form $f\in K[x_0,\cdots ,x_n,u_1\cdots ,u_m]$ of the following type:
        $$ f=x_0p_0+x_1p_1+\cdots +x_np_n+g$$
        where
        \begin{itemize}
            \item $n,m$ are positive integer numbers with $n+m=N$ and $n,m\ge 2$;
            \item $p_i\in K[u_1,\cdots ,u_m]_{d-1}$ are algebraically dependent, but linearly independent;
            \item $g\in K[u_1,\cdots ,u_m]_{d}$.
        \end{itemize}
        
	    The usefulness of this class of forms is justified by Gordan-Noether's theorem in $\PP^4$. Indeed, every form $f\in K[x_0,x_1,x_2,u,v]_d$ with vanishing Hessian, not a cone, can be written, up to a change of variables, as an element of the extension $K[u,v][\Delta]$, where $\Delta=x_0p_0+x_1p_1+x_2p_2$ is a Perazzo hypersurface. In the case $3\le d\le 5$ then every form with vanishing Hessian is a Perazzo form; but this is not true for degrees greater than $6$.
	    
	    \pause
	    Another important reinterpretation of Gordan-Noether's work was given by A. Franchetta in \cite{franchetta}. Here, inspired by the work of U. Perazzo \cite{perazzo} on cubics with vanishing Hessians, he studied the geometry of hypersurfaces with vanishing Hessian in $\PP^4$ giving them a precise description: every such hypersurface is the union of a one-dimensional family of planes satisfying certain properties. Later on, his work has been studied by mathematicians like R. Permutti, R. Gondim, F. Russo, G. Staglianò, C. Ciliberto, and A. Simis trying to further generalise this general topic in many ways (see for example \cite{permutti1}, \cite{permutti2}, \cite{permutti3}, \cite{gondimrusso}, and \cite{gondimrussostagliano}). In particular the work of Ciliberto, Russo, and Simis \cite{cilibertorussosimis} captures our attention: in this paper they studied also the geometry of the dual variety using the polar map and the Gauss map. They proved that, given a hypersurface $X\subset\PP^4$ with vanishing Hessian not a cone, then the dual variety $X^*$ is a scroll variety satisfying certain suitable properties.
	    
	    \pause
	    The aim of this thesis is to study the Artinian  Gorenstein algebras associated to Perazzo hypersurfaces in $\PP^4$, i.e. the Perazzo $3$-folds; so that, we consider forms of the type $f=p_0x_0+p_1x_1+p_2x_2+g\in K[x_0,x_1,x_2,u,v]_d$, with $K$ an algebraically closed field of characteristic zero, and we want to study the intrinsic properties of the Gorenstein algebra $A=S/\Ann_S(f)$. With this aim in mind, we want first to compute the Hilbert vector associated to such algebra. Using the \emph{catalecticant matrices} of every form $p_0,p_1$ and $p_2$, in Proposition \ref{hilbert function}, a closed formula for the Hilbert vector of $A$ is given. In this way we are able to give, in Theorems \ref{upper} and \ref{lower}, respectively an upper and a lower bound on the possible Hilbert vectors, coinciding only if $d=3,4$. In Example \ref{esempiodellavita}, an example of a Perazzo $3$-fold with maximal $h$-vector is given in every degree. We see that the most expected $h$-vector is the maximal one, and the less expected is the minimal one. If we define the plane $\pi$ generated by $p_0,p_1,p_2$ viewed as elements of $\PP(K[u,v]_{d-1})$, then $A$ has minimal $h$-vector if and only if $\pi$ meets the rational normal curve $C_{d-1}$ at exactly three points counted with multiplicity, i.e.
	    \begin{itemize}
	        \item $\pi$ intersects $C_{d-1}$ in three distinct points;
	        \item $\pi$ is tangent at $C_{d-1}$ at a certain point and intersects $C_{d-1}$ at another point;
	        \item $\pi$ is a osculating plane for $C_{d-1}$.
	    \end{itemize}
	    
	    In this way we are able to use the theory of Waring rank for binary symmetric tensors and obtain in Theorem \ref{lowcharacterisationreprise} a classification of Perazzo $3$-folds with minimal $h$-vector up to projectivity. 
	    
	    \pause
	    Secondly, we study the Weak Lefschetz property in the case of Perazzo $3$-folds with minimal or maximal Hilbert vector. We prove the following: if the $h$-vector is minimal then $A$ has the WLP, on the contrary if the $h$-vector is maximal then $A$ always fails the WLP. Then, we try to generalise these properties to general hypersurfaces with zero Hessian in low degrees.
	    
	    \pause
	    Finally, we list Franchetta's theorem and Ciliberto, Russo, and Staglianò's theorems, and we give some computations about the class of Perazzo $3$-folds with minimal $h$-vector and their dual varieties, verifying the cited results.
	    
	    \pause
	    Next we outline the structure of this thesis. In \textbf{Chapter 1}, we discuss Watanabe and de Bondt paper \cite{watanabe2019theory}. In particular, in Section \ref{selfvanishingsystems} we define analytically the self-vanishing systems giving some properties like the relation between the partial derivatives of a given form and the self-vanishing systems (Theorem \ref{definizionedih}), and the Gordan-Noether identity in Theorem \ref{gnidentity}. In Section \ref{classificationofformswithzerohessian} we discuss the classification of hypersurfaces with vanishing hessian: Theorem \ref{proofgordannoether} proves that, up to a Cremona transformation, every form with zero Hessian is a cone. Then in Theorems \ref{binary}, \ref{ternary}, \ref{quaternary}, and \ref{quinary} we present Gordan-Noether classification up to linear transformations in the cases $n=1,2,3,4$. For the seek of completeness we give also the classification of cubic hypersurfaces in $\PP^5$ and $\PP^6$. The chapter ends with an historical note about vanishing Hessian's problem.
 	    
	    \pause
	    In \textbf{Chapter 2}, we introduce the theory of Artinian Gorenstein algebras and the Lefschetz properties. For Artinian Gorenstein algebras the principal result is the Double Annihilator theorem of Macaulay, stated in Theorem \ref{gorensteinannullatore} which states equivalence between these algebras and hypersurfaces. We introduce the Hilbert vector in Definition \ref{hilbertfunction}, and Lefschetz properties in Definitions \ref{wlpdefinition}, \ref{slp}, and \ref{slpins}. A particular attention is made to Proposition \ref{checkweak}, where various tools to check whenever the WLP holds are stated. In Section \ref{higherhessians}, we introduce the higher Hessians which give a necessary and sufficient condition for an Artinian Gorenstein algebra to have the SLP. This theorem, due to Watanabe, is stated in Theorem \ref{watanabe}.
	    
	    \pause
	    \textbf{Chapter 3} is the principal part of this thesis. As we have said, in Section \ref{catal}, we introduce the catalecticant matrices in the more general frame of the Waring problem. To this aim we introduce also the \emph{symmetric border rank}. We show the connection between the catalecticant matrices of a certain form $h\in K[u,v]_t$ and the position of $[h]$ in $\PP(K[u,v]_t)$ with respect to the rational normal curve $C_t$ and its secant varieties $\sigma_r(C_t)$. The most noticeable one is Proposition \ref{rank} where forms with catalecticant matrices of low rank are classified. As a result, we are able to prove Propositions \ref{upper} and \ref{lower} giving an upper and a lower bound to the possible Hilbert vectors of a Perazzo $3$-fold. We are also able to give a characterisation of Perazzo $3$-fold with minimal Hilbert vector in Theorem \ref{lowcharacterisationreprise}; instead Example \ref{esempiodellavita} gives examples of Perazzo $3$-folds with maximal $h$-vector. In Section \ref{wlp}, we consider the Artinian Gorenstein algebra $A$ associated to a Perazzo $3$-folds, which always lacks the SLP, and we prove the following: if the Hilbert vector is maximal, then $A$ fails the WLP (Proposition \ref{main1}); instead if the Hilbert vector is minimal, then $A$ always has the WLP (Proposition \ref{main2}). At the end of the chapter, we briefly treat the general case of hypersurfaces with or without zero Hessian with general degree $d$.
	    
	    \pause
	    In \textbf{Chapter 4}, we recall some arguments treated in Chapter 1 regarding the geometry of hypersurfaces with vanishing Hessian. In particular, we give some definitions about the dual variety of a  variety $X$ and its connection with $X$. Results like Corollary \ref{corollariotuttodiversonienteuguale} and Proposition \ref{proposizionetecnica} try to give a geometrical meaning to the condition of zero Hessian. The main result is Theorem \ref{cfrs} by Ciliberto, Franchetta, Russo, and Simis where we summarise the results about the geometry of hypersurfaces with vanishing Hessian in $\PP^4$. Eventually, we specialise Theorem \ref{cfrs} to Perazzo $3$-folds; in particular, we give explicit computations about Perazzo $3$-folds with minimal Hilbert vector.

	\pagebreak
    \cleardoublepage\phantomsection
    \addcontentsline{toc}{chapter}{Acknowledgements}
    \chapter*{Acknowledgements}
    
    I would like to thank my supervisor Professor Emilia Mezzetti for her patience, professionalism, assistance, and inspiration during my University journey, starting from my very first exam, and especially in this last year. I thank her for trusting me, supporting me, and for all the learning possibilities she has offered me. I am grateful to be able to conclude my Master degree with such an amazing teacher.
   
    \pause
    I would like to thank Professor Rosa Maria Mirò-Roig as well, for the two months' work here in Trieste: it has been a great honour to be involved in her and Professor Mezzetti's research, and to be able to publish my first article alongside them.

    \pause
    Eventually, I thank Lisa Nicklasson for her comments on Theorem \ref{main2}.
	    
	
	\chapter{Hypersurfaces with vanishing Hessian}

	This first chapter concerns the study of hypersurfaces with vanishing Hessian. The principal paper that we follow is due to J. Watanabe and M. de Bondt \cite{watanabe2019theory}. Their principal goal was to rediscover and analyse the classical paper of Gordan and Noether \cite{gordannoether}. Using the so called \emph{self-vanishing systems}, that we will define in Definition \ref{selfvanishingsystem}, it is possible to simplify Gordan-Noether's theorems. In this thesis we focus only on homogeneous polynoamials because we have a correspondence with hypersurfaces in the projective space. We simply name homogeneous polynomials as \emph{forms}. A different approach can be found also in \cite{pirolafavalebricalli}, \cite{debondt1} and \cite{debondt2}. A more complete bibliography and historical notes can be found in the last section of this chapter. We start now giving some definitions and some properties that characterise forms with vanishing Hessian. We remark that we are working only with algebraically closed fields of characteristic zero. 
	
	
	
	\begin{definition}
	    Let $R=K[x_0,\dots,x_n]$ be the polynomial ring in $n+1$ variables over the field $K$. Given $f\in R$, a homogeneous polynomial of degree $d$, we define the \emph{Hessian matrix of $f$} to be
	    $$\Hess(f):=\left(\frac{\partial^2f}{\partial x_i\partial x_j}\right)_{1\le i,j\le n}\in (R_{d-2})^{(n+1)\times (n+1)}.$$
	    We say that $f$ has \emph{vanishing Hessian}, or equivalently $f$ is a form with \emph{zero Hessian}, if $$\hess(f):=\det\Hess(f)\in R_{(d-2)(n+1)}$$
	    is the zero polynomial.
	    
	    As an abuse of notations, we say that a hypersurface $X=V(f)\subseteq\PP^{n}$ has vanishing Hessian if $f$ has vanishing Hessian.
	\end{definition}
	
	We can now characterise this type of forms using the following Proposition:
	\begin{proposition}\label{partialalgebraic}
	    Let $f\in K[x_0,\dots,x_n]$ be a form. Then the following are equivalent:
	    \begin{itemize}
	        \item $f$ has vanishing Hessian;
	        \item the partial derivatives of $f$ are algebraically dependent.
	    \end{itemize}
	\end{proposition}
	
	This fact can be proved using the following more general theorem concerning the relation between the Jacobian matrix and algebraic relations of forms of the same degree.
	\begin{theorem}\label{jacobian}
	    Let $R=K[x_0,\dots,x_n]$. Let $\mathbf{f}=(f_0,f_1,\dots,f_n)$ be a $(n+1)$-tuple of forms of the same degree. Then $\rank (\frac{\partial f_i}{\partial x_j})=\text{tr.deg}_K K(f_0,\dots,f_n)$. In particular the following conditions are equivalent for any fixed integer $r\in\{1,2,\dots,n\}$:
	    \begin{itemize}
	        \item $f_{i_1},\dots,f_{i_r}$ are algebraically dependent for every $i_1,\dots,i_r\in \{0,1,\dots,n\}$;
	        \item the rank of the Jacobian matrix $(\frac{\partial f_i}{\partial x_j})$ is strictly less than $r$;
	        \item $\text{tr.deg}_K K(f_{i_1},\dots,f_{i_r})<r$ for every $i_1,\dots,i_r\in \{1,2,\dots,n\}$.
	    \end{itemize}
	\end{theorem}
	\begin{proof}
	    See \cite[Theorem 8]{M-2013}.
	\end{proof}
	
	In the case $\mathbf{f}=(\frac{\partial f}{\partial x_0},\frac{\partial f}{\partial x_1},\dots,\frac{\partial f}{\partial x_n})$ we can see that Prop \ref{partialalgebraic} is a simple corollary of Theorem \ref{jacobian}. In particular, the transcendence degree of the extension of $K$ generated by partial derivatives says the dimension of the null space of the Hessian matrix. As a notation, we will write $f_i$ to refer to the partial derivative of $f$ with respect to $x_i$. Later in this chapter we will see some theorems that are true for a generic $(n+1)$-tuple of polynomials, as in the previous theorem. In this case $\mathbf{f}=(f_0,f_1,\dots,f_n)$ will refer simply to a \emph{system} of forms in $K[x_0,\dots,x_n]$.
	
	\pause
	Considering that the algebraic relations play a key role, it is useful to give the following definition.
	
	\begin{definition}\label{If}
	    Let $f\in K[x_0,\dots,x_n]$ be a form with zero Hessian. Let $y=(y_0,\dots,y_n)$ be a coordinate system independent of $x=(x_0,\dots,x_n)$ and let
	    $$\phi:K[y_0,\dots,y_n]\to K[x_0,\dots,x_n]$$
	    be the $K$-homomorphism defined by $$\phi(g(y))=g\left(\frac{\partial f}{\partial x_0}(x),\dots,\frac{\partial f}{\partial x_n}(x)\right).$$ We denote by $\I(f)$ the kernel of $\phi$. $\I(f)$ is a non zero prime ideal containing all possible algebraic relations among the partials of $f$ (by Proposition \ref{partialalgebraic}).
	\end{definition}
	
	
	It is an open problem to classify all possible forms with vanishing Hessian. Over many years, a lot of mathematicians tried to give a description of this type of forms. The first attempt was made by O. Hesse in \cite{hesse1} and \cite{hesse2}. Here he claimed that "any hypersurface $X\subset \PP^{n}$ with vanishing hessian is a cone". 
	
	\begin{definition}
	    A hypersurface $X=V(f)\subset\PP^{n}$ is a \emph{cone} if the following equivalent conditions hold:
	    \begin{itemize}
	        \item there exists a point $Q\in X$ such that the line joining $Q$ and another point $P\in X$, $P\ne Q$, is always contained in $X$,
	        \item a variable can be eliminated from $f$ by means of a linear transformation of the variables.
	    \end{itemize}
	    
	    With an abuse of notation, when speaking of a cone, we will refer either to the hypersurface or to the form defining it.
	\end{definition}
	
	It is clear that every form that defines a cone has vanishing Hessian, but the opposite implication is not so clear. In fact, Gordan-Noether disproved O. Hesse claim proving that the condition of the partial derivatives to be algebraically dependent is not sufficient for the form to be a cone. In fact, we have the following proposition.
	
	
	\begin{proposition}\label{basicelimination}
	    Suppose that $f\in K[x_0,\dots,x_n]$ is a form with zero Hessian. Let $\I(f)$ be the ideal of $K[y_0,\dots,y_n]$ as defined in Definition \ref{If}. Then the following conditions are equivalent:
	    \begin{enumerate}
	        \item the ideal $\I(f)$ contains a linear form,
	        \item the partial derivatives of $f$ are linearly dependent,
	        \item a variable can be eliminated from $f$ by means of a linear transformation of the variables.
	    \end{enumerate}
	\end{proposition}
	\begin{proof}
	    Knowing that $\I(f)$ contains the algebraic relations between the partials of $f$, the equivalence between 1 and 2 is clear. Suppose now there exists a non-trivial relation
	    $$a_0f_0+a_1f_1+\dots+a_nf_n=0$$
	    where $f_i=\frac{\partial f}{\partial x_i}$ and $a_i\in K$. Without loss of generality, we can suppose that $a_0\ne 0$. We now define a set of linearly independent linear forms $l_0,\dots,l_n$ as follows: $l_0=\frac{1}{a_0}x_0$ and $l_i=x_i-\frac{a_i}{a_0}x_0$ for $i=1,\dots,n$. Then $\frac{\partial x_i}{\partial l_0}=a_i$ and
	    $$\frac{\partial f}{\partial l_0}=\frac{\partial f}{\partial x_0}\frac{\partial x_0}{\partial l_0}+\frac{\partial f}{\partial x_1}\frac{\partial x_1}{\partial l_0}+\dots+\frac{\partial f}{\partial x_n}\frac{\partial x_n}{\partial l_0}=0.$$
	    This shows that if $f$ is expressed in terms of $l_i$, then $f$ does not contain $l_0$. Thus 2 implies 3. The same argument shows the opposite implication as well.
	\end{proof}
	
	\begin{example}\label{perazzo5}
    	The first example that Gordan-Noether gave in their paper was $f=u^2x+uvy+v^2z\in K[x,y,z,u,v]$. Later, this hypersurface was named "Perazzo hypersurface". The partial derivatives and the Hessian matrix can be easily computed as follow:
    	$$\mathbf{f}=(u^2,uv,v^2,2ux+vy,uy+2vz)$$
    	$$\Hess(f)=\left(\begin{array}{ccccc}
    	     0 & 0 & 0 & 2u & 0 \\
    	     0 & 0 & 0 & v & u \\
    	     0 & 0 & 0 & 0 & v^2 \\
    	     2u & v & 0 & 2x & y \\
    	     0 & u & 2v & y & 2z
    	\end{array}\right)\ .$$
    	We see that $f$ has zero Hessian, but no non-trivial linear relation is satisfied by the partial derivatives. So, as we have state above, $V(f)$ is a hypersurface, not a cone, with vanishing Hessian.
    	
    	\pause
    	We notice that $\Hess(f)$ has a $3\times 3$ all zero sub-matrix. This remark will inspire Theorem \ref{main1} to prove that a certain matrix has zero determinant.
	\end{example}
	
	\begin{example}["Un esempio semplicissimo"]
	    \label{unesempiosemplicissimo}
	    This class of examples was first studied by U. Perazzo in \cite{perazzo}. We can generalise Example \ref{perazzo5} with the form $f=x_0x_3^2+x_1x_3x_4+x_2x_4^2+x_5^3+\dots+x_n^3$ which involves $n+1$ variables, $n\ge4$. The partials are again linearly independent, but they vanish in the form $g(y)=y_0y_2-y_1^2$.
	\end{example}
	
	We now ask ourselves: is it possible to classify all forms with vanishing Hessian up to linear transformations? Gordan-Noether answered to this question up to the dimension $n=4$. In particular they proved that for binary, ternary and quaternary forms Hesse's claim is true. Instead, for a form $f\in K[x,y,z,u,v]$ with zero Hessian, they proved that $f$ can be written as an element of $K[u,v][\Delta]$ where $\Delta$ is a homogeneous polynomial in five variables of the form $\Delta=p_1(u,v)x+p_2(u,v)y+p_3(u,v)z$ . On the other hand, they also proved that every hypersurface with zero Hessian turns out to be birational equivalent to a cone. Furthermore, it is possible to explicitly construct such a birational map which is often tricky for general classification problems.
	
	\pause
	Example \ref{perazzo5} is going to be crucial in the third chapter, where we will extend the definition of Perazzo hypersurface, and in the classification of hypersurfaces with vanishing Hessian in five variables.
	
	\pause
	In order to prove all these statements we have to talk about "self-vanishing systems". They will be described and studied in a general way in the next section. We will not write down all the details, which can be found in \cite{watanabe2019theory}, but only the most noticeable ones for our purposes.
	
	\section{Self-vanishing systems}\label{selfvanishingsystems}
	Let $x=(x_0,\dots,x_n)$ and $y=(y_0,\dots,y_n)$ be two sets of indeterminates and let $K(x,y)=K(x_0,\dots,x_n,y_0,\dots,y_n)$ denote the rational function field. We want to introduce a differential operator that acts as the canonical scalar product between the set of indeterminates $y$ and the gradient of $f$ with respect to $x$. So, we introduce the differential operator
	$$\D_x(y):K(x,y)\to K(x,y)$$
	which is defined by
	$$\D_x(y)f(x,y):=\sum_{j=0}^n y_j\frac{\partial f}{\partial x_j}(x,y)$$
	for $f(x,y)\in K(x,y)$. This definition is for a general rational function, but we are mostly interested in homogeneous polynomials in only the $x$ variables. In fact, for a homogeneous polynomial $f(x)=f(x_0,\dots,x_n)\in K[x_0,\dots,x_n]$ and for $j\ge 0$, we define $f^{(j)}(x,y)$ to be the polynomial in $K[x_0,\dots,x_n,y_0,\dots,y_n]$ given by
	\begin{equation}\label{fj}
	    f^{(j)}(x,y)=\frac{1}{j!}\D_x(y)^j(f(x)).
	\end{equation}
	It is remarkable that $f^{(j)}(x,y)$ is homogeneous in both systems of coordinates, but in general with different degrees. Moreover, if it is non-zero we get $\deg_y f^{(j)}=j$.
	\begin{proposition}\label{proprietaSVS}
	    Let $f(x)\in K[x]$ be a homogeneous polynomial of degree $d$, then the following equalities hold:
	    \begin{itemize}
	        \item $f^{(j)}(x,x)=\binom{d}{j}f(x)$ for every $j=1,\dots,d$,
	        \item $f^{(d)}(x,y)=f(y)$,
	        \item $f(x+ty)=\sum\limits_{j=0}^d t^jf^{(j)}(x,y)$ for an indeterminate $t$.
	    \end{itemize}
	\end{proposition}
	\begin{proof}
	        See \cite[Proposition 3.2 and before]{watanabe2019theory}.
	\end{proof}
	
	\begin{definition}
	    We now consider $\mathbf{h}=(h_1(x),\dots,h_n(x))$ a non-trivial system of homogeneous polynomials $h_i\in K[x_0,\dots,x_n]$ of the same degree. We now define the differential operator $\D_x(\mathbf{h}):K(x)\to K(x)$ associate to $\mathbf{h}$ by
	    $$\D_x(\mathbf{h})f(x):=\sum_{j=0}^n h_j(x)\frac{\partial f}{\partial x_j}(x)=f^{(1)}(x,\mathbf{h}(x)).$$
	    Furthermore, we denote by $\sol(\mathbf{h};S)$ the set of solutions in the subring $S\subseteq K[x]$ of the partial differential equation
	    \begin{equation}\label{equazioneuno}
	        \D_x(\mathbf{h})f(x)=0.
	    \end{equation}
	    More precisely, $\sol(\mathbf{h};S)=\left\{ f(x)\in S:\,\sum\limits_{j=0}^n h_j(x)\frac{\partial f}{\partial x_j}(x)=0 \right\}$. 
	    
	    Note that $\sol(\mathbf{h};K[x])$ is a graded subalgebra of $K[x]$, a fact which is not true for a general subring $S$. For this reason we are most interested in the case $S=K[x]$ and the set of solutions are simply named $\sol(\mathbf{h})=\sol(\mathbf{h};K[x])$.
	\end{definition}
	\begin{definition}\label{selfvanishingsystem}
	    A system $\mathbf{h}=(h_0(x),\dots,h_n(x))$ of forms is called \emph{self-vanishing} if $h_i\in\sol(\mathbf{h})$ for all $j=0,1,\dots,n$. In addition we also say that $\mathbf{h}$ is a \emph{reduced self-vanishing} system if $\gcd(h_0,\dots,h_n)=1$.
	\end{definition}
	
	$\sol(\mathbf{h})$ is a prime ideal of $K[x]$. Furthermore, if $f,g$ are two elements such that their product is in $\sol(\mathbf{h})$, then both $f,g\in\sol(\mathbf{h})$ (see \cite[Corollary 3.12]{watanabe2019theory}). For this reason, it is always possible to obtain a reduced self-vanishing system from a general self-vanishing system dividing by the greatest common divisor.
	
	\begin{example}
	    A constant vector $\mathbf{h}=(a_0,\dots,a_n)\in K^{n+1}$ is clearly a self-vanishing system.
	\end{example}
	\begin{example}[GN type self-vanishing system]
	    We now suppose that $\mathbf{h}$ satisfies the following conditions:
	    \begin{enumerate}
	        \item $h_0=h_1=\dots=h_r=0$ for some integer $r$, $0\le r\le n$,
	        \item the polynomials $h_{r+1},\dots,h_n$ involve only the variables $x_0,\dots,x_r$.
	    \end{enumerate}
	    Then clearly $\mathbf{h}$ is a self-vanishing system of forms. We name them "Gordan-Noether type" (GN type for short).
	\end{example}
	
	\begin{notation}
	    Given a system of forms $\mathbf{f}=(f_0,f_1,\dots,f_n)$ and an index $j$, $0\le j\le n$, we denote by $\partial \mathbf{f}\mathbin{/}{\partial x_j}$ the system of forms defined by
	    $$\frac{\partial \mathbf{f}}{\partial x_j}:=(\frac{\partial f_0}{\partial x_j},\frac{\partial f_1}{\partial x_j},\dots,\frac{\partial f_n}{\partial x_j}).$$
	\end{notation}
	
	\begin{proposition}\label{definizionedih}
	    Let $\mathbf{f}=(f_0,f_1,\dots,f_n)$ be a system of forms in $K[x]$, in which the components are algebraically dependent. Let $y=(y_0,\dots,y_n)$ be a coordinate system independent of $x$. Let $\phi:K[y]\to K[x]$ be the $K$-homomorphism defined by
	    $$y_j\longmapsto f_j, \, (j=0,1,\dots,n).$$
	    Let $g=g(y)\in \ker\phi$ be a non-zero element which we can suppose to be homogeneous and of minimal degree (as in Def. \ref{If}). We now define $h_j\in K[x]$ by
	    $$h_j(x):=\phi\left(\frac{\partial g}{\partial y_j}\right)=\frac{\partial g}{\partial y_j}(f_1(x),\dots,f_n(x)).$$
	    Let $\mathbf{h}=(h_0(x),\dots,h_n(x))$. Then:
	    \begin{itemize}
	        \item $\mathbf{h}$ is a syzygy of $\mathbf{f}$, i.e. $\sum_{j=1}^n h_jf_j=0$,
	        \item $\mathbf{h}$ is a syzygy of $\partial\mathbf{f}\mathbin{/}\partial x_j$ for every $j=1,2,\dots,n.$
	    \end{itemize}
	\end{proposition}
	\begin{proof}
	        See \cite[Proposition 3.8]{watanabe2019theory}.
	\end{proof}
	\begin{remark}\label{unoèdiversodazero}
	    In Proposition \ref{definizionedih}, it is possible that $h_j=0$ for some $j$. This means that $g$ does not depend on $x_j$, by the minimality of $g$ itself. Hence, $\mathbf{h}$ is not trivial. If we drop the condition of minimality of the degree of $g$, it could be that $\mathbf{h}$ is trivial. Proposition \ref{definizionedih} is in any case true, but we are not interested in this case.
	    
	    Moreover, the choice of $g$ is unique if and only if $\text{tr.deg}_K K(f_0,\dots,f_n)=n$. The strong fact is that this proposition is always true for every possible choice of such a polynomial $g$.
	\end{remark}
	\begin{example}
	    In Proposition \ref{definizionedih}, $\mathbf{h}$ is not in general a self-vanishing system. Consider for example $\mathbf{f}=(x_0^2,x_1^2,ix_0x_1)$ a system of forms in $\C[x_0,x_1,x_2]$ where $i$ is the imaginary unit. Then the polynomial $g$ can be taken to be $g(y_0,y_1,y_2)=y_0y_1+y_2^2$. This is homogeneous and of minimal degree by the linear independence of the $f_i$s. Then we obtain the system $\mathbf{h}=(x_1^2,x_0^2,2ix_0x_1)$ which is not a self-vanishing system. In fact, computing $\D_x(\mathbf{h})h_0(x)$ we obtain
	    $$\D_x(\mathbf{h})h_0(x)=x_1^2\cdot 2x_1\ne 0.$$
	    One can easily check that the relations stated in Proposition \ref{definizionedih} are nevertheless true.
	\end{example}
	
	 In contrast with this example, we now state a sufficient condition for $\mathbf{h}$ to be a self-vanishing system.
	 
	\begin{theorem}\label{svs}
	    Let $f=f(x)\in K[x]$ be a homogeneous polynomial and define $f_j=\frac{\partial f}{\partial x_j}$. Assume that $f_0,\dots,f_n$ are algebraically dependent (namely, $f$ has vanishing Hessian). Let $\mathbf{f}=(f_0,f_1,\dots,f_n)$ and let $\mathbf{h}=(h_0(x),\dots,h_n(x))$ be a system of forms defined as in Proposition \ref{definizionedih}. Then:
	    \begin{enumerate}
	        \item $f(x)\in\sol(\mathbf{h})$,
	        \item $f_j(x)\in\sol(\mathbf{h})$ for $j=0,\dots, n$,
	        \item $f(x)\in\sol(\partial\mathbf{h}\mathbin{/}\partial x_j)$ for $j=0,\dots, n$,
	        \item $\mathbf{h}$ is a self-vanishing system.
	    \end{enumerate}
	    Moreover, it is possible to choose $\mathbf{h}$ to be a reduced self-vanishing system.
	\end{theorem}
	\begin{proof}
        The principal tool is Proposition \ref{definizionedih}. Indeed, 1 is a direct consequence of it. Moreover, being $\mathbf{h}$ a syzygy of $\partial\mathbf{f}\mathbin{/}{\partial x_j}$ and knowing that the partials satisfy the Schwarz equality, we obtain
        $$\sum_{j=0}^n h_j(x)\frac{\partial f_k}{\partial x_j}(x)=\sum_{j=0}^n h_j(x)\frac{\partial f_j}{\partial x_k}(x)=0.$$
        So, 2 is true. We now consider the product $\mathbf{h}\cdot\mathbf{f}$ which is zero by Proposition \ref{definizionedih}. We now compute the derivative with respect to $x_j$ applying the chain rule:
        $$\frac{\partial}{\partial x_j}(\mathbf{h}\cdot\mathbf{f})=\frac{\partial \mathbf{h}}{\partial x_j}\cdot\mathbf{f}+\mathbf{h}\cdot\frac{\partial\mathbf{f}}{\partial x_j}=0.$$
        As we said, the second addend is zero, so $f(x)\in\sol(\partial\mathbf{h}\mathbin{/}\partial x_j)$. The last statement can be argued as follow: the polynomials $h_i$ are each defined as an algebraic combination of the elements $f_0,\dots,f_n$, so they are all elements of the extension $K[f_0,\dots,f_n]$. But, $K[f_0,\dots,f_n]\subseteq\sol(\mathbf{h})$ because $\sol(\mathbf{h})$ is a $K$-algebra and all $f_i$ are solution of the PDE (\ref{equazioneuno}). Finally, $\mathbf{h}$ is, by definition, a self-vanishing system.
        
        If we divide each $h_i$ by $\gcd(h_0,\dots,h_n)$, then clearly, all statements are again true. Furthermore, as we discussed above, the resulting system of forms is a reduced self-vanishing system.
	\end{proof}
	
	The next results are going to be stated for a general self-vanishing system. In the next section we will combine them with Theorem \ref{svs} to obtain our classification.
	
	\begin{theorem}[Gordan-Noether identity]\label{gnidentity}
	    Suppose that $\mathbf{h} = (h_0(x),\dots,h_n(x))$ is a self-vanishing system of forms in $K[x]$. Then, for a homogeneous polynomial $f(x)\in K[x]$ and considering $f^{(j)}(x,\mathbf{h}(x))$ as defined in (\ref{fj}), the following conditions are equivalent:
	    \begin{enumerate}
	        \item $f(x)\in\sol(\mathbf{h})$,
	        \item $f^{(j)}(x,\mathbf{h}(x))=0$ for $j=1,2,\dots$,
	        \item $f$ satisfies the identity of polynomials $f(x+t\mathbf{h}(x))=f(x)$ for any $t\in K'$, where $K'$ is any extension field of $K$. This identity is called \emph{Gordan-Noether identity}.
	    \end{enumerate}
	\end{theorem}
	\begin{proof}
	    The equivalence between 2) and 3) is a direct consequence of Proposition \ref{proprietaSVS}. To prove that 1) is equivalent to 2) we need two facts. The first one is the following link between $\D_x(\mathbf{h})f^{(j)}$ and $f^{(j+1)}$ which is proved in \cite[Proposition 3.10]{watanabe2019theory}:
	    $$\D_x(\mathbf{h})(f^{(j)}(x,\mathbf{h}))=(i+1)f^{(j+1)}(x,\mathbf{h}).$$
	    So, it is enough to show the case $j=1$, but this is trivial knowing the following second fact
	    \begin{align*}
	        f^{(1)}(x,\mathbf{h}(x)) &= \D_x(y)f(x)|_{y=\mathbf{h}}\\
	        &= \D_x(\mathbf{h})f(x)\\
	        &= h_1(x)\frac{\partial f}{\partial x_1}(x)+h_2(x)\frac{\partial f}{\partial x_2}(x)+\dots+h_n(x)\frac{\partial f}{\partial x_n}(x). \qedhere
	    \end{align*}
	\end{proof}
	
	We are now able to state the most important result which is the purpose of all this theory.
	
	\begin{corollary}\label{proofgordannoether}
	    Suppose that $\mathbf{h}=(h_0(x),\dots,h_n(x))$ is a self-vanishing system in $K[x]$ such that $h_n(x)\ne 0$. Put
	    $$s_i(x)=x_i-\frac{h_i(x)}{h_n(x)}x_n, \text{ for } 0\le i\le n-1.$$
	    Let $f(x)\in\sol(\mathbf{h},K[x])$. Then
	    \begin{equation}\label{equazioneproof}
	        f(x)=f(s_0,\dots,s_{n-1},0)
	    \end{equation}
	    and, in particular,
	    $$\sol(\mathbf{h})=K[s_0,\dots,s_{n-1}]\cap K[x].$$
	\end{corollary}
	\begin{proof}
	    First of all we prove that the polynomial $f$ can be written as in (\ref{equazioneproof}). From Proposition \ref{gnidentity} and choosing $t=-\frac{x_n}{h_n(x)}\in K(x)$ we obtain
	    $$f(x)=f(x+t\mathbf{h}(x))=f\left(x-\frac{x_n}{h_n(x)}\right)=f(s_0,\dots,s_{n-1},0).$$
	    This also proves that $K[s_1,\dots,s_{n-1}]\cap K[x]\supseteq\sol(\mathbf{h})$.
        To prove the other implication, we can proceed as follow. Being $\sol(\mathbf{h})$ a $K$-algebra, we have just to prove that $s_j\in\sol(\mathbf{h})$. Let $A$ the $2\times 2$ matrix, with entries in $K[x]$, defined as $A=\begin{bmatrix}
        x_j & x_n\\
        h_j(x) & h_n(x)
        \end{bmatrix}$. We now compute $\D_x(\mathbf{h})s_j(x)= \D_x(\mathbf{h})\left(\frac{1}{h_n(x)}\det A\right)$ for every index $j$:
	    \begin{align*}
	        \D_x(\mathbf{h})s_j(x)&=\sum_{i=1}^n h_i(x)\frac{\partial}{\partial x_i}\left(\frac{1}{h_n(x)}\det A\right)\\
	        &=\sum_{i=1}^n \frac{h_i(x)}{h_n(x)^2} \left(\tr(\text{adj}A\cdot \frac{\partial A}{\partial x_i}) h_n(x)-\det A \frac{\partial h_n(x)}{\partial x_i}\right)\\
	        &=\sum_{i=1}^n \frac{h_i(x)}{h_n(x)^2}\left[x_n\left(h_j\frac{\partial h_n}{\partial x_i}-h_n\frac{\partial h_j}{\partial x_i}\right)+\delta_j^ih_n-\delta_n^ih_j\right]\\
	        &=\frac{x_n}{h_n(x)^2}\left[h_j\sum_{i=1}^nh_i\frac{\partial h_n}{\partial x_i}-h_n\sum_{i=1}^nh_i\frac{\partial h_j}{\partial x_i}\right]=0.\qedhere
	    \end{align*}
	\end{proof}
	
	\section{Classification of forms with zero Hessian}\label{classificationofformswithzerohessian}
	\subsection{Birational classification}
	
	We now proceed with the classification of forms with zero Hessian. As we discussed above, if we consider birational transformations then we obtain the following theorem.
	\begin{theorem}[Gordan-Noether]\label{gntheorem}
	    Let $f\in K[x_0,\dots,x_n]$ be a form with zero Hessian. Then $f$ is Cremona equivalent to a cone, i.e. there exist a set of variables $\{s_0,\dots,s_n\}$, and a birational map $\psi:K(x_0,\dots,x_n)\to K(s_0,\dots,s_n)$ such that $\psi(f)$ does not involve the variable $s_n$.
	\end{theorem}
	\begin{proof}
	    By Proposition \ref{partialalgebraic}, the partial derivatives of $f$ are algebraically dependent. So, we can apply Proposition \ref{definizionedih} and Remark \ref{unoèdiversodazero} to obtain a self-vanishing system $\mathbf{h}=(h_0,\dots,h_n)$ such that, up to a change of variables, $h_n\ne 0$. We can now apply Corollary \ref{proofgordannoether} to obtain a rational map $\psi^{-1}:K(s_0,\dots,s_n)\to K(x_0,\dots,x_n)$ defined as
	    $$\psi^{-1}(s_i)=x_i-\frac{h_i(x)}{h_n(x)}x_n, \text{ for } 0\le i\le n-1,\quad s_n=x_n.$$
	    Moreover, by Corollary \ref{proofgordannoether}, $f$ written with the $s$ variables does not depend on $s_n$. To conclude, we explicit construct the inverse map. Noticing that $\mathbf{h}$ is a self-vanishing system, Corollary \ref{proofgordannoether} can be applied also to each $h_i$ so that $h_i(x)=h_i(s)$. So we can write
	    \begin{equation*}
	        \psi(x_i)=s_i+\frac{h_i(x)}{h_n(x)}x_n=s_i+\frac{h_i(s)}{h_n(s)}s_n.\qedhere
	    \end{equation*}
	\end{proof}
	\begin{example}
	    As an application of Theorem \ref{gntheorem}, we study the case of the Perazzo hypersurface $f=u^2x+uvy+v^2z$. For the computations of Example \ref{perazzo5} we can consider $g(y)=y_0y_2-y_1^2$ which defines the system $\mathbf{h}=(v^2,-2uv,u^2,0,0)$. Noticing that $h_0\ne 0$, the Cremona transformations $\psi:K(x,y,z,u,v)\to K(x',y',z',u',v')$ can be defined as
	    \begin{equation*}
	        \begin{cases}
	                x'=x\\
	                y'=y+\frac{2u}{v}x\\
	                z'=z-\frac{u^2}{v^2}x\\
	                u'=u\\
	                v'=v
	        \end{cases}
	    \end{equation*}
	    Thus $\psi(f)=v'(y'u'+z'v')$. Instead, if we choose a different element $h_1=-2uv\ne 0$ we obtain another Cremona map $\Tilde{\psi}: K(x,y,z,u,v) \to  K(x',y',z',u',v')$ defined by
	    \begin{equation*}
	        \begin{cases}
	                x'=x+\frac{v^2}{2uv}y\\
	                y'=y\\
	                z'=z+\frac{u^2}{2uv}y\\
	                u'=u\\
	                v'=v
	        \end{cases}
	    \end{equation*}
	    In this case we obtain a different cone, indeed $\Tilde{\psi}(f)=x'(u')^2+z'(v')^2$. So, different choices of $h_i$ can give different cones birational to $f$.
	\end{example}
	\subsection{Classification of the cases $n=1,2,3,4$}
	
	We now focus on linear transformations, so that the problem can be rephrased in this way. Given a form $f\in K[x_0,\dots,x_n]$ of degree $d$, the main goal is to understand whenever there exists or not an invertible matrix $A\in GL(K^{n+1})$ such that the new form $f'(x_0',\dots,x_n')=f(A^{-1}x)$ does not depend on $x_{k}',\dots,x_{n}'$ (for some index $k$). If this happens we say that $n+1-k$ variables can be eliminated from $f$. Clearly as we said above about the Gordan-Noether's paper, this is not possible in general and other kinds of classes will occur. The main basic tools are summarised in the following two lemmas.
	\begin{lemma}
	    Let $f$ be a form in $n+1$ variables. Let $f_i$ be its partial derivatives, which define a system of forms $\mathbf{f}=(f_0,\dots,f_n)$. As we have seen in Proposition \ref{basicelimination}, if the partials are linearly dependent, then one variable can be eliminated. Furthermore, if $\dim_K \bigoplus_{i=0}^n Kf_i=s$, then $n+1-s$ variables can be eliminated from $f$.
	\end{lemma}
	\begin{proof}
	    See \cite[Lemma 2.2]{watanabe2019theory}.
	\end{proof}
	\begin{lemma}\label{elimintionsuperpower}
	    Let $\mathbf{f}=(f_0,\dots,f_n)$ be a any system of forms. Put $g$ and $\mathbf{h}=(h_0,\dots,h_n)$ as defined in Proposition \ref{definizionedih}. Furthermore let $\dim\bigoplus_{i=0}^{n}Kh_j=s$. Then $n+1-s$ variables can be eliminated from $g$. 
	\end{lemma}
	\begin{proof}
	    See \cite[Proposition 2.3]{watanabe2019theory}.
	\end{proof}
	
    We proceed studying the easiest cases:
    \begin{itemize}
        \item for $n=0$ we just have a form with one variable and vanishing Hessian. It is a basic fact that $f$ is of the form $f=ax_0$ which is indeed a cone.
        \item Another obvious case is when $f$ defines a quadric hypersurface, i.e. $d=2$. In that case, $f$ can be written as $f=\mathbf{x}^{T}\cdot B\cdot \mathbf{x}$ for some symmetric matrix $B\in K^{(n+1)\times(n+1)}$. One can notice that $B$ is indeed the Hessian matrix of $f$ and, via the classification of quadric hypersurfaces, we can write $f=x_0^2+\dots+x_r^2$ where $r=\rank B<n+1$.
    \end{itemize}
    
    For $n=1,2$ Hesse's claim is true and we have the following Theorems \ref{binary} and \ref{ternary}.
    \begin{theorem}[binary forms]\label{binary}
        Let $f\in K[x_0,x_1]$ be a form of degree $d$ with vanishing Hessian. Then $f$ is a pure power of a linear form, i.e. $f=l^d$ for some linear form $l\in K[x_0,x_1]_1$.
    \end{theorem}
    \begin{proof}
        Let $g\in\I(f)$ be an algebraic relation of the partials of $f$. It can be taken homogeneous and of minimal degree; then, being $g$ a binary form in a algebraically closed field $K$, it must be a linear form. By Proposition \ref{basicelimination}, the proof is complete.
    \end{proof}
    
    For the case $n=2$ we have to introduce a lemma whose proof can be found in \cite[Lemma 5.2]{watanabe2019theory}.
    
    \begin{lemma}\label{lemmatrinaryform}
        Let $\mathbf{h}=(h_0,\dots,h_n)$ be a self-vanishing system. Then it holds
        $$\text{tr.deg}_K K(h_0,\dots,h_n)\le n.$$
        In particular, if $n\ge 2$, we get the stronger relation
        $$\text{tr.deg}_K K(h_0,\dots,h_n)\le n-1.$$
    \end{lemma}
    \begin{theorem}[ternary forms]\label{ternary}
        Let $f\in K[x_0,x_1,x_2]$ be a form with vanishing Hessian. Then a variable can be eliminated from $f$ by a linear transformation of the variables.
    \end{theorem}
    \begin{proof}
        From Theorem \ref{svs} we get $\mathbf{h}=(h_0,h_1,h_2)$ a reduced self-vanishing system of forms. We want to prove that it is a constant vector. In fact, it will imply that the partial derivatives of $f$ are linearly dependent and thus a variable can be eliminated from $f$ by Proposition \ref{basicelimination}. By contradiction, we suppose that $\mathbf{h}$ is not a constant vector which is equivalent to say $\text{tr.deg}_K K(h_0,h_1,h_2)\ge 1$. By Lemma \ref{lemmatrinaryform} it must be an equality. So, every possible pair of elements in $\{h_0,h_1,h_2\}$ is algebraically dependent and, as we argued in Theorem \ref{binary}, it must be a linear relation. We recall that $\{h_0,h_1,h_2\}$ are also coprime, so the only possibility is that $\mathbf{h}$ is a constant vector.
    \end{proof}
    
    To treat the cases of quaternary and quinary forms we have to introduce a new argument regarding system of forms and rational maps.
    
    \begin{definition}\label{mapZ}
        Given a system of forms $\mathbf{h}=(h_0,\dots,h_n)$ of the same degree, we can define a rational map
        $$\Phi:\PP^n \dashrightarrow\PP^n$$
        defined as $\Phi([x_0:\dots:x_n])=[h_0(x):\dots:h_n(x)]$. We name $T$ and $W$ respectively the fundamental locus and the image of $\Phi$. Gordan and Noether called this function \emph{"die Functionen"}.
    \end{definition}
    
    Being the forms of the same degree, $\Phi$ is a well defined rational map. Moreover, if $\mathbf{h}$ is a self-vanishing system, it does not matter if the system is reduced or not: the defined map is the same. $T$ is the set of points that annihilate simultaneously all $h_i$, thus it is defined by the equations $h_0(x)=\dots=h_n(x)=0$. On the contrary, the image $W$ is defined by the algebraic relations among $h_0,\dots,h_n$.
    
    \pause
    With Thms. \ref{binary} and \ref{ternary} we begin to understand that the classification is really linked to the nature of the reduced self-vanishing systems. In particular, when it turns out to be a constant vector, then a variable can be eliminated. This inspires the following lemmas concerning the behaviour of $\Phi$ depending on whether $\mathbf{h}$ is constant or not.
    
    \begin{lemma}\label{constanth}
        Fixing the above notations, then the following are equivalent:
        \begin{itemize}
            \item $\mathbf{h}$ is constant,
            \item $\dim W=0$, or equivalently the image is a single point,
            \item $T=\emptyset$, i.e. $\Phi$ is defined everwhere in $\PP^n$.
        \end{itemize}
    \end{lemma}
    \begin{proof}
        See \cite[Proposition 6.1]{watanabe2019theory}.
    \end{proof}
    \begin{lemma}\label{notconstanth}
        If $\mathbf{h}$ is not constant, then $\dim W\ge 1$ and $$1\le \dim T\le n-2$$
    \end{lemma}
    \begin{proof}
        See \cite[Proposition 6.2]{watanabe2019theory}.
    \end{proof}
    
    \begin{theorem}[quaternary forms]\label{quaternary}
        Let $f\in K[x_0,x_1,x_2,x_3]$ be a homogeneous polynomial with zero Hessian. Then, up to a linear change of variables, $f$ involves at most three variables. 
    \end{theorem}
    \begin{proof}
        As in the proof of Theorem \ref{ternary}, let $\mathbf{h}=(h_0,h_1,h_2,h_3)$ be a reduced self-vanishing system constructed from $f$. With the notations of Definition \ref{mapZ}, we can consider $L(W)$ the linear closure of $W$, namely $$L(W)=\bigcap\{L:\,L \text{ is a linear space and }W\subseteq L\}.$$
        It is possible to prove that $L(W)\subseteq T$ (\cite[Theorem 6.7]{watanabe2019theory}). So, if we set $s=\dim L(W)$, then by Lemma \ref{notconstanth} $s\le 1$.
        \begin{itemize}
            \item[$s=0$:] $W\subset L(W)$ is zero dimensional, so by Lemma \ref{constanth} the theorem is true.
            \item[$s=1$:] knowing that $\dim_K (Kh_0\oplus\dots\oplus Kh_4)=s+1=2$ and by Lemma \ref{elimintionsuperpower}, we can suppose that the form $g\in\I(f)$, which defines $\mathbf{h}$, does not depend on $x_0$ and $x_1$. So $g$ results to be a form in two variables and thus linear. By Proposition \ref{basicelimination} the result is again true.
            \qedhere
        \end{itemize}
    \end{proof}
    
    \begin{theorem}[quinary forms]\label{quinary}
        Let $K[x_0,x_1,x_2,x_3,x_4]$ be the polynomial ring in five variables and let $\Delta$ be a homogeneous polynomial of the form
        $$\Delta=p_0(x_3,x_4)x_0+p_1(x_3,x_4)x_1+p_2(x_3,x_4)x_2.$$
        Then every element of the ring extension $K[x_3,x_4][\Delta]$ has vanishing Hessian.
        
        Conversely, let $f$ be a form with zero Hessian which does not define a cone, then there exists a homogeneous polynomial $\Delta=p_0(x_3,x_4)x_0+p_1(x_3,x_4)x_1+p_2(x_3,x_4)x_2$ such that $f\in K[x_3,x_4][\Delta]$ after a suitable change of variables.
    \end{theorem}
    \begin{proof}
        Let $f\in K[x_3,x_4][\Delta]$ and let $\mathbf{f}=(f_0,\dots,f_4)$ the system of its partial derivatives. Using the chain rule we can write for $i=0,1,2$
        $$f_i=\frac{\partial f}{\partial x_i}=\frac{\partial f}{\partial \Delta}p_i.$$
        So we obtain the following series of equalities
        $$\frac{f_i}{f_2}=\frac{p_i(x_0,x_1)}{p_2(x_0,x_1)}=\frac{p_i(x_0/x_1,1)}{p_2(x_0/x_1,1)}\in K(x_0/x_1),\quad i=0,1.$$
        We have just proved that $f_0,f_1,f_2\in K(x_0/x_1,f_2)$ which has transcendence degree at most equal to $2$. So, they must be algebraically dependent and, by Proposition \ref{partialalgebraic}, $f$ has zero Hessian.
        
        On the contrary, let $f$ be a form with vanishing Hessian which involves properly all the variables. It is possible to find a reduced self-vanishing system $\mathbf{h}=(h_0,\dots,h_4)$, as defined in Theorem \ref{definizionedih}, such that $h_3=h_4=0$ and $h_0,h_1,h_2$ involve only the variables $x_3,x_4$, i.e. $\mathbf{h}$ is GN type. A proof of this fact can be found in \cite[Proposition 7.2]{watanabe2019theory}. Let $\Tilde{K}$ be the algebraic closure of $K(x_3,x_4)$. By \cite[Theorem 6.4]{watanabe2019theory}, $f$ is a polynomial over $\Tilde{K}$ in
        \begin{equation*}
            A:=\begin{vmatrix}
                \frac{\partial h_0}{\partial x_3} & \frac{\partial h_1}{\partial x_3} & \frac{\partial h_2}{\partial x_3}\\
                \frac{\partial h_0}{\partial x_4} & \frac{\partial h_1}{\partial x_4} & \frac{\partial h_2}{\partial x_4}\\
                x_0 & x_1 & x_2
            \end{vmatrix}.
        \end{equation*}
        It could happen that $A$ is reducible, but in this case only one of its irreducible factors is not in $K[x_3,x_4]$. We name this factor $\Delta$ and it has the form $\Delta=p_0(x_3,x_4)x_0+p_1(x_3,x_4)x_1+p_2(x_3,x_4)x_2$.
        
        As we said, $f=b_0+b_1\Delta+b_2\Delta^2+\dots$ for some elements $b_i\in \Tilde{K}$. So the coefficient of $x_j^i$ in $f\in K[x_3,x_4][x_0,x_1,x_2]$ is $b_ip_j^i\in K[x_3,x_4]$ for $j=0,1,2$. Then we must obtain $b_i\in K(x_3,x_4)$. Since the elements $\{p_0,p_1,p_2\}$ are coprime, from the choice of $\Delta$, we deduce that the denominator of every $b_i$ cannot be cancelled by a factor of $\Delta$, so $b_i\in K[x_3,x_4]$. Finally, we obtain that $f$ is an element of $K[x_3,x_4][\Delta]$.
    \end{proof}
    \begin{remark}\label{importantremark}
        We observe that the polynomials $p_0,p_1,p_2$ must be linearly independent, but algebraically dependent. In fact, if $p_2=ap_0+bp_1$, then we can choose new variables as
        $$y_0=x_0+ax_2,\quad y_1=x_1+bx_2,\quad y_2=x_2$$
        and now $f$ does not depend on $y_2$. In the proof of Theorem \ref{quinary} we have also proved that $p_0,p_1,p_2$ are algebraically dependent. More in general, if $p_0,\dots,p_s$ are $s+1$ homogeneous polynomials in $s$ variables of the same degree, then they are algebraically dependent
    \end{remark}
    
    \subsection{Geometrical approach to cubic forms}\label{gaussmapp}
	
    We now study the case of cubics: homogeneous polynomials of degree 3. U. Perazzo studied this case in \cite{perazzo}, then recently R. Gondim and F. Russo in \cite{gondimrusso} gave a modern treatment of his results. We now report this classification for completeness. Their approach is more geometrical and involves the dual variety and the Gauss map of a hypersurface.
    
    \pause
    Given a hypersurface $X=V(f)\subseteq\PP^n$ we can define the rational map
    $$\nabla_f:\PP^n\dashrightarrow(\PP^{n})^*$$
    defined by
    $\nabla_f([x_0:\dots:x_n])=\left[\frac{\partial f}{\partial x_0}(x):\dots:\frac{\partial f}{\partial x_n}(x)\right]$. This map is the so called \emph{polar map}. We name $Z=\overline{\nabla_f(\PP^n)}$ the polar image of the polar map. If we restrict this map to $X$ we obtain the \emph{Gauss map}
    $$\G_X=\nabla_f|_X:X\dashrightarrow(\PP^{n})^*.$$
    This map is defined on the smooth locus of $X$ and it associates, to every smooth point, its embedded tangent space of $X$. The image $X^*=\overline{\G_X(X)}$ is called the \emph{dual variety} of $X$. Clearly, it holds $X^*\subseteq Z$. It is also possible to define the dual variety of $Z$ which we can think to be contained in $\PP^n$, after the identification $(\PP^{n})^{**}\leftrightarrow\PP^n$. A more complete study of such varieties can be found in Chapter \ref{chapter4}. 
    
    \pause
    Let $f$ be a form with vanishing Hessian and let $\mathbf{h}$ be its self-vanishing system defined by an algebraic relation $g$ (Proposition \ref{definizionedih}). It is noticeable that the map $\Phi$, as defined in Def. \ref{mapZ}, is simply the composition $\Phi=\nabla_g\circ\nabla_f$.
    
    \pause
    One can ask why we have introduced the dual variety. This can be explained by the following classical argument.
    \begin{proposition}
        If $X=V(f)\subset \PP^n$ is a hypersurface, then the following are equivalent:
        \begin{itemize}
            \item $X$ is a cone;
            \item $X^*$ is degenerate, i.e. $L(X^*)\subsetneqq (\PP^{n})^*.$
        \end{itemize}
    \end{proposition}
    \begin{proof}
        $X$ is a cone is equivalent to say that $X$ is made of lines all passing through a fixed point $Q\in X$. The point $Q$ corresponds to a hyperplane $H\subseteq(\PP^{n})^*$ and then $X^*\subseteq H$ and so $X^*$ is degenerate. Conversely, if for a suitable hyperplane $H$ it holds $X^*\subseteq L(X)\subseteq H$ then $X$ is a cone with vertex $Q=\G_X^{-1}(H)$.
    \end{proof}
    
    \begin{theorem}[cubic hypersurfaces in $\PP^4$]\label{cubicsinp4}
        Let $X=V(f)\subseteq\PP^4$ be an irreducible cubic hypersurface with zero Hessian, not a cone. Then $X$ is projectively equivalent to $V(x_3^2x_0+x_3x_4x_1+x_4^2x_2)$.
    \end{theorem}
    \begin{proof}
        Having $f$ vanishing Hessian, and by Theorem \ref{quinary}, there exists a polynomial $\Delta=p_0(x_3,x_4)x_0+p_1(x_3,x_4)x_1+p_2(x_3,x_4)x_2$ such that $X$ is projectively equivalent to $V(g)$ with $g\in K[x_3,x_4][\Delta]$. Now we use the properties of $f$ to determine the polynomials $g$ and $\Delta$: $g$ is of the form $$g=g_0(x_3,x_4)+\Delta g_1(x_3,x_4)+\Delta^2g_2(x_3,x_4)+\dots,$$
        therefore, the degree of $\Delta$ can vary between $0$ and $3$. We analyse the various possibilities:
        
        \begin{itemize}
            \item[$\deg \Delta=0$:] $g$ does not depend on the first three variables, so $f$ must be a cone;
            \item[$\deg \Delta=1$:] with a linear change of variables, we can eliminate a variable from $g$ implying $f$ to be a cone;
            \item[$\deg \Delta=2$:] similarly to the previous case we can generate a contradiction;
            \item[$\deg \Delta=3$:] the polynomials $p_i$ have degree 2 and they form a basis of $K[x_3,x_4]_2$ as $K$-vector space (because of Remark \ref{importantremark}). So, we can choose them to be $p_0=x_3^2, p_1=x_3x_4, p_2=x_4^2$. We are almost done because now $g=g_0(x_3,x_4)+\lambda(x_0x_3^2+x_1x_3x_4+x_2x_4^2)$. The polynomial $g_0$, being homogeneous of degree $3$, must be of the form $g_0=x_3^2(\alpha x_3+\beta x_4)+x_4^2(\gamma x_3+\epsilon x_4)$. The final linear transformation is:
            \begin{equation*}
                \begin{cases}
                    x_0'=x_0+\alpha x_3+\beta x_4,\\
                    x_1'=x_1,\\
                    x_2'=x_2+ \gamma x_3+\epsilon x_4,\\
                    x_3'=\xi x_3,\\
                    x_4'=\xi x_4;
                \end{cases}
            \end{equation*}
            where $\xi$ is a root of the polynomial $q(t)=t^2-\lambda$. Finally, $g$ is in the desired form.
        \end{itemize}
        
        A first proof of this theorem can be attributed to U. Perazzo in \cite{perazzo}. A modern version of his proof can be found also in \cite[Theorem 7.6.7]{russo}.
    \end{proof}
    
    \begin{definition}
        Let $X=V(f)\subseteq\PP^n$ be a hypersurface with vanishing Hessian, not a cone. Let $\nabla_f:\PP^n\dashrightarrow(\PP^{n})^*$ be its polar map and $Z$ be its polar image. Then we can define the \emph{Perazzo map} by
        \begin{align*}
            \mathbf{P}_X:\,&\PP^n\dashrightarrow\Gr(\text{codim}Z-1,n)\\
            &P\,\longmapsto\,\,(T_{\nabla_f(P)}Z)^*
        \end{align*}
        where $\Gr(a,b)$ denotes the set of the linear subspaces of dimension $a$ in the vector space $K^b$.
        We also define the \emph{Perazzo rank} to be $\mu=\dim(\overline{\mathbf{P}_X(\PP^n}))$. The general fiber of this map is a linear space of dimension $n-\mu$. So, knowing that the image has dimension $\mu$, the fibers over the general points form a congruence of linear spaces.
        
        If $X$ is an irreducible cubic hypersurface such that every element of the congruence passes through a fixed $\PP^{n-\mu-1}$, then $X$ is called \emph{Special Perazzo Cubic hypersurface}.
    \end{definition}
    \begin{theorem}[cubic hypersurfaces in $\PP^5$ and $\PP^6$]
        Let $n=5,6$ and let $X=V(f)\subseteq\PP^n$ be an irreducible cubic hypersurface with zero Hessian, not a cone. Then $X$ is a Special Perazzo Cubic hypersurface such that $Z^*$ is either a cone or, if $n=6$, a surface in $\PP^3$ of degree at most $3$.
    \end{theorem}
    \begin{proof}
        See \cite[Section 4]{gondimrusso}.
    \end{proof}
    
	\section{Historical notes}
	
	In this section we give a brief historical overview about the classification of forms with vanishing hessian. The first mathematician who was interested in this argument was O. Hesse. In 1852 and 1859 he published two articles \cite{hesse1} and \cite{hesse2} in which he stated and gave twice a proof that all possible forms with zero Hessian are cones. Subsequently in \cite{pasch}, a paper of 1875, M. Pasch proved Hesse's claim for binary and ternary forms: this seemed to reinforce Hesse theory. Unfortunately, the proofs of Hesse were incorrect and in 1876 P. Gordan and M. Noether wrote \cite{gordannoether}. They gave a counterexample to Hesse's results, and a classification of such forms up to $5$ variables, as we saw in Theorems \ref{binary}, \ref{ternary}, \ref{quaternary}, and \ref{quinary}. Nevertheless, they also proved that up to birational transformations forms with vanishing Hessian are indeed cones.
	
	\pause
	Being written in German and due to the involved arguments it contains, the work of Gordan and Noether remained almost unknown. In 1900, U. Perazzo studied Gordan-Noether's paper and in \cite{perazzo} he described cubic forms with vanishing Hessian in dimension $n=4,5,6$. His paper was the starting point of a series of efforts with the purpose to obtain a geometric interpretation of hypersurfaces with vanishing hessian. In particular, in the middle of 900s, A. Franchetta studied the case of hypersurfaces in $\PP^4$ and published the paper \cite{franchetta}. Some years later in 1957 and in 1964, R. Permutti studied these two papers and rewrote in \cite{permutti1} and in \cite{permutti3} the classic argument of Gordan and Noether trying to generalise it. In the same years, he also studied the case over a field of general characteristic in \cite{permutti2}. It is remarkable that all these papers are written in Italian and for this reason they were not known outside Italy.
	
	\pause
	In 2009, F. Repetto, and A. Garbagnati wrote \cite{garbagniatirepetto}, the first paper in English in which they formalised some of the theory that was constructed in the previous century. Another important step, was made by R. Gondim and F. Russo in \cite{gondimrusso} where all the remaining theory, not treated by Repetto and Garbagnati, was studied and formalised. The geometrical properties of hypersurfaces with zero Hessian were studied also in \cite{gondimrussostagliano} (written by Gondim, Russo, and G. Staglianò), and in \cite{cilibertorussosimis} (written by C. Ciliberto, Russo, and A. Simis).
	
	\pause
	A different approach was taken in 1990 by H. Yamada who was interested in the relation between the Strong Lefschetz Property and the Hessian matrix. Unfortunately, his notes \cite{yamada} were never completed and remained unpublished. For this reason, one of his student J. Watanabe (first alone and then in collaboration with M. de Bondt) completed Yamada's work. He published a first paper \cite{watanabeerror} in 2014, but it still contains some mistakes, then in 2020 Watanabe and de Bondt wrote \cite{watanabe2019theory} the final and right version. The same kind of study, but independent, was also done by C. Lossen in \cite{lossen}.
	
	\pause
	The paper of Gordan and Noether is still alive nowadays and these arguments are object of new efforts. Different proofs of Gordan and Noether's theorems are provided. The most noticeable ones are:
	\begin{itemize}
	    \item \cite{pirolafavalebricalli} written at the beginning of 2022 by G. P. Pirola, D. Bricalli and F. Favale;
	    \item \cite{debondt1}, \cite{debondt2} where M. de Bondt studied the so called "quasi - translations".
	\end{itemize}

 	
	\chapter{Gorenstein algebras and Lefschetz Properties}
	In this chapter we will introduce Gorenstein algebras and the Lefschetz properties. Our principal goal is to introduce various tools, fundamental in checking whenever a Gorenstein algebra satisfies or not one or both the Strong and Weak Lefschetz properties. In fact, these properties can be defined generally for Artinian graded algebras, but the Gorenstein property guarantees a more easier environment to work with.
	
	\pause
	The notion of Gorenstein algebra, and more in general of Gorenstein ring, has been studied since D. Gorenstein in \cite{gorenstein}. We will follow a more suitable formulation for our treatment as in the book \cite{sixauthors} written by T. Harima, T. Maeno, H. Morita, Y. Numata, A. Wachi, and J. Watanabe. We pay our attention to the following scenario.
	
	\begin{definition}\label{gradedArtiniangorenstein}
	    Let $A$ be a $K$-algebra where $K$ is an algebraically closed field of characteristic zero. We say that $A$ is a \emph{graded Artinian $K$-algebra} if it can written as
	    $$A=\bigoplus_{i=0}^d A_i$$
	    where $\{A_i\}_i$ is a grading of $A$ such that $A_0\cong K$. $d$ is called the \emph{socle degree} of $A$ and $A_d$ is the \emph{socle} of $A$. We also define its \emph{codimension} as $\dim_K A_1$ which is always finite. In our discussion we consider only algebras with \emph{standard} grading i.e. $A=K[A_1]$.
	\end{definition}
	
	Definition \ref{gradedArtiniangorenstein} contains various hypotheses, and each of them is necessary. We need $A$ to be graded because it is necessary to define the Lefschetz properties. $A$ is Artinian, so that there are finitely many non-zero homogeneous components, all of finite dimension as $K$-vector fields. Furthermore, with these assumptions, $A$ is isomorphic to the quotient of a polynomial ring with respect to a homogeneous ideal.
	
	\section{Gorenstein algebras}
	We now introduce the concept of Gorenstein graded $K$-algebra and its relation with the Poincaré duality algebra. We then characterise such algebras using the Annihilator of a form. From now on $A$ is a graded Artinian $K$-algebra.
	\begin{definition}\label{definizionedigorenstein}
	    Let $A=\bigoplus_{i=0}^d A_i$ be a graded Artinian $K$-algebra. Let $\mathfrak{m}:=\bigoplus_{i=1}^d A_i$ be the maximal homogeneous ideal. $A$ is said to be \emph{Gorenstein} if the socle ideal
	    $$\Soc(A):=(0:\mathfrak{m})=\{x\in A:\,x\,\mathfrak{m}=0\}$$
	    of $A$ is a $1$-dim vector space over $K$.
	\end{definition}
	
    For a general such an algebra, two different concepts seem to coincide: on one hand we have the socle $A_d$ and in the other the socle ideal $\Soc(A)$. We have that $A_d\subseteq\Soc(A)$, but in general they can be different\footnote{A graded algebra of socle degree $d$ is named \emph{level} if $(0:\mathfrak{m})=A_d$}.
    \begin{example}
        Let $I=(x^2,xy,y^3)$ be an ideal of $R=K[x,y]$. Then the quotient $A=R/I$ is a graded Artinian $K$-algebra, in particular $A=K\oplus\langle x,y\rangle\oplus\langle y^2\rangle$. In this case, being $2$ the socle degree, $A_2$ is not the socle ideal of $A$ because it contains also the element $x$. In particular, $A$ is not Gorenstein.
    \end{example}
    In a Gorenstein algebra they coincide because the socle ideal is one dimensional.
    
    \begin{definition}
        Let $A=\bigoplus_{i=0}^d A_i$ as in Definition \ref{gradedArtiniangorenstein}. $A$ is called a \emph{Poincaré duality} algebra if $\dim_K A_d=1$, and the bilinear pairing
        \begin{equation}\label{pairing}
            A_k\times A_{d-k}\longrightarrow A_d\cong K
        \end{equation}
        is non-degenerate for $k=0,1,\dots,\left[\frac{d}{2}\right]$.
    \end{definition}
    
    \begin{proposition}
        A graded Artinian $K$-algebra is a Poincaré duality algebra if and only if $A$ is Gorenstein.
    \end{proposition} 
    \begin{proof}
        See \cite{geramitaharimamiglioreshin}, or \cite[Proposition 2.1]{maenowatanabe}.
    \end{proof}
        
    \begin{definition}
        Let $S=K[X_0,\dots,X_n]$ be the ring of differential operators which acts on $K[x_0,\dots,x_n]$ (by definition $X_i:=\sfrac{\partial}{\partial x_i}$). Given a form $f\in K[x_0,\dots,x_n]$ of degree $d$ we define the \emph{Annihilator} of $f$ in $S$ as
        $$\Ann_S(f):=\{q\in S:\,q\cdot f=0\}.$$
        $\Ann_S(f)$ is a homogeneous ideal of $S$ such that $S_i\subseteq\Ann_S(f)$ for $i>d$. Clearly, $A_d$ is isomorphic to $K$ under the map
        \begin{align*}
            A_d&\xrightarrow{\,\,\cong\,\,} K\\
            [g]\,&\longmapsto\,\,g\cdot f
        \end{align*}
        and the map (\ref{pairing}) is non-degenerate for every $k$. Thus, the quotient $A=S/\Ann_S(f)$ is a graded Artinian Gorenstein $K$-algebra of socle degree $d$. We name $A$ the Gorenstein $K$-algebra associated to $f$.
    \end{definition}
    
    \begin{theorem}[Double Annihilator theorem of Macaulay]\label{gorensteinannullatore}
        There is a one to one correspondence as follows:
        \begin{equation*}
            \left\{ \parbox{4.5cm}{\centering 
                   Artinian Gorenstein algebras of socle degree $d$ and codimension $n+1$}
            \right\}
             \xrightarrowdbl{1:1}
             \left\{
                \parbox{4.5cm}{\centering 
                    hypersurfaces, not cones, of degree $d$ in $\PP^n$}
             \right\}
        \end{equation*}
        So that, given an Artinian Gorenstein algebra $A$, there exists a form $F$ such that $A\cong S/\Ann_S(F)$. $F$ is called the \emph{dual generator} of $A$.
        In particular, two homogeneous polynomials $F$ and $G$, which differ by a linear change of variables, give two isomorphic Gorenstein algebras $S/\Ann_S(F)\cong S/\Ann_S(G)$, and viceversa.
    \end{theorem}
    \begin{proof}
        There are many references, for instance \cite[Chapter IV]{macaulay}, \cite[Theorem 21.6]{eisenbud}, and \cite[Lemma 2.12]{iarrobinokanev}.
    \end{proof}

	\section{Lefschetz properties and Hilbert vector}
	
	\begin{definition}[Hilbert vector]\label{hilbertfunction}
	    Let $A$ be as in Definition \ref{gradedArtiniangorenstein}. We name \emph{Hilbert function} the map $H_A:\numberset{N}\to \numberset{N}$ given by
	    $$H_A(i)=\dim_K A_i.$$
	    Being this function definitely zero, we can embrace its information in a vector $h=(h_0,h_1,\dots,h_d)$ where $d$ is the socle of $A$ and $h_i=H_A(i).$ This vector is called the \emph{Hilbert vector} of $A$ ($h$-vector for shortly).
	    
	    The Hilbert vector $h$ is said to be \emph{unimodal} if, for some integer $s$,
	    $$h_0\le h_1\le\dots\le h_s\ge h_{s+1}\ge\dots\ge h_d.$$
	    Analogously, it is called \emph{symmetric} if
	    $$h_k=h_{d-k} \text{ for every }k.$$
	\end{definition}
	\begin{remark}\label{gorenstainsymmetric}
	   The Hilbert vector of an Artinian Gorenstein algebra is, by definition, always symmetric, but in general it can also be non-unimodal.
	\end{remark}
	
	\begin{example}\label{stanley}
	    One of the first counterexample to the unimodality of the Gorenstein algebras was given by R. P. Stanley in \cite[Example 4.3]{stanley}. He gave an explicit Gorenstein algebra with $h$-vector $(1,\ 13, \ 12,\ 13,\ 1)$. We start defining the quotient $T=K[x,y,z]/(x,y,z)^4$. This is a graded Artinian $K$-algebra with socle degree $t=3$, but not Gorenstein; in fact, its Hilbert vector is $(1,\ 3,\ 6,\ 10)$. Define $E=\hom_K(T,K)$, we can consider the product $A=T\times E$, endowed with the componentwise addiction and the multiplication $(a,\phi)\cdot(b,\psi)=(ab,a\psi+b\phi)$. With this structure, $A$ can be graded becoming an Artinian $K$-algebra which is also Gorenstein. The Hilbert function of $A$ satisfy the relation
	    \begin{equation*}
	        \begin{cases}
	        H_A(i)=H_T(i)+H_T(t+1-i) &\text{ if } 0\le i \le t+1\\
	        H_A(i)=0 &\text{ otherwise.}
	        \end{cases}
	    \end{equation*}
	    Therefore, $A$ has Hilbert vector $(1,\ 13, \ 12,\ 13,\ 1)$ as we have stated.
	    
	    \pause
	    Similarly, we can consider the graded Artinian $K$-algebra  $$T=K[x,y,z]/(x,y,z)^m$$ with socle degree $t=m-1$, $m\ge 4$, and proceed as above. In this way we obtain infinitely many Gorenstein algebras with non-unimodal Hilbert vector.
	    
	    \pause
	    The mathematical field of non-unimodal Gorenstein algebras is still a research topic. For example see
	    \cite{migliorezanello}.
	\end{example}
	
	At first glance, it seems that, for a general graded Artinian $K$-algebra $A$, the behaviour of the Hilbert vector of $A$ could vary freely. This fact turns out to be false: Macaulay gave an upper bound to the growth of $\dim A_k$ as the parameter $k$ varies. He used the so call \emph{s-th expansion}: let $m$ and $s$ two positive integers, then $m$ can be uniquely written as 
	$$m=\binom{m_s}{s}+\binom{m_{s-1}}{s-1}+\dots+\binom{m_i}{i},$$
	where $m_s>m_{s-1}>\dots>m_i\ge i\ge 1$. This expression is called the s-th expansion of $m$.
	
	\pause
	We also define the quantity
	$$m^{\langle s\rangle}=\binom{m_s+1}{s+1}+\binom{m_{s-1}+1}{s}+\dots+\binom{m_i+1}{i+1}.$$
	
	\begin{example}
	    We give some examples of binomial expansions in the following table.
	    \begin{center}
	        \def\arraystretch{2.8}
	        \hspace*{-0.4cm}
	        \begin{tabular}{|c||c|c|c|c|}
	            \hline
	           $(m, s)$ & $(5,1)$ & $(6,2)$& $(7,2)$ & $(6,3)$\\[0.8em]
	           \hline
	            \stackanchor[8pt]{s-th}{expansion} & $\displaystyle\binom{5}{1}$ & $\displaystyle\binom{4}{2}$ & $\displaystyle\binom{4}{2}+\binom{1}{1}$ & $\displaystyle\binom{4}{3}+\binom{2}{2}+\binom{1}{1}$\\[0.8em]
	            \hline
	            $m^{\langle s\rangle}$ & $\displaystyle\binom{6}{2}= 15$ & $\displaystyle\binom{5}{3}=10$ & $\displaystyle\binom{5}{3}+\binom{2}{2}=11$ & $\displaystyle\binom{5}{4}+\binom{3}{3}+\binom{2}{2}=7$\\[0.8em]
	            \hline
	        \end{tabular}
	    \end{center}
	\end{example}
	
	\pause
	\begin{theorem}\label{osequenza}
	    Let $h=(h_0,h_1,\dots,h_d,\dots)$ be an infinite sequence of non-negative integers. Then there exists a graded $K$-algebra with $h$ as Hilbert vector if and only if $h$ is an O-sequence, i.e.
	    \begin{itemize}
	        \item $h_0=1$,
	        \item $h_{i+1}\le h_i^{\langle i \rangle}$.
	    \end{itemize}
	\end{theorem}
	\begin{proof}   
	    See \cite{macaulay2}.
	\end{proof}
	
	\begin{definition}[Weak Lefschetz property]\label{wlpdefinition}
	    Let $A=\bigoplus_{i=0}^d A_i$ be an Artinian graded $K$-algebra. Let $L\in A_1$ be an element of degree $1$. We say that $L$ is a \emph{Weak Lefschetz element} if the multiplication map
	    \begin{equation*}
	        \times L:A_i\to A_{i+1}
	    \end{equation*}
	    has full rank for all $i$, $0\le i\le d-1$. We say that $A$ has the Weak Lefschetz property (WLP for short) if it admits a Weak Lefschetz element.
	\end{definition}
	
	\begin{definition}[Strong Lefschetz property]\label{slp}
	    Let $A=\bigoplus_{i=0}^d A_i$ be an Artinian graded $K$-algebra. Let $L\in A_1$ be an element of degree $1$. We say that $L$ is a \emph{Strong Lefschetz element} if the multiplication map
	    \begin{equation*}
	        \times L^c:A_i\to A_{i+c}
	    \end{equation*}
	    has full rank for all $i$, $0\le i\le d-1$, and for all $c$, $1\le c\le d-i$. We say that $A$ has the Strong Lefschetz property (SLP for short) if it admits a Strong Lefschetz element.
	\end{definition}
	
	\begin{definition}[Strong Lefschetz property in the narrow sense]\label{slpins}
	    Let $A=\bigoplus_{i=0}^d A_i$ be an Artinian graded $K$-algebra. We say that $A$ has the \emph{Strong Lefschetz property in the narrow sense} if there exists an element $L\in A_1$ such that the multiplication map
	    \begin{equation*}
	        \times L^{d-2i}:A_i\to A_{d-i}
	    \end{equation*}
	    is an isomorphism for all $i$, $0\le i\le \lfloor \frac{d}{2}\rfloor$.
	\end{definition}
	
	Given $A$, as we will see in Proposition \ref{unimodalsymmetric}, the SLP in the narrow sense implies the symmetry of $A$'s Hilbert vector. If we know in advice that the $h$-vector is symmetric, then the notion of Strong Lefschetz property coincide with the one in the narrow sense. In particular,
	\begin{center}
	    SLP+symmetric Hilbert vector $\iff$ SLP in the narrow sense.
	\end{center}
	For this reason when we will consider Gorenstein algebras, we will use both notions without distinction.
	
	\begin{proposition}\label{unimodalsymmetric}
	    Let $A$ be an Artinian graded algebra. If $A$ has one among the Lefschetz properties, then its Hilbert vector is unimodal. Furthermore, if $A$ has the Strong Lefschetz property in the narrow sense, then its $h$-vector is also symmetric.
	\end{proposition}
	\begin{proof}
	    By definition, if $A$ has the SLP in the narrow sense, then its $h$-vector is symmetric. Clearly, we have a series of implication between the Lefschetz properties:
	\begin{center}
	    SLP in the narrow sense $\implies$ SLP $\implies$ WLP.
	\end{center}
	    Thus, we just need to check the unimodality property only in the WLP's case.
	    
	    Let $\mathfrak{m}$ be the maximal homogeneous ideal of $A$. Since $A$ has the standard grading, we have that, for every $k$, $\mathfrak{m}^k=\bigoplus_{i=k}^d A_i$, i.e. it is generated by $A_k$. Let $s\ge 0$ be the first integer such that $\dim_K A_s>\dim_K A_{s+1}$. Since $A$ has the WLP, there exists an element $L\in A_1$ such that $\times L:A_s\to A_{s+1}$ is surjective. Being $\times L:A\to A$ a morphism of graded $K$-algebras, we have that $\mathfrak{m}^{s+1}=L\mathfrak{m}^s$. Furthermore, this fact implies that $\mathfrak{m}^{j+1}=L\mathfrak{m}^j$ for every $j\ge s$, or equivalently that $\times L:A_j\to A_{j+1}$ is surjective for every $j\ge s$. Finally,
	    $$\dim_K A_0\le \dim_K A_1\le\dots\le \dim_K A_s\ge \dim_K A_{s+1}\ge\dots\ge \dim_K A_d.\qedhere$$
	\end{proof}
	
	As we have seen in Example \ref{stanley}, there exists Gorenstein algebras, and more in general Artinian graded algebras, with non-unimodal Hilbert vector. The non-unimodality of the Hilbert vector of an algebra $A$ implies that $A$ do not satisfy any Lefschetz property. Therefore, when we treat the Lefschetz properties, it is licit to tacitly suppose the unimodality of the Hilbert vector. 
	
	\pause
	The Weak and Strong Lefschetz properties were inspired by the work of S. Lefschetz about the cohomology ring of a complex manifold. In particular, the \emph{Hard Lefschetz theorem} states that the cohomology ring associated to any compact K\"ahler manifold has the Strong Lefschetz property (more details can be found in \cite{lefschetz}).
	
	\pause
	The set of Strong Lefschetz elements forms a Zariski open subset of $A_1$; in general, it can be also empty. If $A$ has the SLP, then this subset is dense in $A_1$, so we can say that the general element is a Strong Lefschetz element. In the same way, these arguments can also be applied to the Weak Lefschetz property.
	
	\pause
	The problem of understanding whenever an algebra has such properties is still open. In the last decades, even if mathematicians had conjectured that all general graded Artinian algebras have both the Lefschetz properties, many examples of algebras failing one or both Lefschetz properties were found. As an example we give some new results about Artinian Gorenstein algebras.
	
	\begin{theorem}
	    For each pair $(n,d)\not\in \{(3,3),(3,4)\}$ with $n,d\ge 3$ there exist Artinian Gorenstein algebras of codimension $n+1$ and socle degree $d$ failing the Strong Lefschetz property. Furthermore, we can choose arbitrarily at which level the map
	    $$\times L^{d-2k}:A_k\longrightarrow A_{d-k}$$
	    is not an isomorphism.
	\end{theorem}
	\begin{proof}
	    See \cite[Corollary 3.2]{gondim}.
	\end{proof}
	
	\begin{theorem}
	    For each pair $(n,d)\not\in \{(3,3),(3,4),(3,6),(4,4)\}$ with $n,d\ge 3$ there exist Artinian Gorenstein algebras of codimension $n+1$ and socle degree $d$ with unimodal Hilbert vector, but failing the Weak Lefschetz property.
	\end{theorem}
	\begin{proof}
	    See \cite[Theorem 3.8]{gondim}.
	\end{proof}
	
    We now group together a series of results that are useful to check if an Artinian Gorenstein algebra has or not the Weak Lefschetz property.
    
    \begin{proposition}\label{checkweak}
        Let $A=\bigoplus_{i=0}^d A_i$ be a graded Artinian Gorenstein algebra and let $L\in A_1$. We call $l_k:A_k\to A_{k+1}$ the restriction of the multiplication by $L$ at the $k$-th homogeneous part.
        Then the following statements hold true.
        \begin{enumerate}
            \item If $l_{k_0}$ is injective for some $k_0$, then $l_{k}$ is injective for every $k\le k_0$.
            \item If $l_{k_0}$ is surjective for some $k_0$, then $l_{k}$ is surjective for every $k\ge k_0$.
            \item If $\dim A_{k_0}=\dim A_{k_0+1}$ for some $k_0$, then $L$ is a Weak Lefschetz element if and only if $l_{k_0}$ is an isomorphism.
            \item If $\dim A_{k_0-1}<\dim A_{k_0}>\dim A_{k_0+1}$ for some $k_0$, then $L$ is a Weak Lefschetz element if and only if $l_{k_0-1}$ is injective.
        \end{enumerate}
    \end{proposition}
    \begin{proof}
        For statements 1.,2., and 3. see \cite[Proposition 2.1]{miglioremironagel}. Statement 4. is a direct consequence of the the first two and the following fact: if $l_{k_0-1}$ is injective then the map $l_{k_0}:A_{k_0}\to A_{k_0+1}$ is surjective. Since $l_{k_0-1}$ is injective, its dual map $l_{k_0-1}^*:A_{k_0}^*\to A_{k_0-1}^*$ is surjective. Since $A$ is a Gorenstein algebra, we can view this map as $l_{k_0-1}^*:A_{k_0}\to A_{k_0+1}$. This map coincide with the map $l_{k_0}:A_{k_0}\to A_{k_0+1}$, thus it is surjective.
    \end{proof}
    
    By Proposition \ref{checkweak} we also deduce the following property for $A$, an Artinian Gorenstein algebra: a necessary condition for $A$ to have one or both the Lefschetz properties is that 
    $$h_0<h_1<\dots<h_s=h_{s+1}=\dots=h_{\left[\frac{d}{2}\right]}=\dots=h_{d-s}>h_{d-s+1}>\dots>h_d$$
    for some integer $s$. The quantity $h_s$ is called the \emph{Sperner number} of $A$. So, if we want to check the WLP the conditions 3., and 4. of Proposition \ref{checkweak} satisfy all possible cases.
    
    As a final note, the Hilbert vector of an Artinian Gorenstein algebra, which has one or both the Lefschetz properties, has restrictions on its value. In fact, a similar statement to Theorem \ref{osequenza} can be proved for Gorenstein algebras involving the Lefschetz properties.
    
    \begin{theorem}\label{sisequenza}
        Let $h=(h_0,h_1,\dots,h_d)$ be a tuple of non-negative integers. Then there exists an Artinian Gorenstein $K$-algebra satisfying the WLP with $h$ as Hilbert vector if and only if $h$ is an SI-sequence, i.e.
	    \begin{itemize}
	        \item $h$ is symmetric,
	        \item $h$ is unimodal,
	        \item $\Delta h=(h_0,h_1-h_0,\dots,h_t-h_{t-1})$ is an O-sequence, where $t=\min\{i\vert h_i\ge h_{i+1}\}$ (the $t$ where the elements of $\Delta h$ are strictly positive).
	    \end{itemize}
	    
	    The same statement holds true also replacing SLP with WLP.
    \end{theorem}
    \begin{proof}
        See \cite[Theorem 1.2]{harima} and \cite[Theorem 3.2]{altafi}.
    \end{proof}
	

	\section{Higher Hessians}\label{higherhessians}
	This property was first studied by J. Watanabe in \cite{watanabequeen}, and then, together with T. Maeno, in \cite{maenowatanabe}.
	
	\begin{definition}[higher Hessians]\label{higherhessiansdefinizione}
	    Let $G\in K[x_0,\dots,x_n]$ be a form of degree $d$, and let $A=S/\Ann_S(G)$ be the Gorenstein $K$-algebra associated to $G$ with Hilbert vector $h=(h_0,\dots,h_d)$. Given $k\le \left[\frac{d}{2}\right]$, and $\mathcal{B}=\{\alpha_i\}_i$ a $K$-basis of $A_k$ we define the \emph{k-th Hessian matrix} with respect to $\mathcal{B}$ to be the square matrix
	    $$\Hess^k_G(x_0,\dots,x_n):=\left(\alpha_i\alpha_jG(x_0,\dots,x_n)\right)_{i,j}\in M(h_k\times h_k;K[x_0,\dots,x_n]_{d-2k}).$$
	    Moreover, we define the \emph{k-th Hessian} with respect to $\mathcal{B}$
	    $$\hess^k_G(x_0,\dots,x_n):=\det(\Hess^k_G(x_0,\dots,x_n))\in K[x_0,\dots,x_n]_{(d-2k)h_k}.$$
	\end{definition}
	
	The $0$-th Hessian is just the polynomial $G$ and, in the case $\dim A_1=n+1$, the 1st Hessian, with respect to the standard basis, is the classical Hessian. The definition of $k$-th Hessian and $k$-th Hessian matrix depends on the choice of a basis of $A_k$, but the vanishing of the $k$-th Hessian is independent of this choice. Thus, the condition of vanishing of the $k$-th Hessian is well-posed.
	
	\pause
	The higher Hessians are useful to fully characterise the Strong Lefschetz elements, and, in general, to understand if an Artinian Gorenstein algebra has or not the Strong Lefschetz property.

	\begin{theorem}[Watanabe] \label{watanabe}
        Let $F\in K[x_0, \cdots, x_n]_d$ be a homogeneous polynomial of degree $d$ and let $ A= S/\Ann_S(F)$ be the associated Artinian Gorenstein algebra. Let also $S=K[X_0,\dots,X_n]$ be the ring of the differential operators.
        
        $L=a_0X_0+\cdots +a_nX_n\in A_1$ is a strong Lefschetz element of $A$ if and only if $\hess _F^k(a_0,\cdots, a_n)\ne 0$ for $k=0,\cdots,[d/2]$. More precisely, fixed a basis $\mathcal{B}$ of $A_k$, since $A$ is Gorenstein, we can canonically define a basis in $A_{d-k}$. Then, up to a multiplicative constant, $\Hess_F^k(a_0,\cdots, a_n)$ is the matrix of the dual map of the multiplication map $\times L^{d-2k}: A_{k}  \to  A_{d-k}$.
    \end{theorem}
    \begin{proof}
        Since $A$ is a Poincaré duality we have that $A_d\cong K$ through the map that associates $[q(X)]$ to $q(X)F(x)$. Moreover, an element $L=a_0 X_0+\cdots +a_n X_n\in A_1$ is a Strong Lefschetz element if and only if the bilinear pairing
        \begin{equation*}
            \begin{split}
            A_k\times A_k&\longrightarrow A_d\cong K\\
            ([\xi],[\nu])\,&\longmapsto L^{d-2k}[\xi\nu]=L^{d-2k}\xi\nu(X)F(x)
        \end{split}
        \end{equation*}
        is non-degenerate for $k=0,1,\dots,[d/2]$. This property can be checked on a fixed basis $\mathcal{B}=\{\alpha_i\}_i$ of $A_k$. Recalling that the formula $(a_0X_0+\cdots +a_nX_n)^eG(x_0,\dots,x_n)=e!\,G(a_0,\dots,a_n)$ holds for every form $G$ of degree $e$,
        we obtain
        \begin{equation*}
           \begin{split}
                \left(L^{d-2k}\alpha_i\alpha_j(X)F(x)\right)_{i,j}&=(d-2k)!\left(\alpha_i\alpha_j(X)F(x)_{|(a_0,\dots,a_n)}\right)_{i,j}\\
                &=(d-2k)!\Hess^k_F(a_0,\dots,a_n).\qedhere
            \end{split} 
        \end{equation*}
    \end{proof}
	
	\begin{corollary}\label{watanabecorollario}
	    Let the notations be as in Theorem \ref{watanabe}. Then the following statements hold true.
	    \begin{enumerate}
	        \item $A$ has the SLP if and only if $\hess _F^k\ne 0$ for every $k=1,2,\dots,\left[\frac{d}{2}\right]$.
	        \item If $d\le 4$, $A$ has the SLP if and only if $F$ has not vanishing Hessian. Furthermore,
	        if $n\le 4$ and $V(F)$ is not a cone, then $A$ always has the SLP.
	        \item If $F$ has vanishing Hessian, then $A$ fails the SLP.
	    \end{enumerate}
	\end{corollary}
	\begin{proof}
	    It is a direct consequence of Theorems \ref{binary}, \ref{ternary}, \ref{quaternary}, and \ref{watanabe}.
	\end{proof}
	
	\begin{example}
	We now give an example where we apply Theorem \ref{watanabe} and Corollary \ref{watanabecorollario} to study the Lefschetz properties of an Artinian Gorenstein algebra of socle degree $d=5$. This example is known as Ikeda's example (see \cite[Example 4.4]{ikeda}): it gives a form of degree $5$ in four variables with non-vanishing Hessian, but with vanishing $2$-nd Hessian. This will imply that the related Gorenstein algebra fails to have both the Strong and the Weak Lefschetz property. 
	
	\pause
	 Let $S=K[X,Y,Z,T]$ be the ring of differential operators on $R$. We consider $f=xz^3t+yzt^3+x^3y^2\in K[x,y,z,t]$ a form of degree five and with four variables. A simple computation shows the minimal generators of $\Ann_S(f)$:
	 \begin{equation*}
	    \begin{split}
	        \Ann_S(f)=\langle  XT^2, Y^2T, XYT, X^2T, YZ^2, XZ^2-YT^2, Y^2Z ,XYZ,\\ X^2Z, Y^3, T^4, Z^2T^2,Z^4, X^2Y^2-2Z^3T,X^3Y-2ZT^3, X^4 \rangle.
	    \end{split}
	\end{equation*}
    In this way, the $h$-vector of $A=S/\Ann_S(f)$ can be easily computed as $h_A=(1, \ 4, \ 10, \ 10, \ 4, \ 1)$. Since the Hessian of $f$ is
    \begin{equation*}
	    \begin{split}
	        \hess(f)=8 (9 t^7 x y^3 z + 8 t^6 z^6 - 45 t^5 x^2 y^2 z^3 + 27 t^4 x^4 y^4 - 27 t^3 x^3 y z^5 -\\- 54 t^2 x^5 y^3 z^2 + 3 t x^4 z^7 + 27 x^6 y^2 z^4)\ne 0,
	    \end{split}
	\end{equation*}
	by Theorem \ref{watanabe}, the general element $L\in A_1$ makes the map
	$\times L^3:A_1\longrightarrow A_4$
	an isomorphism. We are now interested in the second Hessian. With this aim, we fix a $K$-basis of $A_2$: $\{X^2, Y^2, Z^2, T^2, XY, XZ, XT, YZ, YT, ZT\}$. We obtain
    
    $$
    \Hess_f^2=\begin{pmatrix}
     0 & 12x & 0 & 0 & 6y & 0 & 0 & 0 &  0 & 0  \\
     12x & 0 & 0  & 0 & 0 & 0 & 0 & 0 & 0 &  0  \\
      0 & 0 & 0 & 0 & 0 & 6t & 6z & 0 &  0 & 6x  \\
       0 & 0 & 0 & 0 & 0 & 0 & 0 &  6t & 6z & 6y \\
       6y & 0 & 0 & 0 & 6x & 0 & 0 & 0 &  0 & 0  \\
       0 & 0 & 6t & 0 & 0 & 0 & 0 &  0 &  0 & 6z  \\
       0 & 0 & 6z & 0 & 0 & 0 & 0 & 0 &  0 & 0  \\
       0 & 0 & 0 & 6t & 0 & 0 & 0 & 0 &  0 & 0  \\
        0 & 0 & 0 & 6z & 0 & 0 & 0 & 0 &  0 & 6t  \\
         0 & 0 & 6x & 6y & 0 & 6z & 0 & 0 &  6t & 0  
    \end{pmatrix}.$$
    
    A straightforward computation proves that $\hess_f^2=0$. By Theorem \ref{watanabe}, for any $L\in A_1$, the multiplication map $
    \times L: [A]_{2}  \longrightarrow  [A]_{3}
    $ has zero determinant; therefore, it is not an isomorphism.
    
    \pause
    In conclusion,
    \begin{itemize}
        \item $\hess(f)\ne 0\,\implies$ for the general element $L\in A_1$, the multiplication map
        $$\times L^3:A_1\longrightarrow A_4$$
        is an isomorphism;
        \item $\hess_f^2=0\,\implies$ for all elements $L\in A_1$, the map
        $$\times L:A_2\longrightarrow A_3$$
        is not an isomorphism. 
    \end{itemize}
    Therefore, the Gorenstein algebra $A$ fails to have both the Strong and the Weak Lefschetz properties. Thus, the second point of Corollary \ref{watanabecorollario} gives an optimal bound to the link between the Strong lefschetz Property and the Hessian of $F$.
    \end{example} 
	
	
	\chapter{Perazzo hypersurfaces in $\PP^4$}
    In this chapter we will see the application of our previous theory to a class of hypersurfaces, the so called "Perazzo hypersurfaces", that are a special class of hypersurfaces with vanishing Hessian, but in general not cones. Such hypersurfaces were first studied, in the case of cubic hypersurfaces in low dimension, by U. Perazzo \cite{perazzo}. He was trying to determine their geometry using the relation between the dual variety and the vanishing Hessian condition. We give a general definition of Perazzo hypersurface, and we study, in particular, the Gorenstein $K$-algebras associated to Perazzo 3-folds, trying to understand whenever they have or not the Lefschetz properties. With this goal in mind, we look for the minimum and maximum $h$-vector that these algebras can reach. We get the answer using the notions of catalecticant matrix and Waring rank. We will subsequently prove that if the $h$-vector is minimum then WLP holds, on the contrary, if it is maximum then the WLP fails. We will then conclude the chapter trying to generalise such properties to hypersurfaces in $\PP^4$ with or without vanishing Hessian.

    \begin{definition}\label{perazzo} Fix $N\ge 4$. A {\em Perazzo} hypersurface  $X\subset \PP^N$ of degree $d\ge 3$ is a hypersurface defined by a form $f\in K[x_0,\cdots ,x_n,u_1\cdots ,u_m]$ of the following type:
    $$ f=x_0p_0+x_1p_1+\cdots +x_np_n+g$$
    where
    \begin{itemize}
        \item $n,m$ are positive integer numbers with $n+m=N$, $n,m\ge 2$;
        \item $p_i\in K[u_1,\cdots ,u_m]_{d-1}$ are algebraically dependent, but linearly independent;
        \item $g\in K[u_1,\cdots ,u_m]_{d}$.
    \end{itemize}
    We will sometimes refer to a form $f$ as Perazzo hypersurface if $V(f)$ is a Perazzo hypersurface.
    \end{definition}
    
    The assumptions on $d,m$ and $n$ are necessary because, otherwise, the forms $p_0,\dots,p_n$ cannot be simultaneously algebraically dependent and linearly independent. Some authors add the assumptions "reduced" and "irreducible" to this definition (see, for instance, \cite[Definition 3.12]{gondim}); but, for our treatment, these hypotheses are not so necessary, therefore, we will add them when needed.
    
    \begin{example}
    We have seen a first example of this kind of hypersurface in Example \ref{perazzo5} with the polynomial $f$ defined as $f(u,v,x,y,z)=u^2x+uvy+v^2z$.
    \end{example}
    
    As we stated in Proposition \ref{partialalgebraic} and in Proposition \ref{basicelimination}, a hypersurface defined by a polynomial $f$ has vanishing hessian if and only if the partial derivatives of $f$ are algebraically dependent, and it is a cone if and only if they are linearly dependent. It follows that Perazzo hypersurfaces, introduced in Definition \ref{perazzo}, have all vanishing hessian. We will see in Proposition \ref{cone} that, in general, they are not cones. In particular, by Corollary \ref{watanabecorollario} and setting $S=K[X_0,\dots,X_n,U_1,\dots,U_m]$, we have that the Gorenstein $K$-algebra $A=S/\Ann_S(f)$, associated to a Perazzo hypersurface $f$, fails to have the Strong Lefschetz property.
    
    \pause
    We characterise now Perazzo hypersurfaces that are cones. This case will be excluded in the rest of this chapter.
    
    \begin{proposition}\label{cone}
        Let $N=n+m$, with $n,m\ge 2$. 
        
        We consider $f\in K[x_0, \cdots , x_n, u_1\cdots ,u_m]$ a Perazzo hypersurface of degree $d$ such that $X=V(f)\subseteq\PP^N$ is a cone. Then $X$ is projectively equivalent to a cone of vertex a point over a Perazzo hypersurface $Y\subseteq\PP^{N-1}$ of degree $d$.
    \end{proposition}
    \begin{proof}
        By assumption $f$ can be written as 
        $$f=x_0p_0+x_1p_1+\cdots +x_np_n+g$$
        where, $p_0,\dots,p_n\in K[u_1,\cdots ,u_m]_{d-1}$ are linearly independent and algebraically dependent forms, and $g\in K[u_1,\cdots ,u_m]_{d}$. Thus, we can explicitly compute the partial derivatives of $f$ as follows:
        \begin{equation*}
            \nabla f=
            \begin{pmatrix}
            p_0\\
            p_1\\
            \vdots\\
            p_n\\
            x_0\frac{\partial}{\partial u_1}p_0+ x_1\frac{\partial}{\partial u_1}p_1+ \dots+ x_n\frac{\partial}{\partial u_1}p_n+ \frac{\partial}{\partial u_1}g\\
            \vdots\\
            x_0\frac{\partial}{\partial u_m}p_0+x_1\frac{\partial}{\partial u_1}p_1+ \dots+ x_n\frac{\partial}{\partial u_m}p_n+ \frac{\partial}{\partial u_m}g
            \end{pmatrix}.
        \end{equation*}
        By Proposition \ref{basicelimination}, since $X$ is a cone, there exists a linear relation among the partial derivatives of $f$. So, there exist coefficients $\lambda_0,\dots,\lambda_n,\mu_1,\dots,\mu_m\in K$, not all zero, such that
        \begin{equation}\label{fconeequation}
            \lambda_0p_0+\dots+\lambda_np_n+\mu_1\left(x_0\frac{\partial p_0}{\partial u_1}+ \dots+ \frac{\partial g}{\partial u_1}\right)+\dots+\mu_m\left(x_0\frac{\partial p_0}{\partial u_m}+\dots\right)=0.
        \end{equation}
        
        Grouping together all the coefficients of the variables $x_0,\dots,x_n$ we obtain an equivalent system of equations
        
        \begin{equation*}
            \begin{cases}
                \mu_0\frac{\partial}{\partial u_1}p_0+ \mu_1\frac{\partial}{\partial u_2}p_0+ \dots+ \mu_m\frac{\partial}{\partial u_m}p_0=0\\
                \mu_0\frac{\partial}{\partial u_1}p_1+ \mu_1\frac{\partial}{\partial u_2}p_1+ \dots+ \mu_m\frac{\partial}{\partial u_m}p_1=0\\
                \dots\\
                \mu_0\frac{\partial}{\partial u_1}p_n+ \mu_1\frac{\partial}{\partial u_2}p_n+ \dots+ \mu_m\frac{\partial}{\partial u_m}p_n=0\\
                \lambda_0p_0+\dots+\lambda_np_n+ \mu_0\frac{\partial}{\partial u_1}g+ \mu_1\frac{\partial}{\partial u_2}g+ \dots+ \mu_m\frac{\partial}{\partial u_m}g=0
            \end{cases}.
        \end{equation*}
        Looking at the proof of Proposition \ref{basicelimination}, if one of the coefficients in (\ref{fconeequation}) is different from zero, then the related variable can be eliminated with a linear transformation. Being the $p_i$ linearly independent, one of the $\mu_i$ has to be different from $0$. Up to a renumbering of the variables, we can suppose $\mu_m\ne 0$. Following Proposition \ref{basicelimination}, we can now introduce a change of variables that will eliminate the variable $u_m$ from $f$
        \begin{equation}\label{changeofvariablescone}
            \begin{split}
                &y_0=x_0-\frac{\lambda_0}{\mu_m}u_m,\dots,y_n=x_n-\frac{\lambda_n}{\mu_m}u_m,\\
                w_1=u_1-&\frac{\mu_1}{\mu_m}u_m,\dots,w_{m-1}=u_{m-1}-\frac{\mu_{m-1}}{\mu_m}u_m, w_m=u_m.
            \end{split}
        \end{equation}
        
        One can notice that the coefficients $\mu_1,\dots,\mu_m$ are also coefficients of a vanishing linear combination of the partial derivatives of each of the forms $p_0,\dots,p_n$. Moreover, the new variables $w_1,\dots,w_m$ depend only on the variables $u_1,\dots,u_m$. Thus the same change of variables also eliminates from $p_0,\dots,p_n$ the same common variable $u_m$. We name $q_i(w_1,\dots,w_m):=p_i(u)_{|(w_1,\dots,w_m)}$ the polynomials $p_i$ after the change of variables (\ref{changeofvariablescone}). Clearly, $q_0,\dots,q_n\in K[w_1,\dots,w_{m-1}]_{d-1}$ and they continue to be linearly independent and algebraically dependent. Using the change of variables (\ref{changeofvariablescone}), $X$ is projectively equivalent to a cone of vertex a point over $Y=V(h)\subseteq\PP^{N-1}$ where $h\in K[y_0,\dots,y_n,w_1,\dots,w_{m-1}]_d$.
        We just need to prove that $h$ is a Perazzo hypersurface. Defining
        \begin{equation*}
            \begin{split}
                g'(w_1,\dots,w_{m}):=g(w_1,\dots,w_{m})+ \frac{w_m}{\mu_m}(\lambda_0q_0(w_1,\dots,w_{m-1})+\dots\\
                \dots +\lambda_n q_n(w_1,\dots,w_{m-1}))
            \end{split}
        \end{equation*}
        we can write $h$ as
        \begin{equation*}
            \begin{split}
                h(y_0,\dots,y_n,w_1,\dots,w_{m-1}) =y_0q_0(w_1,\dots,w_{m-1})+\dots\\
                \dots+y_n q_n(w_1,\dots,w_{m-1})+ g'(w_1,\dots,w_{m}).
            \end{split}
        \end{equation*}
        Finally, we can notice that $g'$ is a homogeneous polynomial of degree $d$ that does not depend on $w_m$. Thus, we have just proved that $h$ is a Perazzo hypersurface.
    \end{proof}
    \begin{important}
        We are interested only in hypersurfaces which are not cones, so from now on, when we refer to a hypersurface, we tacitly assume its property of not being a cone.
    \end{important}

    As we have seen in Theorem \ref{quinary}, in the case $N=4$, the Perazzo 3-folds are the building block for all hypersurfaces with vanishing Hessian. In fact, every possible form $f$ with vanishing Hessian is in the extension $K[x_3,x_4][\Delta]$ where $\Delta$ is a Perazzo hypersurface defined as $\Delta=x_0p_0+x_1p_1+x_2p_2$ with $n=2$ (see Remark \ref{importantremark}).
    
    \pause
    In general, in the form $f$, also the powers of $\Delta$ could appear. In the case of cubic hypersurfaces this fact does not occur (see Theorem \ref{cubicsinp4}), so the Perazzo 3-folds fully characterise this case. Also when $\Delta$ appears in $f$ at degree $1$, $f$ is, after a renaming of the $p_i$, a Perazzo hypersurface with $n=2$ and $m=2$. Therefore, if we want to study Gorenstein $K$-algebras of quinary forms with vanishing Hessian, a natural way is first to consider the case of Perazzo hypersurfaces.
    
    \pause
    Summarising, the aim of this chapter is to consider Perazzo hypersurfaces $f$ of the form $f=x_0p_0(u,v)+x_1p_1(u,v)+x_2p_2(u,v)+g(u,v)$ trying to understand if the associated Gorenstein $K$-algebra, which lacks the SLP, has or not the Weak Lefschetz property. To this end we first introduce the concept of catalecticant matrix and Waring rank and we will then use them to compute the Hilbert vector of such algebras.
    
    \section{Catalecticant matrices, Waring rank and binary forms}\label{catal}
    
    For this section we will use \cite{ottaviani}, \cite{iarrobinokanev} and \cite{bernardigimiglianoida} as principal references. We start considering a generic homogeneous polynomial of degree $t$ in two variables, $h\in K[u,v]_t$. We ask ourselves if it is possible to write $h$ as sum of $t$-th powers of some linear forms and, in the affirmative case, how to find such expression. Furthermore, we look for a combination which involves the lowest number of distinct linear forms\footnote{Nowadays, this problem is viewed in a more general context as a problem of decomposition of tensors. The name of this classic problem is Big Waring Problem.}. We name \emph{Waring rank} of $h$ the minimum integer $e$ such that $h=l_1^t+\dots+l_e^t$ where the $l_i$ are distinct binary linear forms. We now discuss about this problem and its relation with the catalecticant matrices that will be defined later.
    
    \pause
    We can interpret the polynomial $h$ as a point of the projective space $\PP^t$ which can be canonically identified as $\PP(K[u,v]_t)$. In fact, if we write 
    \begin{equation}\label{formabinariaeq}
        h(u,v)=\sum _{i=0}^{t} \binom{t}{i} h_i u^{t-i} v^{i}
    \end{equation}
    then there is a 1-1 correspondence with the element $[h]:=[h_0:\dots:h_t]\in\PP^t$ up to a non constant scalar. This correspondence has a good behavior with respect to change of variables and projectivity. If fact, a change of the variables $u$ and $v$ corresponds to a projectivity in $\PP^t$ and viceversa. 
    
    \begin{remark*}
       We have made this choice because the catalecticant matrices of $h$, as defined later in Definition \ref{catalecticant}, are much simpler to be written. The arguments of this section are still true also defining $h$ without the binomial coefficients.
    \end{remark*}
    
    In particular, forms that can be written as pure power of a linear form are exactly the image's elements of the Veronese embedding $\nu:\PP^1\to\PP^t$ defined as $$\nu([a:b])=[a^t:a^{t-1}b:\dots:b^t].$$
    
    This image $C_t$ is a smooth variety also called \emph{rational normal curve}. In this variety, there are exactly the elements with $1$ as Waring rank. If we consider the line generated by two different elements $p,q\in C_t$, then all the elements in this line have Waring rank at most $2$. With this idea in mind we consider, for all $r\ge 1$, the $r$-secant variety of $C_t$ defined as
    $$\sigma_r(C_t)=\overline{\bigcup\left\{\langle p_1,\ldots,p_r\rangle:\,p_1,\ldots,p_r\in C_t \text{ distinct}\right\}}.$$
    Clearly, we have $\sigma_1(C_t)=C_t$. We now recall some basic properties of the rational normal curve and its secant varieties.
    \begin{proposition}
        Let $\sigma_r(C_t)$ be the $r$-secant variety of the rational normal curve $C_t$ as defined above, then:
        \begin{enumerate}
            \item if $\{P_0,\dots,P_s\}$ are $s+1$ different points in $C_t$, $s\le t$, then they are linearly independent;
            \item the dimension of $\sigma_r(C_t)$ is $\min\{2r-1, t\}$;
            \item there are a series of inclusions $C_t\subsetneq\sigma_2(C_t)\subsetneq\dots\subsetneq\PP^t$;
            \item $\sigma_r(C_t)$ is the singular locus of $\sigma_{r+1}(C_t)$, whenever $\sigma_{r+1}(C_t)\subsetneq\PP^t$.
        \end{enumerate}
    \end{proposition}
    \begin{proof}
        Point 1. can be proved in the general case when $s=t$ and $P_i\ne [0:\dots:0:1]$. In this case the matrix which involves the coordinates of that points is
        $$V=\begin{pmatrix}
            1 & b_0 & b_0^2 & \dots & b_0^{t}\\
            1 & b_1 & b_1^2 & \dots & b_1^{t}\\
            1 & b_2 & b_2^2 & \dots & b_2^{t}\\
            \vdots & \vdots & \vdots & &\vdots \\
            1 & b_t & b_t^2 & \dots & b_t^{t}\\
        \end{pmatrix}.$$
        V is a Vandermonde matrix, and its determinant $\prod_{i<j}(b_j-b_i)$ is different from $0$ by assumption. If $P_0=[0:\dots:0:1]$, then surely the first coordinate of the other points is different from zero. Arguing as in the previous case we obtain that the points $P_1,\dots,P_t$ are linearly independent being the matrix of their coordinate non-singular. Finally, the elements $P_0,\dots,P_t$ are still linearly independent, because the first coordinate of $P_0$ is zero.
        
        For points 2.- 4. a proof can be found in \cite[Proposition 1.2.2 and Corollary 1.2.3]{russo}.
    \end{proof}
    
    By construction, the general element of $\sigma_r(C_t) \setminus\sigma_{r-1}(C_t)$ has Waring rank $r$. But, in $\sigma_r(C_t)$ there are also elements that are limit position of elements in the linear span of $r$ distinct points; thus, they have Waring rank greater than $r$. For example, as a consequence of Proposition \ref{lemmarank1}, $\sigma_{2}(C_t)$ contains the tangential surface of $C_t$, and $\sigma_{3}(C_t)$ contains the osculating $3$-fold of $C_t$, where the tangential surface of $C_t$, $TC_t$, is the closure of the union of the embedded tangent lines to $C_t$. With some computations, is possible to see that the tangent line at some point $l_1^t\in C_t$ is the set of forms that can be written as $l_1^{t-1}l_2$, where $l_2$ is a linear form. Similarly, the osculating $3$-fold of $C_t$, $T^2C_t$ is the closure of the union of the embedded osculating planes to $C_t$, and the osculating plane at $l_1^t$ is the set of forms that can be written as $l_1^{t-2}l_2l_3$, with $l_1,l_2$ forms of degree $1$ not necessarily distinct.
    
    \pause
    Being $\sigma_{r}(C_t)$ the smallest variety containing forms with Waring rank $r$, it is natural, given a form $h$, to seek the smallest integer $r$ such that $h\in\sigma_r(C_t)$. The number $r$ is called \emph{symmetric border rank} of $h$. This quantity can be characterised using the so called "catalecticant matrices" of a given form $h$.
    
    \begin{definition}\label{catalecticant}
        Let $h\in K[u,v]_t$ be a homogeneous polynomial. For every $k$, $0\le k\le t$, we name $\Cat_{k}(h)$ the \emph{$k$-catalecticant} matrix, or \emph{$k$-Hankel} matrix, of $h$ a $(t-k+1)\times (k+1)$ matrix $\Cat_k(h)$ with coefficients in $K$ defined as follows. These matrices are chosen in such a way that the following relation is true:
        \begin{equation*}
            \Cat_{t-k}(h)\cdot
            \begin{pmatrix}
                u^{t-k}\Bstrut{}\\
                \binom{t-k}{1}u^{t-k-1}v\Bstrut{}\\
                \binom{t-k}{2}u^{t-k-2}v^2\Bstrut{}\\
                \vdots\Bstrut{}\\
                v^{t-k}
            \end{pmatrix}=
            \begin{pmatrix}
                \frac{\partial^k h}{\partial u^k}\Bstrut{}\\
                \frac{\partial^k h}{\partial u^{k-1}\partial v}\Bstrut{}\\
                \frac{\partial^k h}{\partial u^{k-2}\partial v^2}\Bstrut{}\\
                \vdots\Bstrut{}\\
                \frac{\partial^k h}{\partial v^k}
            \end{pmatrix}\quad \forall k=0,1,\dots,t.
        \end{equation*}
        If we write $h$ as in (\ref{formabinariaeq}), then the $k$-catalecticant matrix of $h$ is
        \begin{equation*}
            \Cat_k(h)=
            \begin{pmatrix}
                h_0 & h_1 & \cdots & h_{k} \\
                h_1 & h_2 & \cdots & h_{k+1} \\
                \vdots & \vdots & \ddots & \vdots \\
                h_{t-k} & h_{t-k+1} & \cdots & h_{t} 
            \end{pmatrix}.
        \end{equation*}
    \end{definition}
    
    Equivalently, the matrix $\Cat_k(h)$ can be interpreted in the following way. Consider $K[U,V]$ the ring of the differential operators and, for every $k$, $0\le k\le t$, consider the map
    $$A(h)_k:K[U,V]_k\to K[u,v]_{t-k}$$
    that acts derivating the form $h$. If we fix the standard basis of the monic monomials on $K[U,V]_k$, and the basis $\left\{ \binom{d}{k}\binom{t-k}{i}k!\,\, u^{t-k-i} v^{i}:\,0\le i\le t-k\right\}$ on $K[u,v]_{t-k}$, then the matrix associated to $A(h)_k$ is exactly $\Cat_k(h)$. We also notice that $\Cat_k(h)^T=\Cat_{t-k}(h)$.
    
    \pause
    The first useful meaning of the catalecticant matrices is that they can be used to explicit compute the linear forms that can be used to write a form $h$ as $h=l_1^t+\dots+l_r^t$. We remark that the following statement holds because $K$ is an algebraically closed field of characteristic zero.
    
    \begin{proposition}\label{lemmarank4}
        Let $h\in K[u,v]_t$ and let $l_i$ be distinct linear form for $i=1,\dots,e$. Let $r\le t-e$. Then $h=\sum_{i=1}^e l_i^t$ if and only if $\Imm \Cat_r(h)\subseteq\langle (l_1)^{t-r},\dots,(l_e)^{t-r}\rangle$.
    \end{proposition}
    \begin{proof}
        See \cite[Corollary 1.2]{ottaviani}.
    \end{proof}
    
    For the seek of completeness, we introduce also the algorithm which allows to explicitly compute the linear forms that appear in the minimal decomposition of a given form. This algorithm is known as \emph{Sylvester algorithm}; it was first discovered by J. J. Sylvester in \cite{sylvester}.
    
    \pause
    Let $h\in K[u,v]_t$ be a form and let $\Cat_k(h)$ be its catalecticant matrices as defined in Definition \ref{catalecticant}. Suppose now that the Waring rank of $h$ is $e$ and so $h=\sum_{i=1}^e l_i^t$. We will see in Proposition \ref{lemmarank3} how it is possible to compute $r$, the symmetric border rank of $h$. The Waring rank $e$ and the linear forms $l_i$ are computed as follows:
    
    \pause
    1> INPUT: a form $h$ and its symmetric border rank $r$
    
    2> START
    
    3> set $k=r$
    
    4> take an element $g$ in the kernel of $A(h)_k$ viewed as an element in $K[u,v]_k$
    
    5> decompose $g=g_1^*\dots g_k^*$ as product of linear forms
    
    6> set $g_i$ as the apolar linear form of $g_i^*$\footnote{If $g=au+bv\in K[u,v]$ is a binary linear form, then its polar form is $g^\perp=bu-av$}
    
    6> IF $g_i$ are all distinct then:
    
    7>\hspace{1.2cm}compute the coefficients $\lambda_1,\dots,\lambda_k$ such that $f=\sum_{i=1}^k \lambda_k g_i^t$
    
    8>\hspace{1.2cm}set $\xi_i$ as a $t$-th root of $\lambda_i$ for all $i$
    
    9>\hspace{1.2cm}set $l_i=\xi_i g_i$ for all $i$
    
    10>\hspace{1cm}set $e=k$
    
    11> ELSE
    
    12>\hspace{1cm}$k=k+1$ and JUMP TO 4
    
    13> STOP
    
    13> OUTPUT: the integer $e$ and the forms $l_1,\dots,l_r$ 
    
    \pause
    As we have said above, we now state some results about Waring rank and catalecticant matrices with the purpose of proving Proposition \ref{rank}. To do that we give the definition of join of two projective varieties.
    
    \begin{definition}
        Given two projective varieties $X$ and $Y$, the \emph{join} of $X$ and $Y$ is the closure of the union of all lines joining a point of $X$ with a different point of $Y$. We refer to this variety as $J(X,Y)$.
    \end{definition}
    
    \begin{proposition}\label{lemmarank1}
        Let $\sigma_{a,b}(C_t)$ denote the set of all the forms of degree $t$ with symmetric border rank $a$ and Waring rank $b$. With the above notations we have that
        $$\forall r,\, 2\le r\le \frac{t+1}{2}:\quad \sigma_r(C_t)\setminus\sigma_{r-1}(C_t)=\sigma_{r,r}(C_t)\cup\sigma_{r,t-r+2}(C_t).$$
        Moreover we have,
        $$\sigma_{2,t}(C_t)=TC_t\setminus C_t;$$
        $$\sigma_{r,t-r+2}(C_t)=J(TC_t,\sigma_{r}(C_t))\setminus\sigma_{r-1}(C_t)\quad\text{for all } r,\,3\le r< \frac{t+2}{2}.$$
    \end{proposition}
    \begin{proof}
        See \cite[Theorem 23, and Corollary 26]{bernardigimiglianoida}.
    \end{proof}
    
    \begin{proposition}\label{lemmarank2}
        Let $k\le e\le\frac{t}{2}$, and $h\in K[u,v]_t$. Then $[h]\in\sigma_k(C_t)\iff \rank \Cat_e(h)\le k\iff$ the symmetric border rank of $h$ is greater or equal to $k$. In particular, the first integer $r$ for which the catalecticant matrix $\Cat_r(h)$ has not the expected rank (i.e. it has rank strictly less than $r+1$) is the symmetric border rank of $h$.
    \end{proposition}
    \begin{proof}
        See \cite[Theorem 1.3]{ottaviani}.
    \end{proof}
    
    \begin{proposition}\label{lemmarank3}
        Let $h\in K[u,v]_t$, and $s=\rank\Cat_{\lfloor\frac{t}{2}\rfloor}(f)$. Then for every $k$, $s\le k\le t-s$, it holds $\rank\Cat_k(h)=s$.
    \end{proposition}
    \begin{proof}
        See \cite[Theorem 1.43]{iarrobinokanev}.
    \end{proof}
    
    The following final proposition of this section is a key result for our treatment. In fact, we will use it in Theorem \ref{lowcharacterisation} to characterise some Perazzo forms using their Waring rank and the rank of their catalecticant matrices.
    
    \begin{proposition}\label{rank}
        We fix an integer $t\ge 4$, and we keep the notations introduced in Definition \ref{catalecticant}. The following statements hold:
        \begin{itemize}
            \item[(1)] If $\rank \Cat_k(h)=1$ for some $1\le k \le \lfloor \frac{t}{2} \rfloor$ (and, hence, for all $k$ in the same interval), then $h=\ell ^{t}$ for some $\ell \in K[u,v]_1$.
            \item[(2)] If $\rank \Cat_k(h)=2$ for some $2\le k \le \lfloor \frac{t}{2} \rfloor$ (and, hence, for all $k$ in the same interval), then either $h=\ell _1 ^{t}+\ell _2 ^{t}$ or $h=\ell _1^{t-1}\ell _2$ for some $\ell_1, \ell _2 \in K[u,v]_1$.
            \item[(3)] If $\rank \Cat_k(h)=3$ for some $3\le k \le \lfloor \frac{t}{2} \rfloor$ (and, hence, for all $k$ in the same interval), then either $h=\ell _1 ^{t}+\ell _2 ^{t}+(\lambda \ell _1+\mu \ell _2)^{t}$ or $h=\ell _1^{t}+\ell _2^{t-1}(\lambda \ell _1+\mu \ell_2)$ for some $\ell_1, \ell _2 \in K[u,v]_1$ and $\lambda, \mu \in K^*$.
        \end{itemize}
    \end{proposition}
    \begin{proof}
    Let $r$ be any integer such that $r+1\leq k$. By Proposition \ref{lemmarank2}, it follows that all the minors of order $r+1$ of $\Cat_k(h)$ vanish if and only if $[h]\in\sigma_r(C_{t})$. For $r=1$, this gives immediately (1). For $r=2$, we get that if $\Cat_k(h)$ has rank $2$, then  $[h]\in\sigma_2(C_{t})\setminus C_t$. By Proposition \ref{lemmarank1}, it follows that either $h$ has Waring rank $2$ or $h\in TC_{t}$; this proves (2). Similarly, for $r=3$, $\rank \Cat_k(h)=3$ implies that $h\in \sigma_3(C_{t})\setminus\sigma_2(C_{t})$. So, again by Proposition \ref{lemmarank1}, either the Waring rank of $h$ is $3$, or $h$ belongs to the join of $C_{t}$ and its tangential surface. This proves (3) and concludes the proof.
    \end{proof}
    
    \begin{example}
        1.- We want to study the polynomial $h_1=u^3+3uv^2\in \C[u,v]_3$. Its catalecticant matrices are
        \begin{equation*}
            \Cat_0(h_1)=\begin{pmatrix}
            1\\
            0\\
            1\\
            0
            \end{pmatrix},\,
            \Cat_1(h_1)=\begin{pmatrix}
            1&0\\
            0&1\\
            1&0
            \end{pmatrix},\,
        \end{equation*}
        \begin{equation*}
            \Cat_2(h_1)=\Cat_1(h_1)^T,\,
            \Cat_3(h_1)=\Cat_0(h_1)^T.
        \end{equation*}
        
        By Proposition \ref{lemmarank2} we deduce that $[h_1]\not\in C_3$. Since the dimension of $\sigma_2(C_3)\subset\PP^3$ is three, we also deduce that $[h_1]\in \sigma_2(C_3)$. By Proposition \ref{lemmarank1}, we know that $\sigma_2(C_3)=\sigma_{2,2}(C_3)\cup\sigma_{2,3}(C_3)$. So there are only two possibilities:
        $$h_1=l_0^3+l_1^3 \text{, or } h_1=l_0l_1^2$$
        for a suitable choice of two linear forms $l_0$, and $l_1$. In the second case $h_1$ can be written only as sum of three third powers which is indeed the maximal number. We now use the Sylvester algorithm to compute such forms. Being $2$ the symmetric border rank of $h_1$, we consider the kernel of $A(h_1)_2$ that is generated by $$g=U^2-V^2=(U-V)(U+V),$$
        where $U=\partial_{u}$ and $V=\partial_v$. A trivial computation shows that $(u-v)^\perp=-u-v$, $(u+v)^\perp=u-v$ and so
        $$h_1=\frac{1}{2}(u+v)^3+\frac{1}{2}(u-v)^3.$$
        Since we are considering an algebraically closed field of characteristic zero, we can write $h_1=l_0^3+l_1^3$ taking the $3$-rd root of $\frac{1}{2}$.
        
        2.-We now study the polynomial $h_2=u^4-2u^3v+2uv^3-v^4\in \C[u,v]_4$. Its catalecticant matrices are
        \begin{equation*}
            \Cat_0(h_2)=\begin{pmatrix}
            1\\
            -\frac{1}{2}\\
            0\\
            \frac{1}{2}\\
            -1
            \end{pmatrix},\,
            \Cat_1(h_2)=\begin{pmatrix}
            1&-\frac{1}{2}\\
            -\frac{1}{2}&0\\
            0&\frac{1}{2}\\
            \frac{1}{2}&-1
            \end{pmatrix},\,
            \Cat_2(h_2)=\begin{pmatrix}
            1&-\frac{1}{2}&0\\
            -\frac{1}{2}&0&\frac{1}{2}\\
            0&\frac{1}{2}&-1
            \end{pmatrix},
        \end{equation*}
        \begin{equation*}
            \Cat_3(h_2)=\Cat_1(h_2)^T,\,
            \Cat_4(h_2)=\Cat_0(h_2)^T.
        \end{equation*}
        
        The matrix $\Cat_2(h_2)$ has rank $2$; thus, by Proposition \ref{rank}, there exist two linear forms $l_0,l_1$ such that one and only one of the following cases occur:
        $$h_2=l_0^3+l_1^3 \text{, or } h_2=l_0l_1^2.$$
        Using again the Sylvester algorithm, we start considering the kernel of $A(h_2)_2$: it is generated by $U^2+2UV+U^2=(U+V)^2$. So we have to skip this case and consider the kernel of $A(h_2)_3$. Again one element in it is $U^3+2U^2V+UV^2=U(U+V)^2$. Finally, we can pick the element $U^3V+UV^3=UV(U+iV)(U-iV)$ in the kernel of $\Cat_4(h_2)$, where $i$ is the imaginary unit. Taking for any linear factor its polar form, we obtain the forms $-v,u,iu-v$ and $-iu-v$. Thus a straightforward computation proves that
        $$h_2=-v^4+u^4+\frac{i}{4}(iu-v)^4-\frac{i}{4}(iu+v)^4.$$
        We notice that the form $U+V$ appears as multiple factor in the elements in the kernel of $A(h_2)_2$ and $A(h_2)_3$. This fact can be explained as follows: the dual element of $U+V$ is $u-v$ and the form $h_2$ can be factorised as $h_2=(u+v)(u-v)^3$.
    \end{example}

    \section{Construction of the Hilbert vector}\label{constructionofthehilbertvector}
    Let $f$ be a Perazzo $3$-fold of degree $d$ of the form $f=x_0p_0(u,v)+x_1p_1(u,v)+x_2p_2(u,v)+g(u,v)\in K[x_0,x_1,x_2,u,v]_d$. Let $S=K[y_0,y_1,y_2,U,V]$ the ring of the differential operators. We now focus on the construction of the Hilbert vector of the $K$-algebra $A=S/\Ann_S(f)$. First of all we notice that the condition of $V(f)$ not being a cone implies that $\dim A_1=\dim A_{d-1}=5$ because these quantities count the number of linearly independent elements among the partial derivatives. In particular, for $d=3$, the only possible $h$-vector is $(1,5,5,1)$. So, from now on, we will assume that $d\ge 4$. Following Section \ref{catal}, we write
      \begin{equation*}
          \begin{array}{rcl} p_0(u,v) & =& \sum _{i=0}^{d-1} {\binom{d-1} {i}}a_iu^{d-1-i}v^{i}, \\
          p_1(u,v) & = & \sum _{i=0}^{d-1}{\binom{d-1} i}b_iu^{d-1-i}v^{i}, \\
          p_2(u,v) & = & \sum _{i=0}^{d-1}{\binom{d-1} i}c_iu^{d-1-i}v^{i}, \text{ and } \\
          g(u,v) & = & \sum _{i=0}^{d}{\binom{d} i}g_iu^{d-i}v^{i}. \end{array} \end{equation*}
    
    For any  $0\le k\le d-1$,  we define four groups of catalecticant matrices with respect to the forms $p_0,p_1,p_2$, and $g$:
    $$
    A_k:=\begin{pmatrix}
     a_0 & a_1 & \cdots & a_{k} \\
     a_1 & a_2 & \cdots & a_{k+1} \\
     \vdots & \vdots & & \vdots \\
      a_{d-1-k} & a_{d-k} & \cdots & a_{d-1} 
    \end{pmatrix},
    \quad 
    B_k:=\begin{pmatrix}
     b_0 & b_1 & \cdots & b_{k} \\
     b_1 & b_2 & \cdots & b_{k+1} \\
     \vdots & \vdots & & \vdots \\
      b_{d-1-k} & b_{d-k} & \cdots & b_{d-1} 
    \end{pmatrix},  
    $$
    
    $$
    C_k=\begin{pmatrix}
     c_0 & c_1 & \cdots & c_{k} \\
     c_1 & c_2 & \cdots & c_{k+1} \\
     \vdots & \vdots & & \vdots \\
      c_{d-1-k} & c_{d-k} & \cdots & c_{d-1} 
    \end{pmatrix} \text{, and }\, G_k=
    \begin{pmatrix}
        g_0 & g_1 & \cdots & g_k\\
        g_1 & g_2 & \cdots & g_{k+1}\\
        \vdots & \vdots & & \vdots \\
        g_{d-k} & g_{d-k+1} & \cdots & g_d
    \end{pmatrix}.$$
    
    The matrices $A_k$, $B_k$, $C_k$ and $G_k$ are the building blocks of the  matrices $M_k$, $N_k$ and $N'_k$ that will play an important role in the proof of our main results. They are defined as follows:
    \begin{equation}\label{mnnk}
        M_k:=
        \begin{pmatrix}
            A_{k-1} \vert B_{k-1} \vert C_{k-1}
        \end{pmatrix}, \quad N_k:=\begin{pmatrix}
            A_{k} \\ \hline B_{k} \\ \hline C_{k}
        \end{pmatrix}
        \text{ and } \quad
        N'_k:= \begin{pmatrix}
            A_{k} \\ \hline B_{k} \\ \hline C_{k} \\ \hline G_{k}
        \end{pmatrix}.
    \end{equation}
    
    
    The matrices $N_k$ and $M_{k+1}$ contain the same 3 blocks of size $(d-k)\times (k+1)$, but they are arranged in a different way; so, they can have a different rank. The matrix $M_k$ has size $(d-k+1)\times(3k)$, while the matrix $N'_k$ has size $(4d-4k+1)\times(k+1)$. Being $A_k$, $B_k$, and $C_k$ catalecticant matrices, we obtain the relation $M_k=N_{d-k}^T$. So, we have $\rank M_k=\rank N_{d-k}$.
    
    \pause
    We now state a proposition that links the values of the Hilbert vector of $A=S/\Ann_S(f)$ and the matrices we have just defined. Before that we recall few facts in the following lemma.
    \begin{lemma}
        The following statements hold true:
        \begin{itemize}
            \item the dimension, as $K$-vector space, of the $i$-homogeneous part of any polynomial ring $K[x_0,\dots,x_n]$ can be computed as
            $$\dim_K K[x_0,\dots,x_n]_k=\binom{n+k}{k}=\binom{n+k}{n};$$
            \item if $h\in K[u,v]_t$ is written as $h(u,v)=\sum_{k=0}^t\binom{t}{k}h_ku^{t-k}v^k$, then for all integers $m,n\in \N^*$ such that $n\le m\le t$ it holds
            $$\frac{\partial^m h}{\partial u^{m-n}\partial v^n}=t\cdots(t-m+1)\sum_{k=0}^{t-m}\binom{t-m}{k}h_{k+t}u^{t-m-k}v^k.$$
        \end{itemize}
    \end{lemma}
    
    \begin{proposition}\label{hilbert function}
    Let $S=K[y_0,y_1,y_2,U,V]$ be the ring of the differential operators. Let $f$ be a Perazzo $3$-fold of degree $d$ in $\PP^4$, and let $h=(h_0, h_1, \ldots, h_d)$ be $h$-vector associated to $A=S/\Ann_S(f)$. Let $M_k$, and $N_k'$ be the matrices, arising from $f$, defined above. Then
    $h_0=h_d=1, h_1=h_{d-1}=5$ and, for $2\leq k\leq d-2$, $h_k=4k+1-m_k-n_k$, where $m_k=3k-\rank M_k$ and  $n_k=k+1-\rank N'_k$ (or equivalently $h_k=\rank M_k+\rank N'_k$).
    \end{proposition}
    \begin{proof}The key argument to compute the value $h_k$ is to determine the number of linear relations among the $k$-partial derivatives of $f$, which is $\dim\Ann_S(f)_k$; in fact,
    $$h_k=\dim A_k=\dim S_k-\dim \Ann_S(f)_k={\binom{4+k}{k}}-\dim \Ann_S(f)_k.$$
    
    As we have noticed in Remark \ref{gorenstainsymmetric}, the $h$-vector of an Artinian Gorenstein algebra is symmetric and, hence, we only have to compute the first half of the vector, so that $h_k$ for $0\le k \le \lfloor \frac{d}{2} \rfloor $.
    
    \pause
    We now compute by hand $\dim \Ann_S(f)_k$ for any $k$,  $0\le k \le \lfloor \frac{d}{2} \rfloor $. The first values $h_0$ and $h_1$ are the simplest to compute. In fact, as we have already argued, their value are $h_0=1$ and $h_1=5$ by the $K$-linear independence of the $p_i$'s.
    
    \pause
    Being $f$ linear in the variables $x_0, x_1$, and $x_2$, we observe that, for any $k\geq 2$, $\Ann_S(f)_k$ clearly contains $(y_0,y_1,y_2)^{k-i}(U,V)^{i}$, for $0\leq i\leq k-2$. Therefore 
    $$\dim A_k\leq {\binom{4+k}{k}}-\sum_{i=0}^{k-2}(i+1){\binom{k-i+2}{2}}=4k+1.$$
    We have to consider now the non trivial linear combination that can arise from the two parts $(y_0,y_1,y_2)(U,V)^{k-1}$ and $(U,V)^k$. Therefore we consider the numbers $m_k$ and $n_k$ as follows:
    $$m_k=\dim(\Ann_S(f)_k\cap (y_0,y_1,y_2)(U,V)^{k-1}),$$
    $$n_k=\dim(\Ann_S(f)_k\cap (U,V)^k).$$
    
    The computation of $m_k,n_k$ will complete the proof. In fact, due to the expression of $f$, if $\varphi\in\Ann_S(f)$ and $\varphi=\varphi_1+\varphi_2$ with $\varphi_1\in(y_0,y_1,y_2)(U,V)^{k-1}$ and $\varphi_2\in(U,V)^k$, then both $\varphi_1$ and $\varphi_2$ are in $\Ann_S(f)$. Therefore we obtain $\dim A_k=4k+1-m_k-n_k$.
    
    \pause
    We enter now in the final part of this proof computing $m_k$, and then $n_k$. For $m_k$, we consider a general polynomial of degree $k$ in $(y_0,y_1,y_2)(U,V)^{k-1}$: 
    \begin{equation*}
        \begin{split}
            (\alpha_0U^{k-1}+\alpha_1U^{k-2}V+\cdots+\alpha_{k-1}V^{k-1})y_0+(\beta_0U^{k-1}+\cdots+\beta_{k-1}V^{k-1})y_1+\\
            +(\gamma_0U^{k-1}+\cdots+\gamma_{k-1}V^{k-1})y_2.
        \end{split}
    \end{equation*}
    It belongs to $\Ann_S(f)_k$ if and only if the coefficients $\alpha_s, \beta_s,$ and $\gamma_s$ solve the equation
    \begin{equation*}
        \begin{split}
            \alpha_0p_{0,u^{k-1}}+\alpha_1p_{0,u^{k-2}v}+\cdots+\alpha_{k-1}p_{0,v^{k-1}}+\beta_0p_{1,u^{k-1}}+\cdots\\
            \cdots+\gamma_0p_{2,u^{k-1}}+\cdots +\gamma _{k-1}p_{2,v^{k-1}}=0.
        \end{split}
    \end{equation*}
    The partial derivatives of $p_0, p_1, p_2$, appearing in the above expression, are all homogeneous of degree $d-k$; bringing together all the coefficients of the same monomial, and setting it zero, we can construct an equivalent homogeneous linear system of $d-k+1$ equations, one for each monomial, in the $3k$ unknowns
    $\alpha_0,\ldots,\alpha_{k-1}, \beta_0, \ldots, \beta_{k-1}, \gamma_0, \ldots, \gamma_{k-1}$. The matrix of the linear system is exactly $M_k$, therefore $m_k=3k-\rank M_k$.
    
    \pause
    We now compute $n_k$ considering a general polynomial of degree $k$ in $U,V$:
    $$\delta_0U^k+\delta_1U^{k-1}V+\cdots+\delta_kV^k.$$
    This element belongs to $\Ann_S(f)_k$ if, as before, the $\delta_s$ are solution of the following equation
    \begin{equation*}
        \begin{split}
            (\delta_0p_{0,u^k}+\delta_1p_{0,u^{k-1}v}+\cdots+\delta_kp_{0,v^k})x_0+(\delta_0p_{1,u^k}+\cdots+\delta_kp_{1,v^k})x_1+\\
            +(\delta_0p_{2,u^k}+\cdots+\delta_kp_{2,v^k})x_2+\delta_0g_{u^k}+\cdots+\delta_kg_{v^k}=0.
        \end{split}
    \end{equation*}
    Bringing out the coefficients of $x_0, x_1, x_2$ and also the coefficients of the monomials in $u,v$ which can be of degree $d-k-1$ or $d-k$, we again get a homogeneous linear system. We have $3(d-k)+(d-k+1)$ equations one for every possible monomial and $k+1$ unknowns. The related matrix is $N'_k$. We then deduce that $n_k=k+1-\rank N'_k$.
    \end{proof}
    
    \begin{example}
        To better understand the role of the matrices $A_k$, $B_k$, $C_k$, and $G_k$ in the rank of the matrices $M_k$ and $N_k'$ we now give an explicit example.
        
        \pause
        Let $f=5u^3v^2x_0+(u^5+v^5)x_1+(2uv^4-3u^3v^2)x_2+u^6-3u^2v^4$ be a Perazzo $3$-fold of degree $d=6$, and let $A=S/\Ann_S(f)$ the related Gorenstein algebra. We are interested to compute the values of $A$'s Hilbert vector $h$. We start by computing the matrices $M_2,M_3,N_2'$ and $N_3'$:
        \begin{equation*}
            M_2=\left(\begin{array}{cc|cc|cc}
                0 & 0 & 1 & 0 & 0 & 0\\
                0 & \sfrac{1}{2} & 0 & 0 & 0 & -\sfrac{1}{10}\\
                \sfrac{1}{2} & 0 & 0 & 0 & -\sfrac{1}{10} & 0\\
                0 & 0 & 0 & 0 & 0 & \sfrac{2}{5}\\
                0 & 0 & 0 & 1 & \sfrac{2}{5} & 0
            \end{array}\right),
        \end{equation*}
        \begin{equation*}
            M_3=\left(\begin{array}{ccc|ccc|ccc}
                0 & 0 & \sfrac{1}{2} & 1 & 0 & 0 & 0 & 0 & -\sfrac{1}{10}\\
                0 & \sfrac{1}{2} & 0 & 0 & 0 & 0 & 0 & -\sfrac{1}{10} & 0\\
                \sfrac{1}{2} & 0 & 0 & 0 & 0 & 0 & -\sfrac{1}{10} & 0 & \sfrac{2}{5}\\
                0 & 0 & 0 & 0 & 0 & 1 & 0 & \sfrac{2}{5} & 0
            \end{array}\right),
        \end{equation*}
        \begin{equation*}
            N_2'=\begin{pmatrix}
                0 & 0 & \sfrac{1}{2}\\
                0 & \sfrac{1}{2} & 0\\
                \sfrac{1}{2} & 0 & 0\\
                0 & 0 & 0\\
                \hline
                1 & 0 & 0\\
                0 & 0 & 0\\
                0 & 0 & 0\\
                0 & 0 & 1\\
                \hline
                0 & 0 & -\sfrac{1}{10}\\
                0 & -\sfrac{1}{10} & 0\\
                -\sfrac{1}{10} & 0 & \sfrac{2}{5}\\
                0 & \sfrac{2}{5} & 0\\
                \hline
                1 & 0 & 0\\
                0 & 0 & 0\\
                0 & 0 & -\sfrac{1}{5}\\
                0 & -\sfrac{1}{5} & 0\\
                -\sfrac{1}{5} & 0 & 0
            \end{pmatrix}\text{, and }
            N_3'=\begin{pmatrix}
                0 & 0 & \sfrac{1}{2} & 0\\
                0 & \sfrac{1}{2} & 0 & 0\\
                \sfrac{1}{2} & 0 & 0 & 0\\
                \hline
                1 & 0 & 0 & 0\\
                0 & 0 & 0 & 0\\
                0 & 0 & 0 & 1\\
                \hline
                0 & 0 & -\sfrac{1}{10} & 0\\
                0 & -\sfrac{1}{10} & 0 & \sfrac{2}{5}\\
                -\sfrac{1}{10} & 0 & \sfrac{2}{5} & 0\\
                \hline
                1 & 0 & 0 & 0\\
                0 & 0 & 0 & -\sfrac{1}{5}\\
                0 & 0 & -\sfrac{1}{5} & 0\\
                0 & -\sfrac{1}{5} & 0 & 0
            \end{pmatrix}.
        \end{equation*}
        By Proposition \ref{hilbert function}, we obtain
        \begin{equation*}
            \begin{split}
                h_0=h_6&=1,\\
                h_1=h_5&=5,\\
                h_2=h_4&=\rank M_2+\rank N'_2=5+3=8,\\
                h_3&=\rank M_3+\rank N'_3=4+4=8.
            \end{split}
        \end{equation*}
        Therefore the Hilbert vector of $A$ is $h=(1,\,5,\,8,\,8,\,8,\,5,\,1)$.
    \end{example}
    
    We have seen in Section \ref{catal} that the nature of a form $h\in K[u,v]_t$ is connected to the position of $[h]$ with respect to the rational normal curve $C_t$ and its secant varieties. Our idea is to consider the projective plane $\pi$ generated by $p_0,p_1$ and $p_2$ trying to understand its position with respect to that curve and its secant varieties. In fact, due to Lemma \ref{generators}, the quotient algebra $S/\Ann_S(f)$ does not depend on the choice of a triplet of generators of $\pi$. So, considering a more suitable triplet $q_0,q_1,q_2$, we obtain an equivalent and explicit Perazzo hypersurface that we can work with. The following Theorem \ref{lowcharacterisation} is a more precise statement.
    
    \begin{lemma}\label{generators} Let $f_1=p_0(u,v)x_0+p_1(u,v)x_1+p_2(u,v)x_2$ and $f_2=q_0(u,v)x_0+q_1(u,v)x_1+q_2(u,v)x_2$ be two Perazzo 3-folds of degree $d$ in $\PP^4$ such that $\langle p_0,p_1,p_2 \rangle = \langle q_0,q_1,q_2 \rangle \subset K[u,v]_{d-1} $.  Then, $S/\Ann _S (f_1)$ is isomorphic to $S/\Ann _S (f_2)$.
    \end{lemma}
    \begin{proof} By Theorem \ref{gorensteinannullatore}, it is enough to prove that $f_1$ and $f_2$ are equal up to a linear change of variables. By assumption, we can write $q_0(u,v)=\lambda _0p_0(u,v)+\lambda _1p_1(u,v)+\lambda _2p_2(u,v)$, $q_1(u,v)=\mu _0p_0(u,v)+\mu _1p_1(u,v)+\mu _2p_2(u,v)$, $q_3(u,v)=\nu _0p_0(u,v)+\nu _1p_1(u,v)+\nu _2p_2(u,v)$. We have 
    $$\begin{array}{rcl} f_2 & = & q_0(u,v)x_0+q_1(u,v)x_1+q_2(u,v)x_2 \\
    & = & (\lambda _0p_0(u,v)+\lambda _1p_1(u,v)+\lambda _2p_2(u,v))x_0+\\
    & & +(\mu _0p_0(u,v)+\mu _1p_1(u,v)+\mu _2p_2(u,v))x_1+ \\
    & & +(\nu _0p_0(u,v)+\nu _1p_1(u,v)+\nu _2p_2(u,v))x_2 \\ & = & (\lambda _0 x_0+\mu _0x_1+ \nu _0x_2)p_0(u,v)+ (\lambda _1 x_0+\mu _1x_1+ \nu _1x_2 ) p_1(u,v)+\\
    & & +  (\lambda _2 x_0+\mu _2x_1+ \nu _2 x_2 ) p_2(u,v).
    \end{array}
    $$
    Therefore, $f_1$ and $f_2$ differ from a linear change of variables and the proof is completed.
    \end{proof}
    
    As an important fact, we now prove that the matrices $M_k$ and $N_k'$ have rank greater or equal to $3$. The important assumption is that the polynomials $p_i$ are linearly independent. To do that we give an explicit characterisation of the Perazzo $3$-folds when these ranks are smaller or equal to $3$.
    
    \begin{theorem}\label{lowcharacterisation}
        Let $X=V(f)\subseteq\PP^4$ be a Perazzo $3$-fold of degree $d\ge 4$. If $\rank M_k\le 3$ and $\rank N_k'\le 3$, then, up to a change of variables,  one of the following cases occurs:
        \begin{itemize}
            \item[(i)]   $f(x_0,x_1,x_2,u,v)=u^{d-1}x_0+u^{d-2}vx_1+u^{d-3}v^2x_2+ au^d+bu^{d-1}v+cu^{d-2}v^2$ with $a,b,c\in K$, or
            \item[(ii)]  $f(x_0,x_1,x_2,u,v)=u^{d-1}x_0+u^{d-2}vx_1+v^{d-1}x_2+au^d+bu^{d-1}v+cv^{d}$ with $a,b,c\in K$, or
            \item[(iii)]  $f(x_0,x_1,x_2,u,v)=u^{d-1}x_0+(\lambda u+\mu v)^{d-1}x_1+v^{d-1}x_2+au^d+b(\lambda u+ \mu v)^d+cv^d$ with $\lambda, \mu\in K^*$ and $a,b,c\in K$.
        \end{itemize}
        
        In particular, we see that in all these cases $\rank M_k,N_k'=3$.
    \end{theorem}
    \begin{proof}
        Let $f=x_0p_0+x_1p_1+x_2p_2+g$. Being $\rank M_k,N_k'\le 3$, the same bound is true also for the rank of the matrices $A_k,B_k,C_k$ and $G_k$. We use Proposition \ref{rank} to compute the polynomials $p_0,p_1$ and $p_2$ in the various possibilities. Then we will see what we can say about the form $g$.
        
        (I) $d\geq 7$ and $\rank A_k=\rank B_k=\rank C_k=3$ for $4\leq k\leq \lfloor\frac{d+1}{2}\rfloor $ and $p_0, p_1, p_2$ have all Waring rank $3$. By Proposition \ref{lemmarank4} and being $\rank M_k\le 3$, the spaces of the columns of $A_k, B_k, C_k$ coincide, so there exist distinct linear forms $l_1, l_2, l_3$ and suitable constants such that 
        $$p_0=\lambda_0l_1^{d-1}+\mu_0l_2^{d-1}+\nu_0l_3^{d-1}$$
        $$p_1=\lambda_1l_1^{d-1}+\mu_1l_2^{d-1}+\nu_1l_3^{d-1}$$
        $$p_2=\lambda_2l_1^{d-1}+\mu_2l_2^{d-1}+\nu_2l_3^{d-1}.$$
        Since $p_0, p_1, p_2$ are linearly independent, the matrix
        $\begin{pmatrix}
            \lambda_0 & \mu_0 & \nu_0\\
            \lambda_1 & \mu_1 & \nu_1\\
            \lambda_2 & \mu_2 & \nu 2
        \end{pmatrix}$
        is invertible, therefore $\langle p_0, p_1, p_2\rangle = \langle l_0^{d-1}, l_1^{d-1}, l_2^{d-1}\rangle$. In view of Lemma \ref{generators} in $f$ we can replace $p_0, p_1, p_2$ with $l_0^{d-1}, l_1^{d-1}, l_2^{d-1}$.
        
        (II) $d\geq 7$, $\rank A_k=3$ for $4\leq k\leq \lfloor\frac{d+1}{2}\rfloor$, but $p_0$ has Waring rank strictly greater than $3$. From Proposition \ref{rank}, it follows that, up to a change of variables, $p_0=u^{d-1}+\alpha uv^{d-2}+\beta v^{d-1}$, for suitable constants $\alpha, \beta$, $\alpha\neq 0$. Then
        \[ M_3=
        \begin{pmatrix}
         1 & 0 & 0 & b_0 & b_1 & b_2 & c_0 & c_1 & c_2\\
         0 & 0 & 0 &b_1&b_2&b_3 &c_1&c_2&c_3\\
         \vdots & \vdots& \vdots& \vdots& \vdots& \vdots&\vdots&\vdots&\vdots\\
         0 & 0 & \alpha & b_{d-4} & b_{d-3} & b_{d-2} & c_{d-4} & c_{d-3} & c_{d-2} \\
         0 & \alpha & \beta & b_{d-3} & b_{d-2} & b_{d-1} & c_{d-3} & c_{d-2} & c_{d-1}
        \end{pmatrix}.
        \]
        From $\rank M_3<4$ it follows $b_1=\cdots=b_{d-3}=c_1=\cdots=c_{d-3}=0$. Therefore
        $$p_1=b_0u^{d-1}+b_{d-2}uv^{d-2}+b_{d-1}v^{d-1},\  p_2=c_0u^{d-1}+c_{d-2}uv^{d-2}+c_{d-1}v^{d-1},$$
        so  we can replace $p_0, p_1, p_2$ with $u^{d-1}, uv^{d-2}, v^{d-1}$.
        
        (III) $\rank A_k=2$ for $3\leq k\leq \lfloor\frac{d+1}{2}\rfloor $ and $p_0$ has Waring rank $2$, so it can be written $p_0=u^{d-1}+v^{d-1}$.
        Then $M_3$ is as in case (II) with $\alpha=0$, $\beta=1$ 
        and 
        \[
        M_2=
        \begin{pmatrix}
         1 & 0 &b_0 & b_1 & c_0 & c_1\\
         0&0&b_1&b_2&c_1&c_2\\
         \vdots & \vdots & \vdots & \vdots & \vdots & \vdots\\
         0& 1 & b_{d-2} & b_{d-1}& c_{d-2} & c_{d-1}
        \end{pmatrix}.
        \]
        From $\rank M_2<4$ we deduce that 
        \[\rank \begin{pmatrix}
         b_1 & b_2\\
         \vdots & \vdots\\
         b_{d-3} & b_{d-2}
        \end{pmatrix}<2, \ 
        \rank \begin{pmatrix}
         c_1 & c_2\\
         \vdots & \vdots\\
         c_{d-3} & c_{d-2}
        \end{pmatrix}<2, \ 
        \rank \begin{pmatrix}
         b_1 & b_2 & \ldots & b_{d-2}\\
         c_1 & c_2 & \ldots & c_{d-2}
        \end{pmatrix}<2.
        \]
        Therefore $$(b_1, \ldots, b_{d-2})=(\lambda^{d-3}, \lambda^{d-4}\mu, \ldots,\mu^{d-3}), \ (c_1, \ldots, c_{d-2})=(\sigma^{d-3}, \sigma^{d-4}\rho, \ldots,\rho^{d-3}),$$ for suitable $\lambda, \mu, \sigma, \rho\in K$. We get:
        $$p_1=b_0u^{d-1}+uv((d-1)\lambda^{d-3}u^{d-3}+{\binom{d-1}{2}}\lambda^{d-4}\mu u^{d-4}v+\cdots)+b_{d-1}v^{d-1},$$
        $$p_2=c_0u^{d-1}+uv((d-1)\sigma^{d-3}u^{d-3}+{\binom{d-1}{2}}\sigma^{d-4}\rho u^{d-4}v+\cdots)+c_{d-1}v^{d-1}.$$
        We can also write
        $$p_1=b_0u^{d-1}+uv\phi_{d-3}+b_{d-1}v^{d-1}, \ 
        p_2=c_0u^{d-1}+kuv\phi_{d-3}+c_{d-1}v^{d-1}$$ where $\phi_{d-3}$ is a form of degree $d-3$ and $k\in K$, because $(b_1, \ldots,b_{d-2})$ and $(c_1, \ldots, c_{d-2})$ are proportional. We  can assume $b_0=c_0=0$ and we get $v^{d-1}\in \langle p_0, p_1, p_2\rangle$, hence $u^{d-1}, uv\phi_{d-3}\in \langle p_0, p_1, p_2\rangle$. Finally, adding to $uv\phi_{d-3}$ suitable multiples of $u^{d-1}, v^{d-1}$, we get $(\lambda u+\mu v)^{d-1}\in \langle p_0, p_1, p_2\rangle$.

        (IV) $\rank A_k=2$ for $3\leq k\leq \lfloor\frac{d+1}{2}\rfloor $ but $p_0$ has Waring rank $>2$, so up to a change of variables $p_0=u^{d-2}v$.
        \[ M_2=
        \begin{pmatrix}
        0&1&b_0&b_1&c_0&c_1\\
        1&0&b_1&b_2&c_1&c_2\\
        \vdots & \vdots &\vdots&\vdots&\vdots&\vdots\\
        0 & 0 & b_{d-2}& b_{d-1} & c_{d-2} & c_{d-1}
        \end{pmatrix}
        \]
        has rank less or equal than $3$, therefore
        \[
        \rank \begin{pmatrix}
        b_2 & b_3 & c_2 & c_3\\
        \vdots & \vdots & \vdots & \vdots\\
        b_{d-2} & b_{d-1} & c_{d-2} & c_{d-1}
        \end{pmatrix}<2,
        \]
        and arguing in a similar way to (III), we conclude that $\langle p_0, p_1, p_2\rangle$ is either of the form $\langle u^{d-1}, u^{d-2}v, (\lambda u+\mu v)^{d-1}\rangle $, or $\langle u^{d-1}, u^{d-2} v, u^{d-3}v^2\rangle.$
        
        (V) $\rank A_k=\rank B_k=\rank C_k=1$, then $p_0, p_1, p_2$ are all pure powers of degree $d-1$.
        
        (VI)  Let $\pi$ be the $2$-plane generated by the polynomials $p_0, p_1, p_2$. 
        
        \begin{enumerate}
            \item[$d=4$] $\pi\subset \PP^3=\PP(K[u,v]_3)$. Then $\pi$ and $C_3$ have to intersect in three points counted with multiplicity, so we are in the case (I), or (II), or (IV);
            \item[$d=5$] $\pi\subset \PP^4=\PP(K[u,v]_4)$. The tangential variety $TC_4$ has codimension $2$, so the intersection $\pi\cap TC_4\neq \emptyset$. If $\pi$ intersects $TC_4$ outside its singular locus $C_4$, up to a change of variables $u^3v\in\pi$ and we conclude as in (IV); otherwise, we are in the situation of (V);
            \item[$d=6$] $\pi\subset \PP^5=\PP(K[u,v]_5)$. Now $\sigma_2(C_5)$ has codimension $2$ and therefore $\pi\cap\sigma_2(C_5)\neq\emptyset$. Therefore we are either in the situation of (III) or of (IV).
        \end{enumerate}
        
        We have proved that for any $d\geq 4$, if the matrices $M_k$ and $N_k'$ associated to $f$ have all rank $\le 3$, then the polynomials $p_0, p_1, p_2$ are as in (i), or (ii), or (iii).
        
        \pause
        A computation proves that in those cases the rank of $M_k$ is always equal to $3$. Being $\rank M_k=\rank N_{d-k}$ and being $N'_k$ an extension of $N_k$, we obtain also that $\rank N_k'=3$.
        
        \pause
        It remains to find out how we can choose the polynomial $g$ in each of the cases. From Proposition \ref{hilbert function} we deduce that the only condition that $g$ has to satisfy is $\Ann_S(p_0x_0+p_1x_1+p_2x_2)_3=\Ann_S(f-g)_3=\Ann_S(f)_3$. In other words, we impose that $g$ is annihilated by a system of generators of $\Ann_S(f)_3$.
        
        (i) If $f(x_0,x_1,x_2,u,v)=u^{d-1}x_0+u^{d-2}vx_1+u^{d-3}v^2x_2+g$, we have that $\Ann_S(f)_3=\langle V^3\rangle$ and so $g=g_0u^d+g_1u^{d-1}v+g_2u^{d-2}v^2$. 
        
        (ii) If $f(x_0,x_1,x_2,u,v)=u^{d-1}x_0+u^{d-2}vx_1+v^{d-1}x_2+g$, we have that $\Ann_S(f)_3=\langle UV^2\rangle$. This gives that $\sum_{i=2}^{d-1}g_i\binom{d}{i}(k-i)i(i-1)u^{k-i-1}v^{i-2}=0$, so $g_2=\dots=g_{d-1}=0$. Thus we get $g=g_0u^d+g_1u^{d-1}v+g_dv^d$.
        
        (iii) If $f(x_0,x_1,x_2,u,v)=u^{d-1}x_0+(\lambda u+\mu v)^{d-1}x_1+v^{d-1}x_2+g$, we have that $\Ann_S(f)_3=\langle \mu U^2V-\lambda UV^2\rangle$. Then we have the condition $$\sum_{i=2}^{d-1}\binom{k-3}{i-1}(\mu g_i-\lambda g_{i+1})u^{d-k-2}v^{i-1}=0 \iff \mu g_i-\lambda g_{i+1}=0,\quad i=1,\dots,k-2.$$
        So we can collect $g_1$ and complete the $d$-th power to obtain $g=au^d+b(\lambda u+\mu v)^d+cv^d$.
    \end{proof}
    
    By Proposition \ref{hilbert function}, we have to look at the matrices $M_k$ and $N_k'$ to fully understand the shape of the Hilbert vector, in particular looking at their rank. As a consequence of Proposition \ref{hilbert function}, in Propositions \ref{upper} and \ref{lower} we will determine the maximal and minimal $h$-vector \footnote{Given two tuple $h=(h_0,h_1,\dots,h_d)$ and $h'=(h'_0,h'_1,\dots,h'_d)$, we say that $h\le h'$ if $h_k\le h'_k$ for every $k$, $0\le k\le \frac{d}{2}$.} that a Gorenstein algebra, associated to a Perazzo $3$-fold, can reach. In both Propositions \ref{upper} and \ref{lower} the key idea is to, respectively, maximise and minimise the rank of the matrices $M_k$ and $N'_k$.
    
    \begin{proposition}\label{upper}  Let $d\ge 4$.
    The maximum $h$-vector of  the Artinian Gorenstein algebras $S/\Ann_S(f)$ associated to $f$, a Perazzo $3$-folds of degree $d$ in $\PP^4$, is:
        \begin{itemize}
        \item[(1)] If  $d=4t+0$  
        then 
            $$h_k= \begin{cases} 4k+1 &\text{ for } 0\le k \le t \\
            d+2 &\text{ for } t+1\le k \le 2t \\
            \text{symmetry.}
             \end{cases}$$
        \item[(2)] If $d=4t+1$ 
        then 
        $$h_k \begin{cases} 4k+1 &\text{ for } 0\le k\le t \\
        d+2 &\text{ for } t+1\le k \le 2t \\
        \text{symmetry.}
         \end{cases}
        $$
        \item[(3)] If $d=4t+2$   
        then 
            $$h_k= \begin{cases} 4k+1 &\text{ for } 0\le k\le t \\
            d+2  &\text{ for } t+1\le k \le 2t+1 \\
            \text{symmetry.}
             \end{cases}$$
        \item[(4)]  If $d=4t+3$ 
        then 
            $$h_k= \begin{cases} 4k+1 &\text{ for } 0\le i\le t+1 \\
            d+2 &\text{ for } t+2\le k \le 2t+1 \\
            \text{symmetry.}
             \end{cases}$$
        \end{itemize}
        \end{proposition}
        
    \begin{example}
        To better understand how the maximal Hilbert vector behaves, we give a bunch of examples in the first degrees.
        \begin{equation*}
            \arraycolsep=2pt\def\arraystretch{1.5}
            \begin{array}{c@{\hspace{1,5cm}}l}
                d=6=4\cdot 1+2 & h=(1,5,8,8,8,5,1) \\
                d=7=4\cdot 1+3 & h=(1,5,9,9,9,9,5,1) \\
                d=8=4\cdot 2+0 & h=(1,5,9,10,10,10,9,5,1) \\
                d=9=4\cdot 2+1 & h=(1,5,9,11,11,11,11,9,5,1) \\
                d=10=4\cdot 2+2 & h=(1,5,9,12,12,12,12,12,9,5,1) \\
                d=11=4\cdot 2+3 & h=(1,5,9,13,13,13,13,13,13,9,5,1)\\
                d=12=4\cdot 3+0 & h=(1,5,9,13,14,14,14,14,14,13,9,5,1)\\
                \vdots & \vdots
            \end{array}
        \end{equation*}
        
    \end{example}
    
    \begin{proof}[Proof (Theorem \ref{upper})]
        Let $f$ be a form that defines a Perazzo $3$-fold and let $S/\Ann_S(f)$ be the associated Artinian graded Gorenstein algebra. By the symmetry of the $h$-vector, we have to compute only the first half, i.e. $h_k$ for $0\le k\le \frac{d}{2}$.
        
        \pause
        By virtue of Proposition \ref{hilbert function}, the maximal $h$-vector can be reached only when the integers $m_k$ and $n_k$ are as small as possible. That is to say, the ranks of the matrices $M_k$ and $N'_k$ are as large as possible, possibly we want them to be full rank.
        
        \pause
        Clearly $\rank M_k\leq \min\{3k, d-k+1\}.$ Therefore 
        \begin{equation*}
            \rank M_k\leq \begin{cases} 3k \,& \text{for} \  k\leq {\frac{d+1}{4}};\\
            d-k+1 \,& \text{for} \ k\geq {\frac{d+1}{4}}.
            \end{cases}
        \end{equation*}
        
        Regarding $N'_k$, we observe that, under the assumption $k\leq\frac{d}{2}$,  $k+1\leq 3(d-k)+(d-k+1)$, so always $\rank N'_k\leq k+1$. Thus we obtain an upper bound on every possible value $h_k$ as $h_k\leq k+1+\min\{3k,d-k+1\}$, i.e.
        \begin{equation*}
            h_k\leq \begin{cases} 4k+1 \,& \text{for} \  k\leq {\frac{d+1}{4}};\\
            d+2 \,& \text{for} \ k\geq {\frac{d+1}{4}}.
            \end{cases}
        \end{equation*}
        
        Distinguishing all possible degrees $d$ modulo $4$, we obtain the formulas stated above. The last fact to prove is that this upper bound is indeed achieved. The following Example \ref{esempiodellavita} will show that in every degree $d$, an example of a Perazzo $3$-fold with maximal Hilbert vector exists.
    \end{proof}
    
    \begin{example}[Maximal Hilbert vector]\label{esempiodellavita}
        For any integer $d\ge 4$, we write $d=3r+\epsilon $ with $0\le \epsilon \le 2$. We take 
        $$\begin{array}{rcl}
            p_0(u,v) & = & \sum _{i=0}^{r}\binom{d-1}{i}\, \frac{1}{1+i}\,u^{d-1-i}v^i,  \\
            p_1(u,v) & = & \sum _{i=r}^{2r-1+\epsilon }\binom{d-1}{i}\, \frac{1}{1+i}\,u^{d-1-i}v^i, \\
            p_2(u,v) & = & \sum _{2r-1+\epsilon }^{d-1}\binom{d-1}{i}\, \frac{1}{1+i}\,u^{d-1-i}v^i, \text{ and } \\
            g(u,v) & = & 0.
        \end{array}$$
        The polynomials $p_i$ are linearly independent because at least one element of the monomial basis of $K[u,v]_{d-1}$ appears once and only once in at least two forms. They are also algebraically dependent due to Remark \ref{importantremark}. We now compute the catalecticant matrices associated to the homogeneous polynomials $p_0,p_1$, and $p_2$.
        
        \pause
        We start computing the catalecticant matrices $A_k$ of $p_0$. The matrix in Figure \ref{matriceuno} occurs whenever $k\le r$; while, the matrix in Figure \ref{matricedue} occurs when $r+1\le k\le \frac{d}{2}$. We see that these matrices are caracterised by a non-zero triangular section.

        \begin{figure}[H]
            \centering
            \includegraphics[scale=0.2]{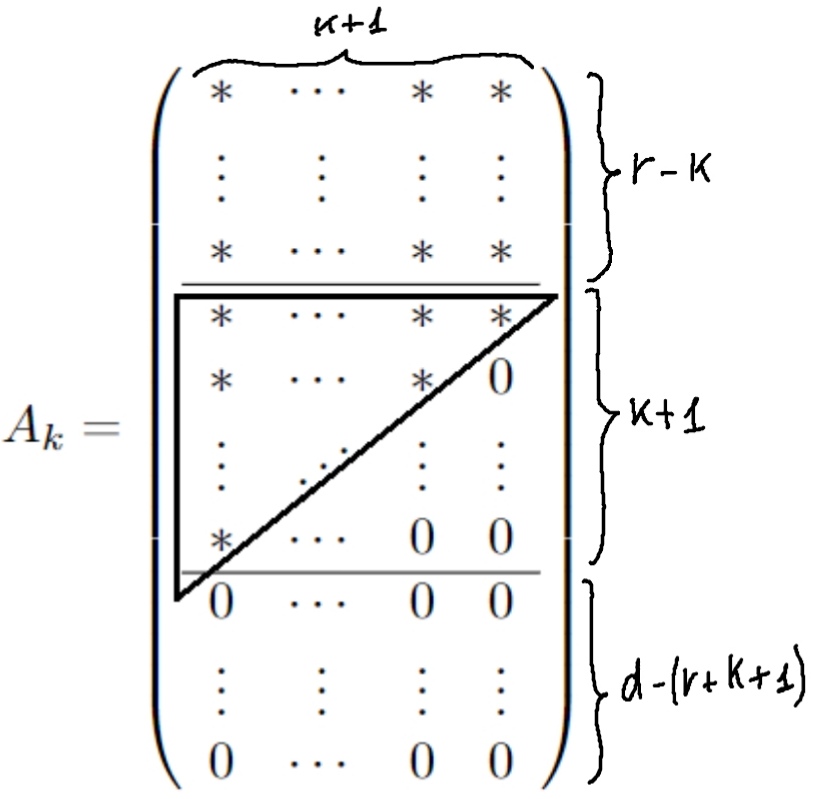}
            \caption{$A_k$ for $k\le r$}
            \label{matriceuno}
        \end{figure}
        
        \begin{figure}[H]
            \centering
            \includegraphics[scale=0.21]{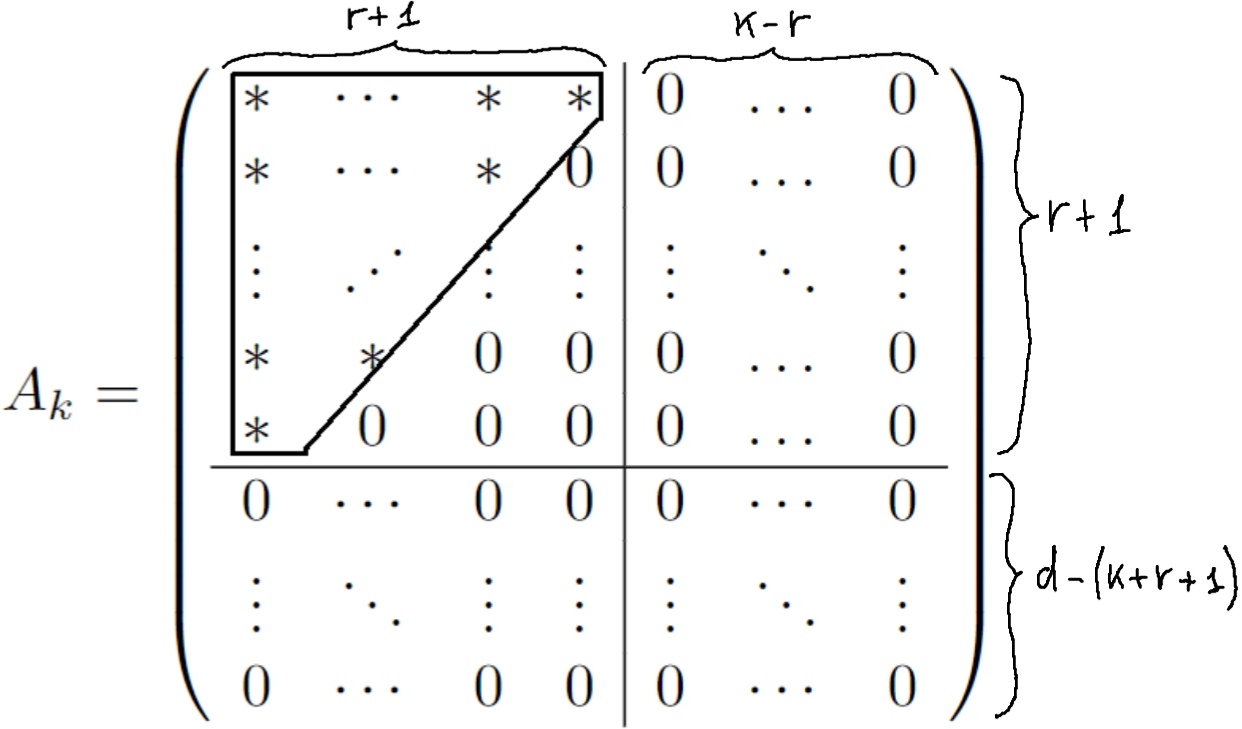}
            \caption{$A_k$ for $k\ge r+1$}
            \label{matricedue}
        \end{figure}
        
        The catalecticant matrices $C_k$ of $p_2$ have a symmetric behaviour with respect to the matrices $A_k$. Finally, for $p_1$, the catalecticant matrices $B_k$ are characterised by a skew-diagonal band which is non-zero. The following four matrices represent the case $k\le r+\epsilon-1$: the matrix in Figure \ref{matricetre} represents the case $k<r+\epsilon-1$, while the three matrices in Figures \ref{matricequattro}, \ref{matriceconque}, and \ref{matricesei} represent the case $k=r+\epsilon-1$ varying the index $\epsilon$.
        
        \begin{figure}[H]
            \begin{minipage}[c]{0.45\linewidth}
                \includegraphics[width=\linewidth]{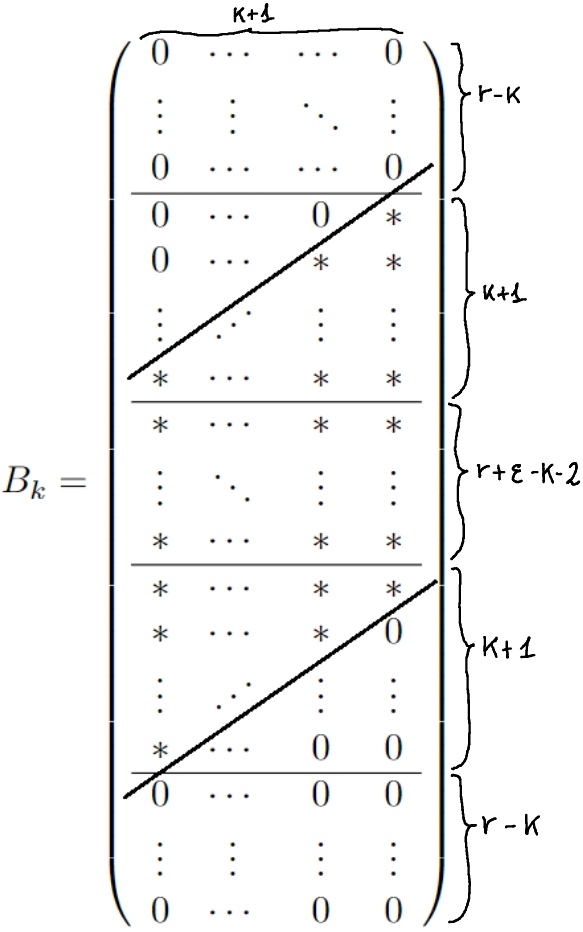}
                \caption{$B_k$ for $k<r+\epsilon-1$}
                \label{matricetre}
            \end{minipage}
            \hfill
            \begin{minipage}[c]{0.45\linewidth}
                \includegraphics[width=\linewidth]{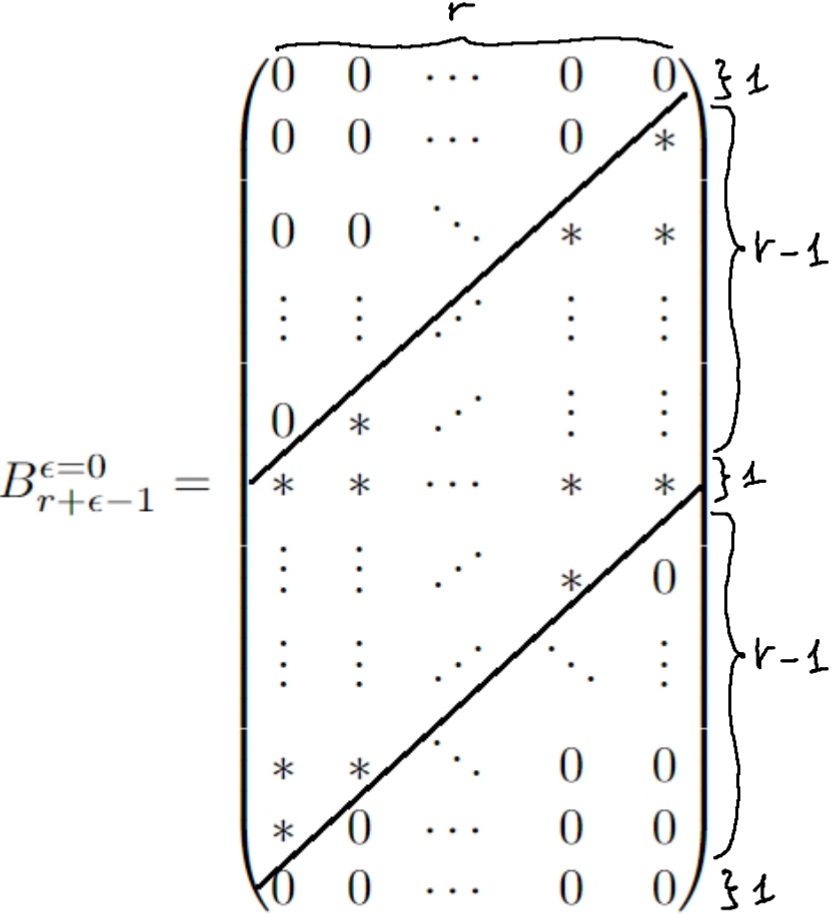}
                \caption{$B_{r+\epsilon-1}$ for $\epsilon=0$}
                \label{matricequattro}
            \end{minipage}%
        \end{figure}
        
        \begin{figure}[H]
            \begin{minipage}[c]{0.49\linewidth}
                \includegraphics[width=\linewidth]{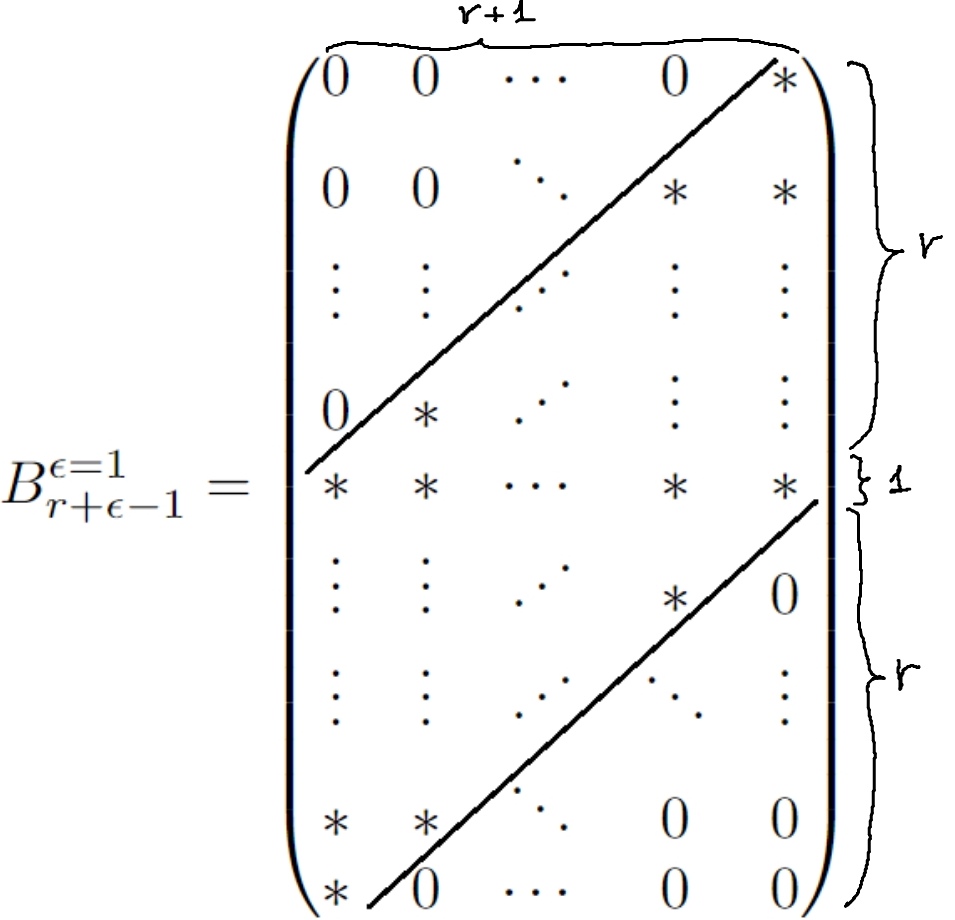}
                \caption{$B_{r+\epsilon-1}$ for $\epsilon=1$}
                \label{matriceconque}
            \end{minipage}
            \hfill
            \begin{minipage}[c]{0.49\linewidth}
                \includegraphics[width=\linewidth]{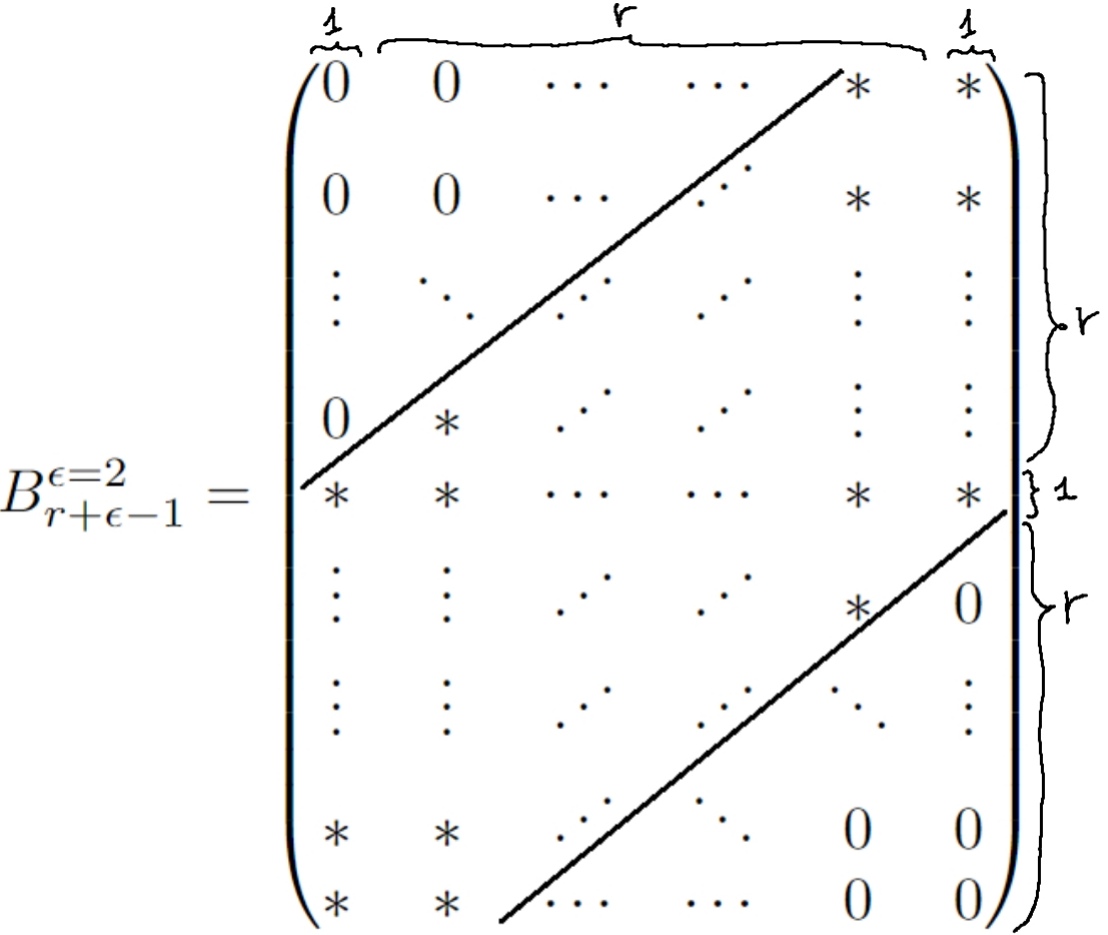}
                \caption{$B_{r+\epsilon-1}$ for $\epsilon=2$}
                \label{matricesei}
            \end{minipage}%
        \end{figure}

        Figure \ref{matricesette} represents the matrix $B_k$ in the case $k\ge r+\epsilon$.
        \begin{figure}[H]
    	    \centering
            \includegraphics[width=\linewidth]{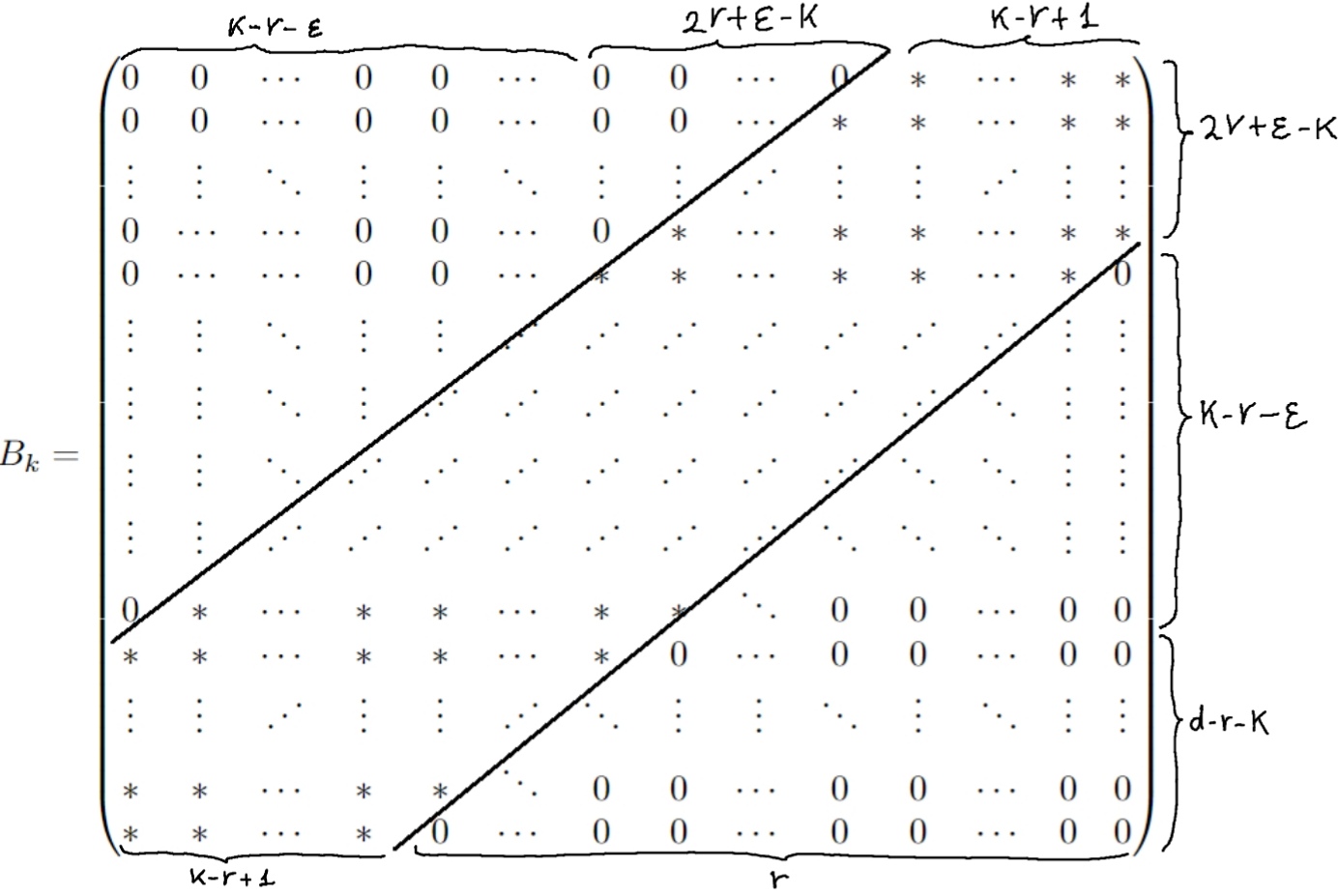}
            \caption{$B_k$ for $k\ge r+\epsilon$}
            \label{matricesette}
        \end{figure}

        We now prove that the matrices $M_k$, and $N_k'$ have all full rank for $k\ge 2$. In particular, we see that
        \begin{equation*}
            \rank N_k'=k+1\text{, and }
            \rank M_k=\begin{cases}
                3k &\text{if } k\le \frac{3}{4}r+\frac{\epsilon-1}{4}\\
                d-k+1 &\text{if } k\ge \frac{3}{4}r+\frac{\epsilon-1}{4}
            \end{cases}.
        \end{equation*}
        
        We remark that these matrices have been defined as follows:
        \begin{equation*}
            M_k=
            \begin{pmatrix}
                A_{k-1} \vert B_{k-1} \vert C_{k-1}
            \end{pmatrix}
            \quad
            N'_k= \begin{pmatrix}
                A_{k} \\ \hline B_{k} \\ \hline C_{k} \\ \hline G_{k}
            \end{pmatrix}.
        \end{equation*}

        The matrices $N_k'$ can be treated as follows:
        \begin{itemize}
            \item if $2\le k\le r$ the matrix $A_k$ has full rank, i.e. it has rank $k+1$. Thus, $N'_k$ has rank $k+1$;
            \item if $r+1\le k\le\frac{d}{2}$, then we can consider the sub-matrix of $N'_k$ given by the first $r+1$ rows and the last $k-r$ rows. This square sub-matrix has non-zero determinant, therefore $N'_k$ has full rank.
        \end{itemize}
        
        Instead, for the matrices $M_k$ we have to proceed as follows: when $k\le \frac{d+1}{4}$ we have to prove that $M_k$'s columns are linearly independent, and, when $k\ge \frac{d+1}{4}$, we have to prove that its rows are linearly independent.
        
        \pause
        First of all, if $k-1\le \frac{d+1}{4}-1=\frac{3}{4}r-\frac{3-\epsilon}{4}$ then $M_k'$ has the desired rank. Indeed, since $k-1<r$, the matrices $A_{k-1},B_{k-1}$ and $C_{k-1}$ have rank $k$. Moreover, due to their shape, when we impose to a linear combination of the columns of $M_k$ to be zero, this is equivalent to a system of three linear combinations one for each matrix $A_{k-1}, B_{k-1}, C_{k-1}$. So, the matrix $M_k$ has rank $3k$.
        
        \pause
        We consider now the cases $\frac{d+1}{4}\le k \le \frac{d}{2}$: we proceed distinguish the various values that the integer $\epsilon$ can take.
        \begin{itemize}
            \item[$\epsilon=0$ :] we suppose for the moment that $r\ge 5$. Then we have the following possibilities:
            \begin{enumerate}
                \item if $\frac{3}{4}r-\frac{3}{4}\le k-1 \le r-2$ then we divide the rows of $M_k$ in five groups with the following number of rows $(r-k,k+1,r-k-2,k+1,r-k)$. These sections have all full rank: the first, the third, and the fifth contains each a sub-matrix of some Hilbert matrix\footnote{The Hilbert matrix $H$ of order $n$ is a $n\times n$ square matrix defined as $H_{ij}=\frac{1}{i+j+1}$. The Hilbert matrix is totally positive, i.e. all minors are strictly positive \cite[Theorem A]{fiedler}.} which has full rank. The second and the fifth contain a triangular matrix with non zero diagonal, so they have full rank. Reasoning as the previous point we obtain our claim.
                \item The case $k-1=r-1$ is similar to the previous point dividing the matrix $M_k$ in two parts with $r+1$ and $r$ rows respectively.
                \item In the final case $r\le k-1\le \frac{d}{2}$, we divide $M_k$ in two parts with $r+1,2r-k$ rows respectively. Again they are full rank and, combined together, they form $M_k$ which turns out to be full rank.
            \end{enumerate}
            If $1\le r\le 4$ then the condition $\frac{3}{4}r-\frac{3}{4}\le r-2$ is no longer true; but, since we are considering $k\ge 2$, for these cases we can still use point $2.$ and $3.$ .
            \item[$\epsilon=1$ :] we start supposing that $r\ge 2$. With this hypothesis we can argue as the case $\epsilon=0$ dividing the value that $k$ can take:
            \begin{itemize}
                \item if $\frac{3}{4}r-\frac{1}{2}\le k-1 \le r-1$, we divide $M_k$ in five groups with $(r-k,k+1,r-k-1,k+1,r-k)$ rows each. Arguing as in $1.$ we obtain our claim;
                \item if $k-1=r$, then we argue as in $2.$;
                \item finally, if $k-1\ge r+1$, then we argue as in $3.$.
            \end{itemize}
            The final case $r=1$ correspond to $d=4$. We have to treat it separately because the inequality $\frac{3}{4}r-\frac{1}{2}\le r-1$ is not true. But again, since we are considering $k\ge 2$, we have only to consider $k=2$. Looking at the general case $k-1 = r$ we can conclude in the same way.
            \item[$\epsilon=2$ :] in this case we do not need to take restrictions to the value of $r$. Again we argue similarly to the case $\epsilon=0$ as follows:
            \begin{itemize}
                \item if $\frac{3}{4}r-\frac{1}{4}\le k-1\le r$, we divide $M_k$ in five groups with $(r-k,k+1,r-k,k+1,r-k)$ rows each. Arguing as in $1.$ we obtain our claim;
                \item if $k-1=r$ the statement is proved as in $2.$;
                \item finally, if $k-1\ge r+1$ then as in the case $3.$, we can conclude.
            \end{itemize}
        \end{itemize}
        
        Summarising, we have just proved that the matrices $N_k'$ and $M_k$ have full rank for $k\ge 2$. This implies that the Artinian Gorenstein algebra $S/\Ann _S(f)$ associated to the Perazzo 3-fold $f=x_0p_0+x_1p_1+x_2p_2$ has $h$-vector as described in Theorem \ref{upper}.
    \end{example}
    
    \begin{proposition}\label{lower} Let $d\geq 4$. Let $f$ be a form that defines a Perazzo $3$-fold of degree $d$. Let $R=K[x_0,x_1,x_2,u,v]$ and let $S=K[y_0,y_1,y_2,U,V]$ be the ring of differential operators on $R$. The minimum possible $h$-vector of  the Artinian Gorenstein algebra $A=S/Ann_S(f)$ is:
    $$ (1, \ 5, \ 6, \ 6, \ \cdots \ 6, \ 6, \ 5, \ 1).
    $$
    \end{proposition}
    \begin{proof} By assumption $f=p_0x_0+p_1x_1+p_2x_2+g$, where $p_0,p_1,p_2$ are linearly independent, but algebraically dependent. Let  $h_A=(h_0, \ h_1, \ h_2, \   \cdots ,\ h_d)$ be the $h$-vector of $A$. By Proposition \ref{hilbert function} we have that $h_0=h_d=1$ and $h_1=h_{d-1}=5$. In a similar way as argued in Proposition \ref{upper}, we are interested to compute $h_k$ for every $k\leq\frac{d}{2}$. Furthermore, its value is minimal if and only if the rank of $M_k$ and $N'_k$ are minimal. In Theorem \ref{lowcharacterisation} we have proved that both these matrices have rank greater or equal to $3$. Then by Proposition \ref{hilbert function} we obtain
    $$h_k=4k+1-m_k-n_k=\rank M_k+\rank N'_k\ge 6.$$
    
    To finish the proof it suffices to give a series of  examples, one for each degree, proving that the lower bound is achieved. This is trivial because every form of Theorem \ref{lowcharacterisation} has the desired $h$-vector.
    \end{proof}
    
    As we have just seen in Propositions \ref{upper} and \ref{lower}, we can expect that a general choice of the $p_i$ and $g$, with the suitable assumptions, will give a Hilbert vector which is maximal. On the other hand, we expect that the choices, which give a minimal Hilbert vector, correspond to very special $p_0,p_1,p_2$. As we have seen in Theorem \ref{lowcharacterisation}, we can classify all possible Perazzo $3$-folds, with minimal $h$-vector, in every degree. In fact, the classes are the cases when the rank of $M_k$ and $N'_k$ is minimal. So, we can rewrite Theorem \ref{lowcharacterisation} as follows.
    
    \begin{theorem}\label{lowcharacterisationreprise}
        Let $X=V(f)\subseteq\PP^4$ be a Perazzo $3$-fold of degree $d\ge 4$. Then the Hilbert vector of $A=S/Ann_S(f)$ is minimal as in Proposition \ref{lower}, if and only if, up to a change of variables, one of the following cases occur:
        \begin{itemize}
            \item[(i)]   $f(x_0,x_1,x_2,u,v)=u^{d-1}x_0+u^{d-2}vx_1+u^{d-3}v^2x_2+ au^d+bu^{d-1}v+cu^{d-2}v^2$ with $a,b,c\in K$, or
            \item[(ii)]  $f(x_0,x_1,x_2,u,v)=u^{d-1}x_0+u^{d-2}vx_1+v^{d-1}x_2+au^d+bu^{d-1}v+cv^{d}$ with $a,b,c\in K$, or
            \item[(iii)]  $f(x_0,x_1,x_2,u,v)=u^{d-1}x_0+(\lambda u+\mu v)^{d-1}x_1+v^{d-1}x_2+au^d+b(\lambda u+ \mu v)^d+cv^d$ with $\lambda, \mu\in K^*$ and $a,b,c\in K$.
        \end{itemize}
    \end{theorem}
    
    \begin{remark}
        As we have noticed in Lemma \ref{generators}, the Hilbert function of the algebra $S/\Ann_S(f)$ depends only on the plane $\pi=\langle p_0, p_1, p_2\rangle\subset\PP^{d-1}$ and not on the choice of a set of generators of $\pi$. In the proof of Theorem \ref{lowcharacterisation} we have proved that the $h$-vector is minimal if and only if the plane $\pi$ is in one of the following special positions:
        \begin{itemize}
            \item $\pi$ intersects $C_{d-1}$ at three distinct points which are necessarily linearly independent, that corresponds to the case (iii);
            \item $\pi$ contains the tangent line of some point at $C_{d-1}$ and also meets $C_{d-1}$ at a second different point, that corresponds to the case (ii);
            \item $\pi$ is an osculating plane of the rational normal curve $C_{d-1}$, that corresponds to the case (i).
        \end{itemize}
    \end{remark}
    
    As a final note Propositions \ref{upper} and \ref{lower} imply that the only possible Hilbert vector for a form defining a Perazzo $3$-fold of degree $4$ is $h=(1,5,6,5,1)$. Being the lower possible Hilbert vector, Theorem \ref{lowcharacterisationreprise} gives us a complete description of them.

    \section{Perazzo 3-folds and the WLP}\label{wlp}
    
    The goal of this section is to analyse whether the Artinian  Gorenstein algebra $A$ associated to a Perazzo 3-fold $X\subset \PP^4$ has the WLP, and the SLP. As we have already seen, by Corollary \ref{watanabecorollario}, since the Perazzo $3$-folds have vanishing first Hessian, it follows that the associated algebras fail the Strong Lefschetz Property. In particular, if $f$ is a Perazzo $3$-fold and $A=S/\Ann_S(f)$, the map
    \begin{equation}\label{formulagnagna}
        \times L^{d-2}: A_{1}  \longrightarrow  A_{d-1}
    \end{equation}
    
    is not an isomorphism for every possible choice of $L\in A_1$.
    
    \pause
    If $d=3$, clearly $A$ fails also WLP. In fact, the Hilbert vector of $A$ is $(1,5,5,1)$; in particular, the map (\ref{formulagnagna}) is just the multiplication $\times L: A_1\to A_2$. So, the algebra $A$ fails to have the Weak Lefschetz property.
    
    \pause
    The case $d=4$, has been studied in more generality by Gondim in \cite{gondim}. He proved that any Gorenstein algebra with codimension $5$ and socle degree $4$ always has the WLP. In particular, Gorenstein algebras associated to Perazzo $3$-folds satisfy this property, and so the WLP holds.
    
    \pause
    In the general case $d\ge 5$, we prove in Theorem \ref{main1} that if the Hilbert vector is maximum, then the WLP fails. On the contrary, if the Hilbert vector is minimum, then the WLP holds true (Theorem \ref{main2}). The middle cases seem to be too general, and nothing can be said. We give in Example \ref{discrepanza} two different Perazzo $3$-folds with an intermediary $h$-vector, but with different behaviour with respect to the WLP. Therefore, the fact that the WLP holds or not for a certain $K$-algebra $A$, could not depend on the Hilbert vector of $A$ only.

    \begin{theorem} \label{main1} Let $X\subset \PP^4$ be a Perazzo 3-fold of degree $d\ge 5$ defined by the form $$f=x_0p_0(u,v)+x_1p_1(u,v)+x_2p_2(u,v)+g(u,v).$$
    If $A=S/\Ann_S(f)$ has maximum h-vector, then $A$ fails WLP.
    \end{theorem}
    
    \begin{proof}
    The key argument is that, with those assumptions, the $\left[\frac{d}{2}\right]$-th Hessian is always zero. When $d$ is odd, then we conclude as in the case $d=3$. Instead, if $d$ is even, then we use that the Hilbert vector is flat. We proceed with the proof distinguishing such two cases.
    
    \noindent {\bf Case 1:} $d$ is odd. Write $d=2r+1$.To show that $A$ fails WLP, we will  prove that for any $L\in [A]_1$,  the multiplication map
    \[
    \times L: A_{r}  \longrightarrow  A_{r+1}
    \]
    is not bijective. By Corollary \ref{watanabecorollario}, it is enough to see the vanishing of the $r$-th Hessian $\hess _f^r$ of $f=x_0p_0(u,v)+x_1p_1(u,v)+x_2p_2(u,v)+g(u,v)$ with respect to a suitable basis $\mathcal{B}$ of $A_r$ (and, hence, all bases).
    First we can notice that a basis $\mathcal{B}$ made of classes with a monomial representative always exists. So, $\Hess _f^r$ is just a submatrix of dimension $h_r\times h_r$ of the following matrix:
    $$\left( \frac{\partial^{2r} f}{\partial u^\alpha\partial v^\beta\partial x_0^\gamma\partial x_1^\delta\partial x_2^\eta}\right)_{\alpha+\beta+\gamma+\delta+\eta=2r}$$
    where monomials are lexicographically ordered (for simplicity). Knowing that $f$ is linear in the variables $x_0,x_1,x_2$, the above matrix can be partially computed as:

    \begin{equation*}
        \hspace*{-1.5cm}\left(
        \begin{array}{cccc|cccc|ccc}
            \frac{\partial^{2r} f}{\partial u^{2r}} & \frac{\partial^{2r} f}{\partial u^{2r-1}\partial v} & \cdots  & \frac{\partial^{2r} f}{\partial u^{r}\partial v^{r}} & \frac{\partial^{2r-1} p_0}{\partial u^{2r-1}} &\frac{\partial^{2r-1} p_0}{\partial u^{2r-2}\partial v} & \cdots & \frac{\partial^{2r-1} p_2}{\partial u^{r}\partial v^{r-1}} & 0 & \cdots & 0 \Bstrut{}\\
            \frac{\partial^{2r} f}{\partial u^{2r-1}\partial v} & \frac{\partial^{2r} f}{\partial u^{2r-2}\partial v^2} & \cdots  & \frac{\partial^{2r} f}{\partial u^{r-1}\partial v^{r+1}} &\frac{\partial^{2r-1} p_0}{\partial u^{2r-2}\partial v} &\frac{\partial^{2r-1} p_0}{\partial u^{2r-3}\partial v^2} & \dots & \frac{\partial^{2r-1} p_2}{\partial u^{r-1}\partial v^{r}} & 0 & \cdots & 0 \Bstrut{}\\
            \vdots & \vdots & \ddots  & \vdots & \vdots & \vdots & \ddots  & \vdots & \vdots & \ddots & \vdots\Bstrut{}\\
            \frac{\partial^{2r} f}{\partial u^{r}\partial v^r} & \frac{\partial^{2r} f}{\partial u^{r-1}\partial v^{r+1}} & \cdots  & \frac{\partial^{2r} f}{\partial v^{2r}} & \frac{\partial^{2r-1} p_0}{\partial u^{r-1}\partial v^r} &\frac{\partial^{2r-1} p_0}{\partial u^{r-2}\partial v^{r+1}} & \dots & \frac{\partial^{2r-1} p_2}{\partial u^{2r-1}} & 0 & \cdots & 0\Bstrut{}\\
            \hline
            \frac{\partial^{2r-1} p_0}{\partial u^{2r-1}} & \frac{\partial^{2r-1} p_0}{\partial u^{2r-2}\partial v} & \cdots & \frac{\partial^{2r-1} p_0}{\partial u^{r-1}\partial v^r}  & 0 & \cdots & \cdots & 0 & 0 & \cdots & 0 \Tstrut{}\Bstrut{}\\
            \frac{\partial^{2r-1} p_0}{\partial u^{2r-2}\partial v} & \frac{\partial^{2r-1} p_0}{\partial u^{2r-3}\partial v^2} & \cdots  & \frac{\partial^{2r-1} p_0}{\partial u^{r-2}\partial v^{r+1}} & 0 & \dots & \dots & 0 & 0 & \cdots & 0 \Bstrut{}\\
            \vdots & \vdots & \ddots  & \vdots & \vdots & \vdots & \ddots  & \vdots & \vdots & \ddots & \vdots \Bstrut{}\\
            \frac{\partial^{2r-1} p_2}{\partial u^r\partial v^{r-1}} & \frac{\partial^{2r-1} p_2}{\partial u^{r-1}\partial v^{r}} & \cdots & \frac{\partial^{2r-1} p_2}{\partial u^{2r-1}}  & 0 & \cdots & \cdots & 0 & 0 & \cdots & 0 \Bstrut{}\\
            \hline
            0 & \cdots & \cdots  & 0 & 0 & \cdots & \cdots  & 0 & 0 & \cdots & 0\Tstrut{}\\
            \vdots & \ddots & \ddots  & \vdots & \vdots & \ddots & \ddots  & \vdots & \vdots & \ddots & \vdots\\
            0 & \cdots & \cdots  & 0 & 0 & \cdots & \cdots  & 0 & 0 & \cdots & 0\\
        \end{array}
        \right).
    \end{equation*}
    
    \vskip 2mm
    The three vertical (respectively, horizontal) blocks  are composed respectively by $r+1,\, 3r, \, \binom{r+4}{4}-(4r+1)$ columns (respectively, rows). Thus every possible choice of a $h_r\times h_r$ submatrix turns out to have at least an all zero sub-submatrix of size $(h_r-(r+1))\times (h_r-(r+1))$. We now use the hypothesis of $A$ to have maximum $h$-vector and Proposition \ref{upper} to obtain that $h_r=2r+3$. We have just proved that $\Hess _f^r$, matrix of dimension $(2r+3)\times (2r+3)$, has an all zero submatrix of dimension $(r+2)\times(r+2)$: this implies that $\hess _f^r$ identically vanishes.
    
    \pause
    \noindent {\bf Case 2:} $d$ is even. Write $d=2r+2$. Note that, since the $h$-vector is maximum, then $h_r=h_{r+1}=h_{r+2}=2r+4$. Using again the hessian criterion of Watanabe (Corollary \ref{watanabecorollario}), we will check that for any $L\in A_1$,  the multiplication map
    \[
    \times L^2: A_{r}  \longrightarrow  A_{r+2}
    \]
    is not bijective. This implies that for any $L\in [A]_1$,  the multiplication maps
    \[
    \times L: A_{r}  \longrightarrow  A_{r+1}\text{, and }\times L: A_{r+1}  \longrightarrow  A_{r+2}
    \]
    are not bijective and, hence, $A$ fails the WLP. 
    
    Same adapted argument of the previous case can be used also here. In fact, the matrix to be considered is $\Hess _f^r$ which is now of size $(2r+4)\times(2r+4)$ which is even bigger of the previous case. Thus, as discussed above, its determinant is always zero.
    \end{proof}
    
    In contrast with the last result we have  that if an Artinian Gorenstein algebra $A$ associated to a Perazzo 3-fold has  minimum $h$-vector, then $A$ has the WLP.
    
    \begin{theorem} \label{main2} Let $X\subset \PP^4$ be a Perazzo 3-fold of degree  $d\ge 5$ and equation $$f=x_0p_0(u,v)+x_1p_1(u,v)+x_2p_2(u,v)+g(u,v)\in R=K[x_0,x_1,x_2,u,v]_d.$$ Let $S=K[y_0,y_1,y_2,U,V]$ be the ring of differential operators on $R$. If $A=S/\Ann_S(f)$ has minimum h-vector, then $A$ has WLP.
    \end{theorem}
    \begin{proof}
        In the following proof we will use the characterisation given in Theorem \ref{lowcharacterisationreprise} to study case by case the WLP. In \cite{fiorindomezzettimiroroig} the same statement is proved using M. Green's theorem \cite[Theorem 1]{green}.
        
        \pause
        Let $S=K[y_0,y_1,y_2,U,V]$ be the ring of differential operators. By Theorem \ref{main2}, we have $h_2=h_3=\cdots=h_{d-2}=6$. By Proposition \ref{checkweak}, $A$ has the WLP if and only if the multiplication map
        $$\times L:A_2\longrightarrow A_3$$
        is an isomorphism for some $L\in A_1$. To prove it, we instead check that the multiplication map
        $$\times L^{d-4}: A_2\longrightarrow A_{d-2}$$
        is an isomorphism for some $L\in A_1$, i.e. the $2$-th Hessian of $f$ does not vanish. Since we have already classified $f$ in Theorem \ref{lowcharacterisationreprise}, we use that classification to explicit compute its second Hessian.
        \begin{itemize}
            \item[(i)]   $f(x_0,x_1,x_2,u,v)=u^{d-1}x_0+u^{d-2}vx_1+u^{d-3}v^2x_2+ au^d+bu^{d-1}v+cu^{d-2}v^2$ with $a,b,c\in K$. A basis of $A_2$ is $\mathcal{B}=\{U^2,UV,V^2,y_0U,y_1U,y_2U\}$. So, the second Hessian is
            $$\hess_f^2=\det\left(
            \begin{array}{cc}
                C & B \\
                B^T & 0
            \end{array}
            \right)=-(\det B)^2,$$
            where $C,B$ are $3\times 3$ matrices. In particular the matrix $B$ has the form
            $$B=(d-3)\left(\begin{array}{ccc}
                (d-1)(d-2)u^{d-4} & \ast &\ast \\
                0 & (d-2)u^{d-4} &\ast \\
                0 & 0 & u^{d-4}
            \end{array}\right).$$
            Since $B$ has non zero determinant, $f$ does not have $2$-th vanishing Hessian. Moreover, the non Lefschetz elements have the form $L=k_0y_0+k_1y_1+k_2y_2+k_3V$ with $k_i\in K$.
            \item[(ii)]  $f(x_0,x_1,x_2,u,v)=u^{d-1}x_0+u^{d-2}vx_1+v^{d-1}x_2+au^d+bu^{d-1}v+cv^{d}$ with $a,b,c\in K$. A basis of $A_2$ is $\mathcal{B}=\{U^2,UV,V^2,y_0U,y_1U,y_2V\}$. So, the second Hessian is
            $$\hess_f^2=\det\left(
            \begin{array}{cc}
                C & B \\
                B^T & 0
            \end{array}
            \right)=-(\det B)^2,$$
            where $C,B$ are $3\times 3$ matrices. In particular the matrix $B$ has the form
            $$B=(d-2)(d-3)\left(\begin{array}{ccc}
                (d-1)u^{d-4} & \ast &\ast \\
                0 & u^{d-4} &\ast \\
                0 & 0 & (d-1)v^{d-4}
            \end{array}\right).$$
            Since $B$ has non zero determinant, $f$ does not have $2$-th vanishing Hessian. Moreover, the non Lefschetz elements have the form $L=k_0y_0+k_1y_1+k_2y_2+k_3U$, or $L=k_0y_0+k_1y_1+k_2y_2+k_4V$ with $k_i\in K$.
            \item[(iii)]  $f(x_0,x_1,x_2,u,v)=u^{d-1}x_0+(\lambda u+\mu v)^{d-1}x_1+v^{d-1}x_2+au^d+b(\lambda u+ \mu v)^d+cv^d$ with $\lambda, \mu\in K^*$ and $a,b,c\in K$. A basis of $A_2$ is $\mathcal{B}=\{U^2,UV,V^2,y_0U,y_1U,y_2V\}$. So, the second Hessian is
            $$\hess_f^2=\det\left(
            \begin{array}{cc}
                C & B \\
                B^T & 0
            \end{array}
            \right)=-(\det B)^2,$$
            where $C,B$ are $3\times 3$ matrices. In particular the matrix $B$ has the form
            $$B=(d-1)(d-2)(d-3)\left(\begin{array}{ccc}
                u^{d-4} & \ast & 0 \\
                0 & \lambda^2\mu(\lambda u+ \mu v)^{d-4} & 0 \\
                0 & \ast & v^{d-4}
            \end{array}\right).$$
            Since $B$ has non zero determinant, $f$ does not have $2$-th vanishing Hessian. Moreover, the non Lefschetz elements have the form $L=k_0y_0+k_1y_1+k_2y_2+k_3U$, or  $L=k_0y_0+k_1y_1+k_2y_2+k_4V$, or $L=k_0y_0+k_1y_1+k_2y_2+k_3(-\mu U+\lambda V)$ with $k_i\in K$. \qedhere
        \end{itemize}
    \end{proof}

    \begin{remark}
       As a consequence of Theorem \ref{main2}, all forms of degree $d$ which define a Perazzo $3$-fold with minimum $h$-vector are examples of forms with zero first order hessian, and all hessians of order $k$ different from zero, for  $2\le k\le\lfloor \frac{d}{2}\rfloor$.
    \end{remark}

    For Gorenstein Artinian algebras associated to Perazzo 3-folds $X$ in $\PP^4$ and with intermediate $h$-vector both possibilities occur.
    
    \begin{example}\label{discrepanza}
    1.-  Let $X\subset \PP^4$ be the Perazzo 3-fold of equation
    $$ f_1(x_0,x_1,x_2,u,v)=(u^5+u^4v)x_0+u^3v^2x_1+v^5x_2\in K[x_0,x_1,x_2,u,v]_6.
    $$
    Let $S=K[y_0,y_1,y_2,U,V]$ be the ring of differential operators on $R$. We can now compute the matrices $M_k$ and $N_k'$ as follows:
    \begin{equation*}
        \hspace{-0.1cm}M_2=\left(\begin{array}{cc|cc|cc}
            1 & \frac{1}{5} & 0 & 0 & 0 & 0\\
            \frac{1}{5} & 0 & 0 & \frac{1}{10} & 0 & 0\\
            0 & 0 & \frac{1}{10} & 0 & 0 & 0\\
            0 & 0 & 0 & 0 & 0 & 0\\
            0 & 0 & 0 & 0 & 0 & 1
        \end{array}\right),\,
        M_3=\left(\begin{array}{ccc|ccc|ccc}
            1 & \frac{1}{5} & 0 & 0 & 0 & \frac{1}{10} & 0 & 0 & 0\\
            \frac{1}{5} & 0 & 0 & 0 & \frac{1}{10} & 0 & 0 & 0 & 0\\
            0 & 0 & 0 & \frac{1}{10} & 0 & 0 & 0 & 0 & 0\\
            0 & 0 & 0 & 0 & 0 & 0 & 0 & 0 & 1
        \end{array}\right)
    \end{equation*}
    \begin{equation*}
        N'_2=\begin{pmatrix}
            1 & \frac{1}{5} & 0\\
            \frac{1}{5} & 0 & 0\\
            0 & 0 & 0\\
            0 & 0 & 0\\
            \hline
            0 & 0 & \frac{1}{10}\\
            0 & \frac{1}{10} & 0\\
            \frac{1}{10} & 0 & 0\\
            0 & 0 & 0\\
            \hline
            0 & 0 & 0\\
            0 & 0 & 0\\
            0 & 0 & 0\\
            0 & 0 & 1\\
            \hline
            0 & 0 & 0\\
            0 & 0 & 0\\
            0 & 0 & 0\\
            0 & 0 & 0\\
            0 & 0 & 0\\
        \end{pmatrix},\text{ and }
        N_3'=\left(\begin{array}{cccc}
            1 & \frac{1}{5} & 0 & 0\\
            \frac{1}{5} & 0 & 0 & 0\\
            0 & 0 & 0 & 0\\
            \hline
            0 & 0 & \frac{1}{10} & 0\\
            0 & \frac{1}{10} & 0 & 0\\
            \frac{1}{10} & 0 & 0 & 0\\
            \hline
            0 & 0 & 0 & 0\\
            0 & 0 & 0 & 0\\
            0 & 0 & 0 & 1\\
            \hline
            0 & 0 & 0 & 0\\
            0 & 0 & 0 & 0\\
            0 & 0 & 0 & 0\\
            0 & 0 & 0 & 0\\
        \end{array}\right)
    \end{equation*}
    
    By Proposition \ref{hilbert function}, he Artinian Gorenstein algebra  $A=S/Ann_S(f)$ has $h$-vector
    $$h=(1, \ 5, \ 7, \ 8, \ 7, \ 5, \ 1).$$
    Using Macaulay2 \cite{mac2} we can compute the generators of $\Ann _S(f_1)$:
    \begin{equation*}
     \begin{split}
        \Ann _S(f_1)=\langle y_2U,y_0U-y_0V-y_1V,y_0^2,y_1^2,y_2^2,y_0y_1,y_0y_2,y_1y_2,\\
        y_0V^2,UV^3,y_1U^3-y_2V^3,U^5-U^4V,U^6,V^6 \rangle;
     \end{split}
    \end{equation*}
    moreover, we also check that for a general linear form $L\in [A]_1$, the multiplication map
     \[
    \times L: A_{2}  \longrightarrow  A_{3}
    \]
    is injective. Hence, by Proposition \ref{checkweak}, $A$ satisfies the WLP.
    
    \pause
    2.- Let $X\subset \PP^4$ be the Perazzo 3-fold of equation
    $$ f_2(x_0,x_1,x_2,u,v)=u^6x_0+u^3v^3x_1+v^6x_2\in K[x_0,x_1,x_2,u,v]_7.
    $$
    Let $S=K[y_0,y_1,y_2,U,V]$ be the ring of differential operators on $R$. We can now compute the matrices $M_k$ and $N_k'$ as follows:
    \begin{equation*}
        \hspace{-0.1cm}M_2=\left(\begin{array}{cc|cc|cc}
            1 & 0 & 0 & 0 & 0 & 0\\
            0 & 0 & 0 & 0 & 0 & 0\\
            0 & 0 & 0 & \frac{1}{20} & 0 & 0\\
            0 & 0 & \frac{1}{20} & 0 & 0 & 0\\
            0 & 0 & 0 & 0 & 0 & 0\\
            0 & 0 & 0 & 0 & 0 & 1
        \end{array}\right),\,
        M_3=\left(\begin{array}{ccc|ccc|ccc}
            1 & 0 & 0 & 0 & 0 & 0 & 0 & 0 & 0\\
            0 & 0 & 0 & 0 & 0 & \frac{1}{20} & 0 & 0 & 0\\
            0 & 0 & 0 & 0 & \frac{1}{20} & 0 & 0 & 0 & 0\\
            0 & 0 & 0 & \frac{1}{20} & 0 & 0 & 0 & 0 & 0\\
            0 & 0 & 0 & 0 & 0 & 0 & 0 & 0 & 1
        \end{array}\right)
    \end{equation*}
    \begin{equation*}
        N'_2=\begin{pmatrix}
            1 & 0 & 0\\
            0 & 0 & 0\\
            0 & 0 & 0\\
            0 & 0 & 0\\
            0 & 0 & 0\\
            \hline
            0 & 0 & \frac{1}{20}\\
            0 & \frac{1}{20} & 0\\
            \frac{1}{20} & 0 & 0\\
            0 & 0 & 0\\
            0 & 0 & 0\\
            \hline
            0 & 0 & 0\\
            0 & 0 & 0\\
            0 & 0 & 0\\
            0 & 0 & 0\\
            0 & 0 & 1\\
            \hline
            0 & 0 & 0\\
            0 & 0 & 0\\
            0 & 0 & 0\\
            0 & 0 & 0\\
            0 & 0 & 0\\
            0 & 0 & 0\\
        \end{pmatrix},\text{ and }
        N_3'=\left(\begin{array}{cccc}
            1 & 0 & 0 & 0\\
            0 & 0 & 0 & 0\\
            0 & 0 & 0 & 0\\
            0 & 0 & 0 & 0\\
            \hline
            0 & 0 & 0 & \frac{1}{20}\\
            0 & 0 & \frac{1}{20} & 0\\
            0 & \frac{1}{20} & 0 & 0\\
            \frac{1}{20} & 0 & 0 & 0\\
            \hline
            0 & 0 & 0 & 0\\
            0 & 0 & 0 & 0\\
            0 & 0 & 0 & 0\\
            0 & 0 & 0 & 1\\
            \hline
            0 & 0 & 0 & 0\\
            0 & 0 & 0 & 0\\
            0 & 0 & 0 & 0\\
            0 & 0 & 0 & 0\\
            0 & 0 & 0 & 0\\
        \end{array}\right)
    \end{equation*}
    
    By Proposition \ref{hilbert function}, the Artinian Gorenstein algebra  $A=S/Ann_S(f)$ has $h$-vector
    $$h=(1, \ 5, \ 7, \ 9, \ 9, \ 7, \ 5, \ 1).$$
    Using Macaulay2 \cite{mac2} we can compute the generators of $\Ann _S(f_1)$:
    \begin{equation*}
        \begin{split}
            \Ann _S(f_2)=\langle y_0^2,y_1^2,y_2^2,y_0y_1,y_0y_2,y_1y_2, y_0v, y_2u,  y_1U^3-y_2V^3,\\
            y_0U^3-y_1V^3, UV^4, U^4V, V^7, U^7 \rangle.
        \end{split}
    \end{equation*}
    
    Setting $\mathcal{B}=\{U^3,U^2V,UV^2,V^3,y_0U^2,y_1U^2,y_1UV,y_1V^2,y_2V^2\}$ as a base of $A_3$, we get that the $3$-rd hessian of $f$ is identically zero. In fact, the matrix
    $$\Hess_3(f)=\left(\begin{array}{cccc|ccccc}
        \ast & \ast & \ast & \ast & \ast & \ast & \ast & \ast & \ast \\
        \ast & \ast & \ast & \ast & \ast & \ast & \ast & \ast & \ast \\
        \ast & \ast & \ast & \ast & \ast & \ast & \ast & \ast & \ast \\
        \ast & \ast & \ast & \ast & \ast & \ast & \ast & \ast & \ast \\
        \hline
        \ast & \ast & \ast & \ast & 0 & 0 & 0 & 0 & 0 \\
        \ast & \ast & \ast & \ast & 0 & 0 & 0 & 0 & 0 \\
        \ast & \ast & \ast & \ast & 0 & 0 & 0 & 0 & 0 \\
        \ast & \ast & \ast & \ast & 0 & 0 & 0 & 0 & 0 \\
        \ast & \ast & \ast & \ast & 0 & 0 & 0 & 0 & 0 
    \end{array}\right)$$
    has vanishing determinant. Thus, by Theorem \ref{watanabe}, for any linear form $L\in A_1$, the multiplication map
     \[
    \times L: A_{3}  \longrightarrow  A_{4}
    \]
    is not bijective and, hence, $A$ fails the WLP.
    \end{example}
    
    \section{Extension to the general case}
    In this final section, we generalise some of the properties, we have seen in the previous two sections, to a general form with or without vanishing Hessian in $5$ variables.
    
    \pause
    First we notice the following: if $\Delta$ is a Perazzo $3$-fold, $g$ is a homogeneous polynomial $g$ in the variables $u$ and $v$, then $g\Delta$ is a Perazzo $3$-fold with $\Tilde{p_i}=gp_i$. Since the degree of a Perazzo $3$-fold has to be at least three, we have that all possible forms $f$ with vanishing Hessian and degree at most $5$ are Perazzo $3$-folds, up to a change of variables. Therefore we obtain an alternative proof of a result of Gondim (\cite[Theorem 3.5]{gondim}) which gives a more strong result with respect to Theorem \ref{main2} in the case $d=4$.
    
    \begin{corollary}
        Let $f\in K[x_0,x_1,x_2,u,v]_4$ be any form of degree $4$ involving exactly $5$ variables, i.e. $V(f)$ is not a cone. Then $S/ \Ann _S(f)$ always has the WLP.
    \end{corollary}
    \begin{proof}
        We distinguish two cases: $f$ has, or not, vanishing Hessian. If $f$ has not vanishing Hessian, then by Corollary \ref{watanabecorollario}, $A$ has the SLP. On the other hand, if $f$ has vanishing Hessian, then $f$ is a Perazzo $3$-fold of degree $4$ with Hilbert vector $(1,5,6,5,1)$. So, we have already classified it as in Theorem \ref{lowcharacterisationreprise} with $d=4$. A straightforward computation proves that for a general element $L\in A_1$, the map
        $$\times L:A_1\longrightarrow A_2$$
        has full rank. Proposition \ref{checkweak} completes the proof.
    \end{proof}
    
    We consider the case of a form $f$ of degree $d=5$ with vanishing Hessian: we have that $f$, up to a change of variables, is a Perazzo $3$-fold. The Hilbert vector $h$ of $A=S/\Ann_S(f)$ satisfy only one of the two following possibilities:
    $$h=(1,5,6,6,5,1)\text{, or } h=(1,5,7,7,5,1).$$
    As a consequence of Theorems \ref{main1}, and \ref{main2}, in the first case $A$ satisfies the Weak Lefschetz property, while in the second case it does not. This facts give us the following result.
    
    \begin{corollary}\label{quintic}
    Let $f$ be a form of degree $d=5$ involving exactly $5$ variables. If $f$ has vanishing Hessian, then the Artinian Gorenstein algebra $S/ \Ann _S(f)$, associated to $f$, has the WLP if and only if, after a possible change of coordinates, one of the following cases holds:
    \begin{itemize}
        \item[(i)]   $f(x_0,x_1,x_2,u,v)=u^4x_0+u^3vx_1+u^2v^2x_2+au^5+bu^4v+cu^3v^2$ with $a,b,c\in K$, or
        \item[(ii)]  $f(x_0,x_1,x_2,u,v)=u^4x_0+u^3vx_1+v^4x_2+au^5+bu^4v+cv^5$ with $a,b,c\in K$, or
        \item[(iii)]  $f(x_0,x_1,x_2,u,v)= u^4x_0+(\lambda u+ \mu v)^4x_1+v^4x_2+au^5+b(\lambda u+ \mu v)^5+cv^5$ with $\lambda, \mu\in K^*$ and $a,b,c\in K$.
    \end{itemize}
    Furthermore, if $f$ does not have vanishing Hessian, then $A$ has both, or none the Lefschetz properties. In particular, $A$ has both the Lefschetz properties if and only if $\hess_2(f)\neq 0$.
    \end{corollary}
    \begin{proof}
    It follows from Corollary \ref{watanabecorollario}, Theorems \ref{lowcharacterisationreprise}, \ref{main1}, and \ref{main2}.
    \end{proof}
    
    If the degree is greater that $5$, the possibilities of the Hilbert vector of $A$ become very numerous. Even in the first case $d=6$, the treatment is complicated and many variables have to be considered. If we consider a form $f\in K[x_0,x_1,x_2,u,v]_6$ with vanishing Hessian, not a cone, then, as we have stated many times, $f\in K[u,v][\Delta]$ with $\Delta$ is a Perazzo $3$-fold with $g=0$. Since $f$ has degree $6$, terms like $\Delta^2$ could appear in $f$. Now, it is not so clear if $f$ is, up to a linear change of variables, a Perazzo $3$-fold. Indeed, we will see in Remark \ref{perazzopersempre} that this statement is false. Moreover, it seems that no relations between the annihilator of $f$ and the annihilator of $\Delta$ can be found.
    
    \pause
    Even if we restrict our attention to Perazzo $3$-folds, we have no clue about the study of the Hilber vector and the WLP. Let now $f=p_0x_0+p_1x_1+p_2x_2+g$ be a Perazzo $3$-fold. For what we have done in Section \ref{wlp}, we know that the WLP holds when the Hilbert vector of $A=S/\Ann_S(f)$ is minimal, i.e. $h=(1,5,6,6,6,5,1)$. On the contrary, $A$ fails the WLP when its $h$-vector is maximal, i.e. $h=(1,5,8,8,8,5,1)$. In the middle cases, a priori, the $h$-vector of $A$ can be one of the following possibilities:
    \begin{equation*}
        \begin{array}{c@{\hspace{1cm}}c}
            \text{unimodal cases} & \text{non-unimodal cases}\\
            h=(1,5,6,7,6,5,1) & h=(1,5,7,6,7,5,1)\\
            h=(1,5,6,8,6,5,1) & h=(1,5,8,6,8,5,1)\\
            h=(1,5,7,7,7,5,1) & h=(1,5,8,7,8,5,1).\\
            h=(1,5,7,8,7,5,1)\\
        \end{array}
    \end{equation*}
    
    We notice that, a priori, we cannot exclude the cases when the Hilbert vector is non-unimodal. The $h$-vector of $A$ depends on the plane $\pi\subseteq\PP^5$ generated by the points $[p_0],[p_1],[p_2]$, and on the form $g$. By the definition of the matrices $N_k'$ and $N_k$ as in (\ref{mnnk}), if we consider the Hilbert vector $h$ of $f$ and the Hilbert vector $h'$ of $p_0x_0+p_1x_1+p_2x_2$ then $h'\le h$. In particular, the difference between $h'_k$ and $h_k$ depends only on $g$ and its catalecticant matrices.
    
    \pause
    As a consequence of Theorem \ref{osequenza} and Theorem \ref{sisequenza} we can deduce the following: 
    \begin{itemize}
        \item for the non-unimodal cases, the WLP is not satisfied by the Artinian Gorenstein algebra by Proposition \ref{unimodalsymmetric};
        \item since $6^{\langle 3\rangle}=7$, the case $h=(1,5,8,6,8,5,1)$ cannot occur by Theorem \ref{osequenza};
        \item if we consider $h=(1,5,6,8,6,5,1)$, we obtain $\Delta h=(1,4,1,2)$. Using Theorem \ref{sisequenza} and the fact $1^{\langle 2\rangle}=1$, we have that the $h$-vector $h=(1,5,6,8,6,5,1)$ cannot be the Hilbert vector of any Gorenstein Artinian algebra having the WLP;
        \item finally, for the remaining unimodal Hilbert vectors, Theorem \ref{sisequenza} does not imply if the $k$-algebra $A$ has or not the WLP.
    \end{itemize}
    
    As a final note, one possible tool to check the WLP is the $2$-nd Hessian of $f$, but in this case it is not fully reliable. Indeed, if $\hess_2(f)\neq 0$ then, for the general element $L\in A_1$, the map
    $$\times L^2:A_2\longrightarrow A_4$$
    is an isomorphism, but this map factorises through the maps
    $$\times L:A_2\longrightarrow A_3 \text{ and }\times L:A_3\longrightarrow A_4$$
    which turns out to be both full rank. Thus, by Proposition \ref{wlp}, $A$ satisfies the WLP. On the contrary, if $\hess_2(f)=0$ nothing can be said about the WLP. In fact, it could happen that $A$ has the WLP, so that the maps 
    $$\times L:A_2\longrightarrow A_3 \text{ and }\times L:A_3\longrightarrow A_4$$
    are both full rank for the general element $L\in A_1$; but, the map
    $$\times L^2:A_2\longrightarrow A_4$$
    could not  be an isomorphism for every possible $L\in A_1$.

    \chapter{The geometry of forms with vanishing Hessian}\label{chapter4}
    In this last chapter, we shortly introduce some notions and results to describe hypersurfaces with vanishing Hessian. In particular, we give a  geometrical description of particular Perazzo $3$-folds with minimal Hilbert vector; namely, the hypersurfaces of Theorem \ref{lowcharacterisationreprise} when $a=b=c=0$. In Section (\ref{gaussmapp}) we have introduced some basic concepts of algebraic geometry like the polar map, and the Gauss map of a hypersurface. We now briefly recall such arguments.
    
    \pause
    Given a projective space $\PP^n$ over a field $K$, we can define the set $(\PP^n)^*$ of all the hyperplanes of $\PP^n$. Over this set we define homogeneous coordinates so that a hyperplane $H$ of equation $a_0x_0+\dots+a_nx_n=0$ corresponds to the point of coordinates $[a_0:\dots:a_n]$. In this way, $(\PP^n)^*$ has a structure of a $n$-dimensional projective space. Moreover, a linear projective subspace $V\subset\PP^n$ of dimension $k$ corresponds to a linear projective subspace $V^*\subset(\PP^n)^*$ of dimension $n-1-k$, and viceversa.
    
    \pause
    Let $X=V(f)\subseteq\PP^n$ be a hypersurface, we denote by 
    $$\nabla _{f}:\PP^n \dashrightarrow (\PP^n)^* $$
    its {\em polar map} defined  by
    $$ \nabla _{f}([p])=\left[\frac{\partial f}{\partial x_0}(p):\frac{\partial f}{\partial x_1}(p): \dots: \frac{\partial f}{\partial x_n}(p)\right],$$
    and by 
    $$\G_X :X \dashrightarrow (\PP^n)^*$$
    the restriction of $\nabla _{f}$ to $X$, i.e. the {\em  Gauss map} of $X$, which associates to each smooth point of $X$ its embedded tangent space. The first map is not defined in $D=V(\frac{\partial f}{\partial x_0},\frac{\partial f}{\partial x_1}, \dots,\frac{\partial f}{\partial x_n})\subset\PP^n$, while the second map is not defined in $\Sing(X)$. Sometimes we will use the notation $X_{reg}$ to denote the open set $X\setminus\Sing(X)$. We define $X^*$, the dual variety of $X$, to be the closure of the image of $\G_X$, so that $X^*=\overline{\G_X(X_{reg})}$. Furthermore, we define $Z=\overline{\nabla _f(\PP^n\setminus D)}$ to be the closure of the image of the polar map. By their definition, we have that
    $$X^*\hookrightarrow Z.$$
    
    In general, it is possible to define the dual variety of an algebraic variety $X\subset\PP^n$ using the \emph{conormal variety} which lies in $X\times(\PP^n)^*$ and then projecting on the second factor (more details are in \cite[Section 1.5]{russo}). Moreover, since $X^*$ is a projective variety, we can repeat the process to obtain the bi-dual variety $X^{**}=(X^*)^*$. Since we can identify $(\PP^n)^{**}$ with $\PP^n$, we can think to $X^{**}$ as a projective variety in $\PP^n$. In general $X$ and $X^{**}$ could differ, but the following statement is well-known in fields of characteristic zero.
    
    \begin{theorem}\label{reflexive}
         Let $X\subset\PP^n$ be an irreducible projective variety over a field of characteristic zero. Then $X$ is reflexive, i.e. $X=X^{**}$.
    \end{theorem}
    \begin{proof}
        See \cite{wallace}.
    \end{proof}
    
    We remark that, from now on, we are working with an algebraically closed field $K$ with characteristic zero. Clearly, we have that the spaces $(\PP^n)^*$ and $\PP((K^{n+1}))$ are naturally isomorphic. So, we can view the map $\nabla_f$ as the quotient of the affine map
    \begin{equation*}
        \begin{split}
            \nabla_f^0:K^{n+1}&\longrightarrow K^{n+1}\\
            p&\longmapsto\left(\frac{\partial f}{\partial x_0}(p),\frac{\partial f}{\partial x_1}(p), \dots, \frac{\partial f}{\partial x_n}(p)\right)
        \end{split}
    \end{equation*}
    with respect to the natural $K^*$-action. This point of view allows to link between the Hessian matrix of $f$, and $Z$. In fact, we have the following formula
    $$\Hess(f)=\Jac(\nabla_f^0).$$
    Furthermore, we can compute the dimension of $Z$. Given a point $\nabla_f(P)$ in $Z$ we obtain
    $T_{\nabla_f(P)}(Z)=\PP(\Imm \Hess(f)(P));$
    thus, we obtain $$\dim Z=\rank\Hess(f)-1.$$
    These considerations have just proved the following statement.
    
    \begin{proposition}\label{znonètutto}
        Let $X=V(f)\subseteq\PP^n$ be a hypersurface. Then
        $$\hess(f)=0\quad \iff\quad Z\subsetneq(\PP^n)^*.$$
        In particular, $Z$ is a hypersurface if and only if $\rank\Hess(f)=n$.
    \end{proposition}
    
    Another possible way to see how the Hessian matrix affects the above construction is considering the Gauss map. In fact, the differential of $\G_X$ at some element $\mathbf{v}$ is related to $\mathbf{v}^T\cdot(\Hess(f)(\mathbf{v}))$. In particular, the following statement holds true.
    
    \begin{proposition}\label{segrerusso}
        Let $X=V(f)\subseteq\PP^n$ be an irreducible hypersurface. Then 
        $$\dim X^*=\rank(d\G_X)=\rank\Hess(f)-2;$$
        in particular, we have that $$\dim X^*=\dim Z-1.$$
    \end{proposition}
    \begin{proof}
        The proof of this statement is a consequence of \cite[Theorem 2]{segre}. A geometric interpretation is stated in \cite[Lemma 7.2.7]{russo}.
    \end{proof}
    \begin{corollary}\label{corollariotuttodiversonienteuguale}
        Let $X=V(f)\subseteq\PP^n$ be an irreducible hypersurface with vanishing Hessian, then
        $$X^*\subsetneq Z\subsetneq(\PP^n)^*.$$
    \end{corollary}
    \begin{proof}
        It is a direct consequence of Proposition \ref{znonètutto} and \ref{segrerusso}.
    \end{proof}
    
    The aim of this chapter is to give a geometric interpretation of hypersurfaces with vanishing Hessian. The first aspect, which we are going to state above, is true in any projective space $\PP^n$. The principal ingredient is the Gordan-Noether identity which has been stated in Theorem \ref{gnidentity}.
    
    \begin{proposition}\label{proposizionetecnica}
         Let $X=V(f)\subseteq\PP^n$ be a hypersurface with vanishing Hessian. Then
         \begin{enumerate}
             \item if $P\in\PP^n$ is a regular point of $X$ such that also $\nabla_f(P)$ is a regular point of $Z$, then
             $$\langle P, (T_{\nabla_f(P)}Z)^*\rangle\subseteq\nabla_f^{-1}(\nabla_f(P));$$
             \item for $P$ general, $\langle P, (T_{\nabla_f(P)}Z)^*\rangle$ is the irreducible component of $\overline{\nabla_f^{-1}(\nabla_f(P))}$ in which $P$ lies. In particular, $\overline{\nabla_f^{-1}(\nabla_f(P))}$ is the union of linear spaces of dimension $n+1-\rank H(f)=\codim Z$ passing through the linear space $(T_{\nabla_f(P)}Z)^*$;
             \item $Z^*\subseteq \Sing(X)$.
         \end{enumerate}
    \end{proposition}
    \begin{proof}
        See \cite[Corollary 7.3.14]{russo}.
    \end{proof}
    
     We remark the following: if $P\in X$ is a regular point of $X$, then $$\nabla_f^{-1}(\nabla_f(P))\cap X\equiv \G_X^{-1}(\G_X(P))$$ represents the set of regular points in $X$ whose tangent space coincides with $T_P X$. In particular, for the general $P\in X_{reg}$, such fibers are linear subspaces of dimension $n-1-\dim X^*$ that is always positive. This fact is just an interpretation of Theorem \ref{reflexive}.
    
    \pause
    We now restrict our attention to hypersurfaces with vanishing Hessian in $\PP^4$. The starting point of this work is Gordan-Noether's article \cite{gordannoether}. A. Franchetta analysed their classification to give a geometric interpretation of such hypersurfaces using the later called \emph{Franchetta hypersurfaces}; then C. Ciliberto, F. Russo and A. Simis gave also the geometry of their dual varieties.
    
    \begin{proposition}
         Let $X=V(f)\subseteq\PP^4$ be an irreducible hypersurface of degree $d\ge 3$ with vanishing Hessian, not a cone. Then $Z^*\subset\PP^4$ is an irreducible plane rational curve or, equivalently, $Z$ is a cone with vertex a line over an irreducible plane rational curve.
    \end{proposition}
    \begin{proof}
        See \cite[Lemma 7.4.13]{russo}.
    \end{proof}
    
    \begin{theorem}[Ciliberto, Franchetta, Russo, Simis]\label{cfrs}
        Let $X=V(f)\subseteq\PP^4$ be an irreducible hypersurface of degree $d\ge 3$, not a cone, with $f\in K[x_0,x_1,x_2,u,v]_d$. Then the following are equivalent:
        \begin{itemize}
            \item $X$ has vanishing Hessian.
            \item $X$ is a Franchetta hypersurface, i.e. $X$ is the union of a one-dimensional family $\Sigma$ of planes. Furthermore, there exists $C\subset\Pi=\PP^2$ a rational plane curve, of degree strictly greater that one and contained in $X$, satisfying the following conditions:
            \begin{enumerate}
                \item the planes of $\Sigma$ are all tangent to $C$;
                \item for the general hyperplane $H=\PP^3\subset\PP^4$, containing $\Pi$, the intersection $H\cap X$, off $\Pi$, is a finite union of planes of $\Sigma$ all tangent to $C$ at the same point.
            \end{enumerate}
            \item $X^*$ is a scroll surface of degree $d$ having a line $L$, of multiplicity greater than $\mu$, as directrix, and sitting in a $3$-dimensional rational cone $W_f$ which has $L$ as vertex. Furthermore, the general plane in $W_f$ intersects $X^*$, off $L$, in $\mu$ lines of the scroll, all passing through the same point of $L$.
        \end{itemize}
        In particular, $X^*$ is smooth variety if and only if one of the following equivalent cases hold:
        \begin{itemize}
            \item $d=3$, $X^*$ is a rational normal scroll, and $X$ contains the plane $L^*$ with multiplicity two;
            \item $X$ is the projection of the Segre variety $\sigma_{1,2}(\PP^1\times\PP^2)\subset\PP^5$ from an external point.
        \end{itemize}
    \end{theorem}
    \begin{proof}
        See \cite{franchetta} and \cite[Theorem 2.17]{cilibertorussosimis}. Optionally, a resume can be found in \cite[Theorem 7.4.14]{russo}.
    \end{proof}
    
    We can specialise the objects that we have just defined in Theorem \ref{cfrs}. In fact, $\Pi=L^*$, $C=Z^*$, and $W_f=Z$, where $Z=\overline{\nabla _f(\PP^n\setminus D)}$. By Theorem \ref{quinary}, up to change of variables, $f\in K[u,v][\Delta]$ with $\Delta=p_0x_0+p_1x_1+p_2x_2$ and $p_i\in K[u,v]_{s-1}$. In particular, $f=\Delta^{\mu}g_0+\Delta^{\mu-1}g_1+\dots+g_{\mu}$ with $\mu=\left[\frac{d}{s}\right]$ and $g_k\in K[u,v]_{ks}$. Thus, we obtain $\Pi=V(u,v)$ which has multiplicity $d-\mu$ in $X$.
    
    \pause
    Given a general hyperplane $H=\PP^3\subset\PP^4$, containing $\Pi$,  we have that $H\cap X$ is a hypersurface of degree $d$ in $H$. In particular, it is formed by the plane $\Pi$ with multiplicity $d-\mu$ and by $\mu$ planes of $\Sigma$ counted with multiplicity $1$ all intersecting the same point $p_H\in Z^*$. Looking at the proof of Theorem \ref{cfrs}, it is possible to see that the family $\Sigma$ is exactly the union of all the possible intersections $H\cap X$ varying the point $p_H$ in the rational curve $Z^*$.
    
    \begin{remark}\label{perazzopersempre}
       Since these objects depend only on the geometry of the hypersurface $X$, if $\Delta$ appears in $f$ not only linearly, then $f$ cannot be a Perazzo $3$-fold up to any change of variables.
    \end{remark}
    
    \section{The geometry of Perazzo 3-folds}
    In this section we are going to study the geometry of the Perazzo $3$-folds. In particular, we will specialise Theorem \ref{cfrs} to such forms, and then we will verify it for a very special class of Perazzo $3$-folds. In fact, we use the classification stated in Theorem \ref{lowcharacterisationreprise} to study the geometry of the Perazzo $3$-folds with minimal Hilbert vector and with $a=b=c=0$.
    
    \pause
    Let $f=p_0x_0+p_1x_1+p_2x_2+g$ be a Perazzo hypersuface as defined in Definition \ref{perazzo}, not a cone. The polar map is defined as
    \begin{equation*}
        \nabla_f([x_0,x_1,x_2,u,v])=\left[p_0(u,v):p_1(u,v):p_2(u,v):\frac{\partial f}{\partial u}:\frac{\partial f}{\partial v}\right].
    \end{equation*}
    
    Due to the shape of the form $f$, the Hessian matrix $\Hess(f)$ has rank four. In fact, we compute the Hessian matrix as
    \begin{equation*}
        \Hess(f)=\begin{pmatrix}
            \ast & \ast & \frac{\partial p_0}{\partial u} & \frac{\partial p_1}{\partial u} & \frac{\partial p_2}{\partial u}\\
            \ast & \ast & \frac{\partial p_0}{\partial v} & \frac{\partial p_1}{\partial v} & \frac{\partial p_2}{\partial v}\\
            \frac{\partial p_0}{\partial u} & \frac{\partial p_0}{\partial v} & 0 & 0 & 0\\
            \frac{\partial p_1}{\partial u} & \frac{\partial p_1}{\partial v} & 0 & 0 & 0\\
            \frac{\partial p_2}{\partial u} & \frac{\partial p_2}{\partial v} & 0 & 0 & 0
        \end{pmatrix}.
    \end{equation*}
    Consider now $B$, the $4\times 4$ sub-matrix formed by the first four columns and rows; then
    $$\det B=-\det\begin{pmatrix}
        \frac{\partial p_0}{\partial u} & \frac{\partial p_1}{\partial u}\\
        \frac{\partial p_0}{\partial v} & \frac{\partial p_1}{\partial v}
    \end{pmatrix}^2.$$
    By Theorem \ref{jacobian}, $\rank\begin{pmatrix}
        \frac{\partial p_0}{\partial u} & \frac{\partial p_1}{\partial u}\\
        \frac{\partial p_0}{\partial v} & \frac{\partial p_1}{\partial v}
    \end{pmatrix}=tr.deg_K K(p_0,p_1)=2$ because the two forms $p_0$ and $p_1$ are linear independent and we are working over an algebraically closed field.
    
    \pause
    Thus we have that $\dim X^*=2$ and $Z$ is a hypersurface. The equation of the variety $Z$ is exactly the irreducible algebraic relation of the forms $p_0,p_1,p_2$. We consider the plane $\Pi=V(u,v)$ which is contained in $X$ with multiplicity $d-1$. If in $(\PP^4)^*$ we set coordinates of the type $[z_0,z_1,z_2,z_3,z_4]$, then $\Pi$ is dual to a line $L=V(z_0,z_1,z_2)\subset X^*$. In particular, $Z$ is a rational cone with $L$ as a vertex. Let now $S=\PP^2\subset(\PP^4)^*$ be a plane containing the line $L$ and contained in $Z$. Then $X^*\cap S$ is the union of the line $L$ and another line contained in $X^*$ passing through a unique point of $L$. Thus, $X^*$ is a scroll surface formed by such lines.
    
    \pause
    On the other hand, we can give a precise description of the variety $X$. Since the rank of $\Hess(f)$ is four, then, for the general regular point $P$ in $X$, we have that $\langle P, (T_{\nabla_f(P)}Z)^*\rangle=\nabla^{-1}_f(\nabla_f(P))$. A general element of the family $\Sigma$ is the unique plane which contains $\nabla^{-1}_f(\nabla_f(P)$ and the tangent space of $Z^*$ at $(T_{\nabla_f(P)}Z)^*$ for some regular point $P$. In that way, given a linear sub-space $H=\PP^3\subset\PP^4$ containing $\Pi$, the intersection $X\cap H$ is the union of the plane $\Pi$ counted $d-1$ times and a unique plane of the family $\Sigma$.
    
    \begin{example}[Perazzo cubic]
        A special attention is given to the hypersurface $X=V(u^2x_0+uvx_1+v^2x_2)$ since it is, up to projectivity, the only Perazzo $3$-fold with degree $d=3$ (see Theorem \ref{cubicsinp4}); thus, it is the only possible choice if we ask to $X^*$ to be a smooth variety. As we have said $L$ has equations $z_0=0,z_1=0,z_2=0$; thus, the linear space $L^*$ has equations $u=0,v=0$, thus $L^*$ is precisely the plane $\Pi$. Clearly, this plane is contained in $X$ and, since the tangent cone at a point $Q=[\Tilde{x_0}:\Tilde{x_1}:\Tilde{x_2}:0:0]\in \Pi$ has equation $p_0(u,v)\Tilde{x_0}+p_1(u,v)\Tilde{x_1}+p_2(u,v)\Tilde{x_2}=0$, it has multiplicity two in $X$. We now see that $X^*$ is a rational normal scroll. Using \cite{mac2}, we have that $Z=V(z_1^2-z_0z_2)$, and the variety $X^*$ has equations
    \begin{equation*}
        \begin{cases}
            z_1z_3+z_2z_4=0\\
            z_0z_3+z_1z_4=0\\
            z_1^2-z_0z_2=0
        \end{cases}
    \end{equation*}
    Since at least one among the $2\times2$ minors of the Jacobian matrix
    \begin{equation*}
        J=\begin{pmatrix}
            0 & z_3 & -z_2\\
            z_3 & z_4 & 2z_1\\
            z_4 & 0 & -z_0\\
            z_1 & z_0 & 0\\
            z_2 & z_1 & 0
        \end{pmatrix}
    \end{equation*}
    does not vanish for every point in the variety $X^*$, we have that $X^*$ is smooth. In fact, $X^*$ is the union of the lines joining the points $[a^2:ab:b^2:0:0]$ and $[0:0:0:b:-a]$ varying the parameters $[a:b]\in\PP^1$. So that, $X^*$ is the scroll with directrix the line $L$ over the rational normal curve $C_2=\{z_1^2-z_0z_2=0\}\subset V(z_3,z_4)=\PP^2$. Eventually, $X$ is the projection of the Segre variety from a certain point. In fact, given the projection
    \begin{equation*}
        \begin{split}
            \pi:\PP^5&\longrightarrow\PP^4\\
            [w_0:w_1:w_2:w_3:w_4:w_5]&\longmapsto[w_1:w_3-w_0:-w_2:w_4:w_5]
        \end{split}
    \end{equation*}
    we have that $X^*$ is parameterized by $\PP^1\times\PP^2$ through the map $\pi\circ\sigma_{1,2}$.
    \end{example}
    
    We now proceed to study the following three kinds of Perazzo forms characterised by having minimal Hilbert vector:
    \begin{itemize}
        \item[(i)]   $f_1(x_0,x_1,x_2,u,v)= u^{d-1}x_0+u^{d-2}vx_1+u^{d-3}v^2x_2$,
        \item[(ii)]  $f_2(x_0,x_1,x_2,u,v)= u^{d-1}x_0+u^{d-2}vx_1+v^{d-1}x_2$,
        \item[(iii)]  $f_3(x_0,x_1,x_2,u,v)= u^{d-1}x_0+(\lambda u+\mu v)^{d-1}x_1+v^{d-1}x_2 \text{ with } \lambda, \mu \in K^*$.
    \end{itemize}
    
    In the case (i), $X_1=V(f_1)$ is not irreducible: in fact, it is the union of the hyperplane $u=0$ counted $d-3$ times, and the cubic Perazzo $3$-fold $u^2x_0+uvx_1+v^2x_2=0$. 
    
    \pause
    In the following two cases we use the following well-know results to study the rationality of a given curve non necessary smooth. The results contained in Theorem \ref{ultimoteorema} are also known as \emph{Pl\"ucker formulas} due to the mathematician who first proves them: J. Pl\"ucker.
    
    \begin{proposition}\label{ultimaproposizione}
        Let $C\in \PP^2$ be a curve of degree $d$. If $C$ has a singular point of multiplicity $d-1$, then $C$ is rational.
    \end{proposition}
    \begin{proof}[Proof (sketch).]
        Name $P$ the point of the curve with multiplicity $d-1$. Given a general line $L$ through the point $P$, we have, from intersection theory, that $L\cap C$ consists of the point $P$, with multiplicity $d-1$, and a different unique point $Q_L$. Parameterizing the lines through $P$ with $\PP^1$, we are able to construct a birational map $\phi:\PP^1\dashrightarrow C$.
    \end{proof}
    
    \begin{theorem}\label{ultimoteorema}
         Let $C\in \PP^2$ be a curve of degree $d$ and let $C^*$ be its dual: a curve of degree $d^*$. We suppose that $C$ has only $\delta$ ordinary double points, and $\kappa$ cusps as singular points. In the same way we suppose that $c^*$ has only $\delta^*$ ordinary double points, and $\kappa^*$ cusps. Then the following hold true:
         \begin{enumerate}
             \item $C$ is rational if and only if $g:=\frac{(d-1)(d-2)}{2}-\#\Sing(C)=0$;
             \item $d^* = d(d-1)-2\delta-3\kappa$;
             \item $\kappa^* = 3d(d-2)-6\delta-8\kappa$.
         \end{enumerate}
    \end{theorem}
    \begin{proof}
        See \cite[Chapter 2, Section 4]{griffithsharris}.
    \end{proof}
    
    We consider now case (ii), so that
    $$f_2(x_0,x_1,x_2,u,v)= u^{d-1}x_0+u^{d-2}vx_1+v^{d-1}x_2.$$
    We use coordinates $z_0,\ldots,z_4$ in $(\PP^4)^*$. A simple computations shows that the equation of $Z_2$ is $z_1^{d-1}-z_0^{d-2}z_2=0$, and the one of $Z_2^*$ is $(d-1)^{d-1}x_0^{d-2}x_2+(d-2)^{d-2}x_1^{d-1}=0$. The curves $Z_2^*$ and $Z_2\cap V(z_3,z_4)$ are both curves of degree $d-1$ with a singular point of multiplicity $d-2$; moreover, in that point, the tangent cone is only a line. By Proposition \ref{ultimaproposizione}, both are rational curves and, in particular, $Z_2$ is a rational cone of dimension three.
    
    \pause
    In case (iii) we have that
    $$ f_3(x_0,x_1,x_2,u,v)=u^{d-1}x_0+(\lambda u+\mu v)^{d-1}x_1+v^{d-1}x_2 \text{ with } \lambda, \mu \in K^*.$$
    
    Since the calculations are difficult to be performed in general, due to  the two free parameters $\lambda, \mu$, we use Macaulay2 \cite{mac2} and compute only the first cases. For low values of the degree $d$ of $f$, $Z_3$ is a cone over a rational curve of degree $d-1$ with $\frac{(d-2)(d-3)}{2}$ distinct nodes. Its dual $Z_3^*$ results to be a rational curve of degree $2d-4$. In particular, due to the rationality of $Z_3^*$ and by Theorem \ref{ultimoteorema}, we have the following:
    \begin{itemize}
        \item if $d=5$, then $Z_3^*$ has degree $6$ and it has $3$ cusps of multiplicity $3$ at the fundamental points $[1:0:0:0:0], [0:1:0:0:0], [0:0:1:0:0]$ and one node;
        \item if $d=6$, then $Z_3^*$ has three cusps of multiplicity $4$ at the same fundamental points, and $3$ nodes;
        \item if $d=7$, then $Z_3^*$ has three cusps of multiplicity $5$ at the same fundamental points, and $6$ nodes.
    \end{itemize}

    \cleardoublepage\phantomsection
    \addcontentsline{toc}{chapter}{References}
    \printbibliography[title={References}]
\end{document}